\documentclass[11pt]{amsart}
\usepackage{amsmath, amsfonts, amsthm, amssymb}
\usepackage[all,cmtip,line]{xy}\usepackage{mathtools}
\usepackage{hyperref} 
\usepackage{array}
\usepackage{enumerate}
\usepackage{srcltx} 

\usepackage{todonotes}
\newtheorem{thm}{Theorem}[section]
\newtheorem{lem}{Lemma}[section]
\newtheorem{prop}[lem]{Proposition}
\newtheorem{cor}[lem]{Corollary}
\newtheorem{q}[lem]{Question}
\newtheorem{defn}[lem]{Definition}
\newtheorem{rem}[lem]{Remark}

\numberwithin{equation}{section}

\newcommand{\bC}{ \mathbb{C}} 

\newcommand{\supp}{ \mbox{supp}}

\newcommand{\parallelsum}{\mathbin{\!/\mkern-5mu/\!}}

\newcommand \eps{\varepsilon}
\newcommand \vphi{\varphi}

\newcommand \om{\omega}

\newcommand{\R}{\mathbb R}

\newlength{\originalbase}
\title{ On the uniqueness of Poincar\'e type cscK metrics}

\author{Yulun Xu}  
\address{Mathematics Department and the Institute for Mathematical Sciences, Stony Brook University, Stony Brook NY, 11794-3651, USA}
\email{yulun.xu@stonybrook.edu}

\author{Kai Zheng}  
\address{University of Chinese Academy of Sciences, Beijing 100190, P.R. China}
\email{KaiZheng@amss.ac.cn}

\date{\today}

\begin{document}

\begin{abstract}
Let $X$ be a compact K\"ahler manifold with a divisor $D$ which may have simple normal crossings. We prove that the $K$-energy is convex on the space of Poincar\'e type K\"ahler metrics. Then, we prove that the Poincar\'e type constant scalar curvature K\"ahler metric is unique up to an automorphism of $X$ preserving $D$.
\end{abstract}
\maketitle

\tableofcontents

\section{Introduction}
\subsection{Background}
Let $(X,\omega_X)$ be a compact K\"ahler manifold. The space of K\"ahler potentials is
\begin{equation*}
    \mathcal{H}=\{\varphi \in C^{\infty}(X):\omega_{\varphi}=\omega_X+\sqrt{-1}\partial \bar{\partial}\varphi >0 \text{ on  }X\}.
\end{equation*}
At any $\varphi\in\mathcal{H}$, the tangent space $T_{\varphi}\mathcal{H}$ is $C^{\infty}(X)$. In a local holomorphic coordinate chart,
\[
\omega_\varphi = \sqrt{-1}\,g_{\varphi,\alpha\bar\beta}\,dz^\alpha \wedge d\bar z^\beta
=\sqrt{-1}\left(g_{\alpha \bar{\beta}}+\varphi_{\alpha\bar\beta}\right)dz^\alpha\wedge d\bar z^\beta.
\]
Its scalar curvature $R_\varphi$ is defined as:
\[
R_\varphi = -g_\varphi^{\alpha\bar \beta} \frac{\partial^2}{\partial z^\alpha \partial \bar{z}^\beta}\log \det(g_{\varphi,i\bar j}).
\]
The central problem in K\"ahler geometry which goes back to Calabi's program \cite{C2} \cite{C3} is to understand the existence and uniqueness problems of constant scalar curvature K\"ahler (cscK) metrics. As shown by Mabuchi \cite{M2}, the cscK metrics are the critical points of a functional called $K$-energy or the Mabuchi functional. By definition, a cscK metric satisfies the following equation:
\[
 R_\varphi=\bar R
 =n\frac{[\operatorname{Ric}(\omega_X)]\cdot[\omega_X]^{n-1}}{[\omega_X]^n}.
\]
The uniqueness of the cscK metric is  apparently a very difficult problem,  given that cscK equation is fourth order in the K\"ahler potential,  so the maximum principle is not applicable.  
The use of geodesics is really indispensable here,  as almost all the arguments for uniqueness of cscK metrics are based on convexity of $K$-energy along geodesics.

It is not hard to check that $K$-energy is convex along smooth geodesics.  However,   it is not always possible to connect two potentials in $\mathcal{H}$ through a smooth geodesic,  and the optimal regularity one can hope for is $C^{1,1}$ (cf.  \cite{CTW}, \cite{DLempert}, \cite{D5}, \cite{LV},\cite{CFH}).  Therefore,  it becomes a crucial question whether $K$-energy is convex along $C^{1,1}$ geodesics.  Moreover,  can one prove uniqueness of cscK based on this? In \cite{BB}, Berman and Berndtsson gave positive answers to these questions. 

Chen-Li-Paun provided in \cite{CLP} another way to prove the convexity of the $K$-energy. Chen-Paun-Zeng used PDE methods in \cite{CPZ} to prove the uniqueness of the cscK metric and extremal K\"ahler metric. 

\subsection{Poincar\'e type K\"ahler metric}
A natural question is: could we prove the uniqueness of cscK metric on more general K\"ahler manifolds? Note that the methods we mentioned above depend on the compactness of the manifold. As a result, it is interesting to consider the following Poincar\'e type K\"ahler metric which is non-compact and complete:
On $\mathbb{C}^n$, we can write down the standard local model for the Poincar\'e type K\"ahler metric: 
\begin{equation}\label{standard cusp}
    \omega_0= \sum_{j=1}^k \sqrt{-1}\frac{ 2 dz^j \wedge d\bar z^j}{|z_j|^2 \log^2(|z_j|^2)} +\sum_{j=k+1}^{n}\sqrt{-1}dz^j \wedge d \bar z^j.
\end{equation}
Writing $z_j$ in the form of polar coordinates, we have that:
\begin{equation*}
    z_j =r_j e^{i\theta_j}.
\end{equation*}
Then we can take $t_j = \log (-2 \log r_j)$. Then we can write $\omega_0$ as:
\begin{equation}\label{standard cusp polar}
    \omega_0 =  \sum_{j=1}^k - 2e^{-t_j} dt_j \wedge d\theta_j +\sum_{j=k+1}^{n}\sqrt{-1}dz^j \wedge d \bar z^j.
\end{equation}
We will use this expression many times in this paper.

We now define Poincar\'e type K\"ahler metrics. Let $D$ be a simple normal-crossing divisor on $X$, let $M=X\setminus D$, and write $D=\sum_{j=1}^N D_j$ as a sum of smooth irreducible components.
\begin{defn}\label{def 1.1}
    We say that $\omega$ is a Poincar\'e type K\"ahler metric of class $\Omega=[\omega_X]_{dR}$ if, for every $p\in D$ and every holomorphic coordinate neighbourhood $U$ of $p$ in which
    \(D\cap U=\{z_1\cdots z_k=0\}\) (called a cusp coordinate neighbourhood), $\omega$ satisfies:
    \begin{enumerate}
        \item There exists a constant $C$ such that $\frac{1}{C}\omega_0 \le \omega \le C \omega_0$ holds in $U$.
        \item There exists a function $\varphi$ defined on $M$ such that $\omega=\omega_X +dd^c \varphi$. For every $l\ge 1$, there exists a constant $C(l)$ such that $|\nabla^l_{\omega_0}\varphi|_{\omega_0}\le C(l)$ on $U$. Moreover, $$\varphi=O\big(\sum_{j=1}^k \log(-\log|z_j|)\big).$$
        \item $\omega$ is a smooth K\"ahler metric on $M$.
    \end{enumerate}
\end{defn}
\begin{rem}
    Note that in the definition of Poincar\'e type K\"ahler metric in \cite[Definition 1.4]{A}, an additional assumption is added: there exists a constant $C(k)$ such that
$
        |\nabla^k_{\omega_0} \omega|_{\omega_0}\le C(k)
$
    for any $k\ge 1$. In fact, this follows directly from part (2) of Definition \ref{def 1.1}, so we omit it.

    A global background metric is used in \cite[Definition 1.4]{A}. Here we use the local background metric $\omega_0$. In fact since the global background metric used there is $C^{\infty}$ quasi-isometric to $\omega_0$ in cusp coordinates, our definition is the same as the traditional definition.
\end{rem}
A lot of progress has been made for the Poincar\'e type cscK metrics. 
Auvray proved in \cite[Theorem 1]{A} the existence of Poincar\'e type geodesic. He also discovered a topological constraint for the Poincar\'e type cscK metrics in \cite{A2}, asymptotic properties of Poincar\'e type extremal K\"ahler metrics in \cite{A3} and the Poincar\'e type Futaki characters in \cite{A4}. See also the construction results obtained by Sektnan \cite{S} and Feng \cite{F} and the existence results on toric varieties by Apostolov, Auvray and Sektnan \cite{MR4215262}.

\subsection{Main results}

\begin{thm}\label{thm convex}
  Suppose that $D$ is a divisor, allowing simple normal crossings. The Poincar\'e type $K$-energy is convex on $\widetilde{\mathcal{PM}}_{\Omega}$.
\end{thm}
Here Poincar\'e type $K$-energy is the $K$-energy where the background metric is a Poincar\'e type metric instead of a smooth metric. See the Section \ref{subsection: Energy functionals}.
Then, we can prove that:

\begin{thm}\label{thm:main}
 Suppose that $D$ is a divisor, allowing simple normal crossings.  Let $\omega_1$ and $\omega_2$ be Poincar\'e type cscK metrics in the same fixed class $\Omega$. Then there exists $g\in\operatorname{Aut}_0^D(X)$ such that
\begin{equation*}
\omega_1=g^*\omega_2.
\end{equation*}
\end{thm}


Note that Auvray proved in \cite[Theorem 2]{A} the uniqueness of the Poincar\'e type cscK metric under the assumption that $K[D]$ is ample. We were informed that Aoi's result \cite{Ao} implies the uniqueness of Poincar\'e type cscK metric under the assumption that $(X,L_X)$ is polarized, $D\in |L_X|$ is smooth, $\mathbf{h}^D=0$ and $Aut_0^D(X)=\{Id\}$, using the uniqueness of conical cscK metric \cite{Z}.

\subsection{Key ingredients}

In the smooth case, Berman and Berndtsson's method consists of two main parts:
\begin{enumerate}\label{BB proof}
    \item First, the restriction of the $K$-energy to a $C^{1,1}$ geodesic can be viewed as an $S^1$-invariant function on a cylinder. They modify the entropy term and express its complex Hessian on the cylinder as the integral of a current. Positivity of this current gives subharmonicity of the modified $K$-energy. Continuity then yields convexity of the $K$-energy.
    \item Second, given two cscK metrics, they add a strictly convex functional to the $K$-energy and use it to show that the two cscK metrics are the same, up to a holomorphic transformation. A key  step in the proof is to solve the Lichnerowicz equation.  
\end{enumerate}

Compared with  the smooth  case, some new difficulties arise for the Poincar\'e type K\"ahler metrics. 
\begin{enumerate}
\item 
Bergman approximation is not available across Poincar\'e type singularities, although it is an important ingredient in Berman and Berndtsson's argument. Our key observation is that it is enough to carry out the approximation away from the singular set.

 When proving subharmonicity of the $K$-energy, we need to show that a current is positive on $X$. Since the Bergman approximation is local, it can be carried out on sets disjoint from $D$, which proves positivity on $M$. To obtain positivity on all of $X$ (see Lemma \ref{posi}), we show that the current has no mass on $D$. The required integrations by parts have no boundary contribution from $D$ (see Lemmas \ref{subharmonic} and \ref{G-S}). We also use the cut-off functions from Lemma \ref{snc-cutoff}, whose gradients and Hessians converge uniformly to zero.

Bergman approximation is also used to prove continuity of the $K$-energy. We take a locally finite cover of $M$ and use a subordinate partition of unity to decompose the $K$-energy into a sum of functionals supported away from $D$. Each functional is upper semicontinuous by the smooth argument. Uniform $C^{1,1}$ bounds control the error terms, so the $K$-energy is the decreasing limit of upper-semicontinuous functions and is therefore upper semicontinuous.

\item Another difficulty is solving the Lichnerowicz equation in the second step of Berman and Berndtsson's proof. On a closed manifold this is straightforward, but for a Poincar\'e type metric the operator acts on a complete noncompact space. The usual H\"older-space Fredholm alternative relies on compactness of the embedding $C^{4,\alpha}\hookrightarrow C^{0,\alpha}$ (see \cite[Theorem 5.3]{GT} and, in the conical setting, \cite[Theorem 5.7]{Z}), which fails here because one may translate a function towards the cusp while keeping its $C^{4,\alpha}$ norm fixed. We instead use a Fredholm alternative in weighted Sobolev spaces; see Proposition \ref{kernel approximation}. Our new observation is that in order to prove the uniqueness of Poincar\'e type cscK metric, it is enough to use such a solution in weighted Sobolev space instead of a solution in weighted H\"older space, which is difficult to get, especially in the case that $D$ has simple normal crossings and there is no additional assumption on the automorphism group of $D$ and $X$.
\end{enumerate}

In Section \ref{section: Preliminaries}, we list some machinery of Poincar\'e type K\"ahler metrics and calculate the gradient of the $K$-energy to show that Poincar\'e type cscK metrics are critical points of the $K$-energy. We also fix notation we will use in this paper. In Section \ref{section: Convexity of the K-energy}, we prove the convexity of the $K$-energy which proves Theorem \ref{thm convex} by showing that it is both weakly subharmonic and continuous. In Section \ref{section: Solvability of the Lichnerowicz operator}, we show how to solve the Lichnerowicz equation for the Poincar\'e type K\"ahler metrics. In Section \ref{section: Reductivity of holomorphic vector fields}, we prove the reductivity of holomorphic vector fields. In Section \ref{section :Uniqueness of the Poincare type cscK metrics}, we prove Theorem \ref{thm:main}.

\section{Acknowledgements}
The first author is partially funded by the Simons Foundation. He thanks his advisors, Jingrui Cheng and Xiuxiong Chen, for their encouragement and discussions regarding this work. 
The second author is partially funded by NSFC grants No. 12171365, 12571091, and 12326426.
He would like to express his deepest gratitude to IHES and the K.C. Wong Education Foundation, also the CRM and the Simons Foundation, when parts of the work were undertaken during his stays.
We also thank Long Li for the discussions about this work.

\section{Preliminaries}\label{section: Preliminaries}

\subsection{Conventions}
 We use the following conventions throughout the paper:
\begin{equation}\label{operator conventions}
\begin{gathered}
d^c\coloneqq \frac{\sqrt{-1}}{2}(\bar\partial-\partial),\qquad
dd^c=\sqrt{-1}\,\partial\bar\partial,\qquad
\omega=\sqrt{-1}\,g_{\alpha\bar\beta}\,dz^\alpha\wedge d\bar z^\beta,\\
\Delta_\omega f\coloneqq-\Lambda_\omega dd^c f
=-g^{\alpha\bar\beta}f_{\alpha\bar\beta},\qquad
\operatorname{Ric}(\omega)\coloneqq-dd^c\log\det(g_{\alpha\bar\beta}),\\
R_\omega\coloneqq\Lambda_\omega\operatorname{Ric}(\omega)
=-g^{\alpha\bar\beta}\partial_\alpha\partial_{\bar\beta}
\log\det(g_{\gamma\bar\delta}).
\end{gathered}
\end{equation}
Thus $\Delta_\omega$ is the nonnegative complex Laplacian (one half of the nonnegative Riemannian Laplacian) and
\begin{equation}\label{basic variation conventions}
n\alpha\wedge\omega^{n-1}=(\Lambda_\omega\alpha)\omega^n,
\qquad
\left.\frac{d}{dt}\right|_{t=0}\omega_{\varphi+tv}^n
=-(\Delta_{\omega_\varphi}v)\omega_\varphi^n.
\end{equation}
We write $\operatorname{tr}_\omega\alpha=\Lambda_\omega\alpha$. All integrals below use the unnormalised volume form $\omega^n$. 

\subsection{Background metric of Poincar\'e type}
First, we construct a Poincar\'e type K\"ahler metric to use as a background metric. Write $D=\sum_{i=1}^k D_i$, where the $D_i$ are smooth divisors. For each $i$, take a holomorphic defining section $\sigma_i$ of $\mathcal{O}([D_i])$ and fix a smooth Hermitian metric $|\cdot|_i$ on this line bundle. Set
$$\rho_i \triangleq -\log(|\sigma_i|_i^2)\ge 1$$ 
on $X\setminus D_i$, after rescaling the Hermitian metric so that $|\sigma_i|_i^2\leq e^{-1}$. Let $\lambda$ be a nonnegative constant to be determined, and set
$$\mathbf{u}_i \triangleq \log(\lambda+\rho_i),\,\,\, \mathbf{u}\coloneqq 1+\sum_{j=1}^{k}\mathbf{u}_j.
$$ 
Auvray shows in \cite[Lemma 1.1]{A} that, for any $A>0$ and all sufficiently large $\lambda$ depending on $A$ and $\omega_X$, the $(1,1)$-form $\omega_X-A\sum_{i=1}^k dd^c\mathbf{u}_i$ is a Poincar\'e type K\"ahler metric.
From now on, we denote 
$$\omega\triangleq\omega_X-A\sum_{i=1}^k dd^c\mathbf{u}_i,\qquad
\rho\triangleq\prod_{i=1}^k(\lambda+\rho_i).$$
We define the space of Poincar\'e type K\"ahler metrics of class $\Omega$ as $ \mathcal{PM}_{\Omega}$. Denote the space of potentials (with respect to the background metric $\omega$) as
$\widetilde{\mathcal{PM}}_{\Omega}$.

\subsection{Quasi coordinates}
Next, the quasi-coordinates (see \cite{TY}) are used to define function spaces associated with Poincar\'e-type K\"ahler metrics.
Let $\Delta$ be a unit disc and let $\Delta^*$ be a punctured unit disc. For any $\delta\in (0,1)$, we can set $$\varphi_{\delta}: \frac{3}{4}\Delta \rightarrow \Delta^*,\quad\xi \mapsto \exp(-\frac{1+\delta}{1-\delta}\frac{1+\xi}{1-\xi}).$$ For any $\delta =(\delta_1,...,\delta_k) \in (0,1)^k$, we can define \begin{align*}
\Phi_{\delta}: \mathcal{P} \triangleq (\frac{3}{4}\Delta)^k \times \Delta^{n-k}  &\rightarrow  (\Delta^*)^k \times \Delta^{n-k} ,\quad \delta \in (0,1)^k,\\
(\xi_{1},...,\xi_{k}, z_{k+1},...,z_n)&\mapsto (\varphi_{\delta_1}(\xi_1),...,\varphi_{\delta_k}(\xi_{k}),z_{k+1},...,z_n).
 \end{align*} 
 We call a holomorphic coordinate system a cusp coordinate system if $D=\{z_1\cdots z_k=0\}$ in that coordinate chart. Let $\omega_0$ be the standard local Poincar\'e type metric given by (\ref{standard cusp}). A direct computation gives
 \begin{equation*}
     \Phi_{\delta}^* \omega_0= \sum_{i=1}^k \frac{2\sqrt{-1}d\xi^i \wedge d \bar \xi^i }{(1-|\xi_i|^2)^2} + \sum_{i=k+1}^n \sqrt{-1}dz^i \wedge d\bar z^i
 \end{equation*}
 which is uniformly quasi-isometric to the Euclidean metric independent of $\delta$. Since any Poincar\'e type metric $\omega$ is quasi-isometric to $\omega_0$ in each cusp coordinate, we have that $\Phi^*_{\delta}\omega$ is quasi-isometric to Euclidean metric in each cusp coordinate.

\subsection{Function spaces}
Auvray \cite{A} used the idea of Tian-Yau \cite{TY} to define the following function spaces:
\begin{defn}
If $U$ is a polydisc neighbourhood of $D$ with $U \cap D=\{z_1\cdots z_k=0\}$, we define, for $f\in C^{p,\alpha}_{loc}(U \setminus D)$ and $(p,\alpha)\in \mathbb{N}\times [0,1)$,
\begin{equation*}
||f||_{C^{p,\alpha}(U \setminus D)} \triangleq \sup_{\delta \in (0,1)^k}||\Phi_{\delta}^* f||_{C^{p,\alpha}(\mathcal{P})},
\end{equation*}
assuming that $U \subset   (c\Delta)^k \times \Delta^{n-k}.$ 

Then given a finite number of such open sets $U\in \mathcal{U}$, covering $D$ and an open set $V \subset \subset X \setminus D$ such that $X=V \cup \bigcup_{U \in \mathcal{U}} U$ and a partition of unity $\{\chi_V\} \cup \{\chi_U: U \in \mathcal{U}\},$ we can define the H\"older space 
\begin{equation*}
C^{p,\alpha}(M) \triangleq \{f\in C_{loc}^{p,\alpha}(M) : ||\chi_V f||_{C^{p,\alpha}(V)}+\max_{U\in 
 \mathcal{U}} ||\chi_U f||_{C^{p,\alpha}(U \setminus D)} < \infty \}.
\end{equation*}
\end{defn}
\begin{defn}
We can define the weighted H\"older norm:
\begin{equation*}
C_{\eta}^{l,\alpha} \triangleq \{f\in C^{l,\alpha}_{loc}(M) : ||\chi_V f||_{C^{l,\alpha}(V)}+\sup_{U\in \mathcal{U}} \sup_{\delta \in (0,1)^k}||\prod_{i=1}^k(1-\delta_i)^{\eta}\Phi_{\delta}^* (\chi_U f)||_{C^{l,\alpha}(\mathcal{P})} < \infty\}.
\end{equation*}
Since  $\prod_{i=1}^k\frac{1}{C(1-\delta_i)} \le \Phi_{\delta}^* \rho \le \prod_{i=1}^k\frac{C}{1-\delta_i}$ for some constant $C$, $||\prod_{i=1}^k(1-\delta_i)^{\eta}\Phi_{\delta}^* (\chi_U f)||_{C^{l,\alpha}(\mathcal{P})} $ is equivalent to $||\Phi_{\delta}^* (\rho^{-\eta} \chi_U f)||_{C^{l,\alpha}(\mathcal{P})} $.
Heuristically, $f\in C_{\eta}^{l,\alpha}$  implies that $f =O(\rho^{\eta})$. We can also define:
\begin{equation*}
    C_{\eta}^{\infty}= \cap_{l=0}^{\infty} C_{\eta}^{l,\alpha}.
\end{equation*}
\end{defn}
Then we can define 
\begin{equation*}
    C^{\infty}(M)= \cap_{k=1}^{\infty}C^{k,\alpha}(M).
\end{equation*}

\begin{defn}
 We can also define the weighted Sobolev space:
\begin{equation*}
W_{\eta}^{k,p} \triangleq \{v \in W^{k,p}_{loc}(M): \int_{M} \sum_{i=0}^k |\nabla_i v|^p \rho^{-\eta p}\omega^n < \infty \}.
\end{equation*}
\end{defn}
Clearly, $W_{\eta}^{k,p} \subset W_{\eta'}^{k,p} $, when $\eta\leq \eta'$.

\subsection{Fiber bundle structure of a neighbourhood of D}\label{subsection: Fiber bundle structure of a neighbourhood of D}
Consider the special case that $D$ is a smooth divisor. According to \cite[Section 3]{A2}, a neighbourhood of $D$, denoted as $\mathcal{N}_A$, can be seen as an $S^1$ bundle over $[A,\infty) \times D$. This fiber bundle can be written as $$q: \mathcal{N}_A \setminus D \xrightarrow{q=(t,p)} [A,\infty) \times D.$$   The function $t$ is defined in \cite{A2}. We have that $t=\mathbf{u}$ up to a perturbation which is a $O(e^{-t})$ term, that is, a $O(\frac{1}{|\log |\sigma||})$ term, as well as its derivatives of any order with respect to Poincar\'e metrics. Denote $p$ as the projection from $\mathcal{N}_A \setminus D$ to $D$. We can also define a connection $\widetilde{\eta}$ in $\mathcal{N}_A \setminus D$ which can be seen as a volume form on each $S^1$ fibre such that $$J dt =2e^{-t} \widetilde{\eta} +O(e^{-t}).$$
In a cusp coordinate $(z_1,...,z_n=r e^{i\theta})$, one has 
\begin{equation}\label{eta dtheta}
\widetilde{\eta} =d \theta +O(1)
\end{equation}
in the sense that $\widetilde{\eta}- d\theta$ and all the derivatives of it of any order with respect to $\omega$ is bounded. 
Given an arbitrary function $f$ supported in a neighbourhood $\mathcal{N}_A$ of $D$, we can decompose $f$ as:
\begin{equation}\label{f decom}
f= f_0(t,p) + f^{\bot},
\end{equation}  
where $$f_0(t,p)=\frac{1}{2\pi}\int_{q^{-1}(t,p)} f \widetilde{\eta}$$ is the $S^1$ invariant part and $f^{\bot}$ is the part that is orthogonal to $S^1$-invariant functions.

\subsection{Cut-off functions}\label{subsection: Cut-off functions}
\begin{lem}
\label{snc-cutoff}
Let \(D=\sum_{j=1}^{N}D_j\) be a simple normal-crossing divisor in a compact K\"ahler manifold \(X\), and let \(\omega\) be a Poincar\'e type K\"ahler metric on \(M=X\setminus D\). There exist a smooth proper function
\[
T:M\longrightarrow [0,+\infty)
\]
and functions \(\eta_m\in C^\infty(M,[0,1])\), \(m\ge1\), with the following properties:
\begin{enumerate}
\item \(T(x)\to+\infty\) as \(x\to D\);
\item if \(\mathcal N_A=\{T>A\}\), then
\[
\eta_m=0\quad\text{on }M\setminus\mathcal N_m,
\qquad
\eta_m=1\quad\text{on }\mathcal N_{2m};
\]
\item \(1-\eta_m\in C_c^\infty(M)\);
\item
$\sup_M\left(
|d\eta_m|_\omega+
|\nabla_\omega^2\eta_m|_\omega
\right)\le \frac{C}{m}$;
\item \(\eta_m\searrow0\) pointwise on \(M\);
\item 
\begin{equation}\label{eq:constructive-cutoff-four}
\lVert\nabla^r\eta_m\rVert_{L^\infty(\omega)}
\le \frac{C_r}{m},\quad \forall r\geq1.
\end{equation}
\end{enumerate}
The same construction applies on \(M\times R\) after pulling \(T\) and \(\eta_m\) back from \(M\).
\end{lem}

\begin{proof}
For every irreducible component \(D_j\), choose a defining section
\(\sigma_j\in H^0(X,\mathcal O(D_j))\) and a smooth Hermitian metric
\(h_j\) on \(\mathcal O(D_j)\). After rescaling \(h_j\), we may assume
\[
\rho_j\coloneqq-\log|\sigma_j|_{h_j}^2\ge1
\]
on \(M\). Put \(t_j=\log\rho_j\) and define
\[
T=\log\left(\sum_{j=1}^{N}\rho_j\right)
 =\log\left(\sum_{j=1}^{N}e^{t_j}\right).
\]

Let \(p\in D\), and suppose that \(D_1,\ldots,D_k\) are the components
passing through \(p\). In holomorphic coordinates centered at \(p\), we
may write
\[
D=\{z_1\cdots z_k=0\}.
\]
In these coordinates,
\[
t_j=\log(-\log|z_j|^2)+O(1)
\]
for \(j\le k\), with the error and all its derivatives controlled with
respect to a product Poincar\'e metric. It follows from the local model that
\[
|dt_j|_\omega+|\nabla_\omega^2t_j|_\omega\le C.
\]
For \(j>k\), the function \(t_j\) is smooth near \(p\), and the same
estimate is immediate.

Set
\[
q_j=\frac{e^{t_j}}{\sum_{\ell=1}^{N}e^{t_\ell}}.
\]
Then
\[
dT=\sum_jq_jdt_j
\]
and
\[
\nabla^2T
 =
 \sum_jq_j\nabla^2t_j
 +\sum_jq_jdt_j\otimes dt_j
 -\left(\sum_jq_jdt_j\right)
  \otimes
  \left(\sum_\ell q_\ell dt_\ell\right).
\]
Since \(0\le q_j\le1\) and \(\sum_jq_j=1\), we obtain
\[
|dT|_\omega+|\nabla_\omega^2T|_\omega\le C.
\]
Moreover, \(T\to+\infty\) whenever at least one \(\sigma_j\) tends to
zero. Since \(X\) is compact, every sublevel set \(\{T\le A\}\) is
therefore compact in \(M\).

Choose a nondecreasing function
\(\chi\in C^\infty(\mathbb R,[0,1])\) such that
\[
\chi=0\quad\text{on }(-\infty,0],
\qquad
\chi=1\quad\text{on }[1,+\infty),
\]
and define
\[
\eta_m=\chi\left(\frac{T}{m}-1\right).
\]
The support properties follow immediately. Furthermore,
\[
d\eta_m
 =\frac1m\chi'\left(\frac Tm-1\right)dT
\]
and
\[
\nabla^2\eta_m
 =
 \frac1{m^2}\chi''\left(\frac Tm-1\right)dT\otimes dT
 +\frac1m\chi'\left(\frac Tm-1\right)\nabla^2T.
\]
The uniform bounds for \(dT\) and \(\nabla^2T\) imply
\[
|d\eta_m|_\omega+|\nabla_\omega^2\eta_m|_\omega
\le\frac{C}{m}.
\]
Finally, for fixed \(x\in M\), the quantity \(T(x)/m-1\) decreases as
\(m\) increases. Since \(\chi\) is nondecreasing, \(\eta_m(x)\) decreases
to zero.

Then we estimate the derivatives of the cutoff function in Lemma~\ref{snc-cutoff} up to fourth order $1\le r\le4$
\begin{equation}\label{eq:constructive-cutoff-support}
\begin{gathered}
\operatorname{supp}\nabla^r\eta_m
\subset\{m\le T\le2m\}.
\end{gathered}
\end{equation}

The proof is a straightforward computation.
We  differentiate $r$ times $\chi(T/m-1)$, in which each term has the form
\[
m^{-k}\chi^{(k)}(T/m-1)
\nabla^{r_1}T\otimes\cdots\otimes\nabla^{r_k}T,
\qquad r_i\ge1,\quad r_1+\cdots+r_k=r,
\]
with universal coefficients. 
Since $T$ has bounded derivatives with respect to the Poincar\'e metric,
\begin{equation*}
\lVert\nabla^rT\rVert_{L^\infty(\omega)}\le C_r.
\end{equation*}
Therefore, we prove that for $m\ge1$,
\begin{equation}\label{eq:constructive-cutoff-four}
\lVert\nabla^r\eta_m\rVert_{L^\infty(\omega)}
\le \frac{C_r}{m},\qquad 1\le r\le4.
\end{equation}
Furthermore, we could prove it for any $r$ in the same way.
\end{proof}

\subsection{The geodesic equation}

Let \(R=[0,1]\times S^1\) be the cylinder. We use the local holomorphic
coordinate
\[
w=t+\sqrt{-1}\,s
\]
on \(R\); although \(w\) is only locally defined, \(dw\) is globally defined.
Set
\[
\mathfrak X=X\times R,\qquad
\mathfrak D=D\times R,\qquad
\mathfrak M=M\times R,
\]
and let \(\pi\colon\mathfrak X\to X\) be the projection. On
\(\mathfrak M\), consider the product Poincar\'e type K\"ahler form
\begin{equation}\label{conventions}
\omega^\ast
\coloneqq
\pi^\ast\omega+\frac{\sqrt{-1}}{2}\,dw\wedge d\bar w.
\end{equation}

Let \(\{\varphi_t\}_{t\in[0,1]}\) be a path in
\(\widetilde{\mathcal{PM}}_{\Omega}\), and let
\(\Phi(x,w)=\varphi_t(x)\) denote its \(S^1\)-invariant extension to
\(\mathfrak M\). Define
\begin{equation}\label{Psi definition}
\Psi\coloneqq\Phi+t(1-t).
\end{equation}
Since
\[
dd^c\bigl(t(1-t)\bigr)
=-\frac{\sqrt{-1}}{2}\,dw\wedge d\bar w,
\]
we have
\[
\omega^\ast+dd^c\Psi
=\pi^\ast\omega+dd^c\Phi.
\]
Writing
\(\omega_{\varphi_t}=\omega+dd^c\varphi_t\), a local determinant
calculation gives
\begin{equation}\label{geo ma}
\begin{split}
\bigl(\omega^\ast+dd^c\Psi\bigr)^{n+1}
&=\frac{n+1}{4}
\left(
\ddot{\varphi}_t
-\lvert\partial_X\dot{\varphi}_t\rvert_{\omega_{\varphi_t}}^2
\right)
\pi^\ast\omega_{\varphi_t}^{\,n}
\wedge\sqrt{-1}\,dw\wedge d\bar w.
\end{split}
\end{equation}
Here a dot denotes \(\partial_t\), and \(\partial_X\) denotes the
Dolbeault operator in the \(X\)-directions.

A cusp geodesic joining \(\varphi_0\) and \(\varphi_1\) is a path for which
\[\omega^\ast+dd^c\Psi\geq 0,\quad \bigl(\omega^\ast+dd^c\Psi\bigr)^{n+1}=0\text{ on }\mathfrak M,\] with
Dirichlet data
\begin{equation}\label{geo ma boundary}
\left.\Psi\right|_{t=0}=\varphi_0,\qquad
\left.\Psi\right|_{t=1}=\varphi_1.
\end{equation}

For the elliptic approximation, put
\[
\widetilde{\varphi}_0(x,t)
\coloneqq(1-t)\varphi_0(x)+t\varphi_1(x).
\]
Choose \(c\in C^\infty([0,1])\) with \(c(0)=c(1)=0\) and \(c''\)
sufficiently large that
\[
\Psi_1\coloneqq\widetilde{\varphi}_0+c(t),\qquad
\omega_1\coloneqq\omega^\ast+dd^c\Psi_1
\]
defines a Poincar\'e type K\"ahler form on \(\mathfrak M\). Such a choice
is possible because the fiber metrics
\((1-t)\omega_{\varphi_0}+t\omega_{\varphi_1}\) are positive and the
mixed terms are uniformly bounded with respect to them.

\begin{rem}\label{choose c}
One may take \(c(t)=A t(t-1)\), where \(A>0\) is sufficiently large and
depends only on the endpoint potentials and the background geometry.
\end{rem}

We record the identities
\begin{lem}
\begin{align*}
\Psi-\Psi_1
&=\Phi-\widetilde{\varphi}_0+t(1-t)-c(t),\\
\partial_t\widetilde{\varphi}_0
&=\varphi_1-\varphi_0,
&
\partial_X\widetilde{\varphi}_0
&=(1-t)\partial_X\varphi_0+t\partial_X\varphi_1,\\
\partial_t(\Psi-\Psi_1)
&=\dot\varphi_t-(\varphi_1-\varphi_0)+1-2t-c'(t),\\
\partial_X(\Psi-\Psi_1)
&=\partial_X(\varphi_t-\widetilde{\varphi}_0).
\end{align*}
\end{lem}

\subsection{Approximate cusp geodesic}
For \(0<\epsilon\leq1\), consider the Dirichlet problem
\begin{equation}\label{approximate cusp equation}
\begin{cases}
\bigl(\omega^\ast+dd^c\Psi^\epsilon\bigr)^{n+1}
=\epsilon\,\omega_1^{n+1}
& \text{on }\mathfrak M,\\[2mm]
\left.\Psi^\epsilon\right|_{t=0}=\varphi_0,
\qquad
\left.\Psi^\epsilon\right|_{t=1}=\varphi_1.
\end{cases}
\end{equation}
Write \(\varphi_t^\epsilon=\Psi^\epsilon-t(1-t)\).
For \(\epsilon=1\), the function \(\Psi^1=\Psi_1\) solves
\eqref{approximate cusp equation}. As \(\epsilon\to0\), a limit of solutions
to \eqref{approximate cusp equation} solves the geodesic equation
\eqref{geo ma} in the weak sense.

For a function \(u\) on \(\mathfrak M\), set
\begin{equation*}
\lVert u\rVert_{C_{\triangle}}
\coloneqq
\sup_{(x,t)\in\mathfrak M}
\left(
\lvert u\rvert
+\lvert\partial_tu\rvert
+\lvert\partial_Xu\rvert_{\omega}
+\lvert\partial_X\bar\partial_Xu\rvert_{\omega}
\right)
\end{equation*}
and
\begin{equation*}
\lVert u\rVert_{C_{\mathbf w}^{1,1}}
\coloneqq
\lVert u\rVert_{C_{\triangle}}
+\sup_{(x,t)\in\mathfrak M}
\left(
\lvert\partial_X\partial_tu\rvert_{\omega}
+\lvert\partial_t^2u\rvert
\right).
\end{equation*}

\begin{defn}\label{ctribw}
A path \(\{\varphi_t\}_{t\in[0,1]}\) is a generalized cusp geodesic if
it is the limit, as \(\epsilon\to0\), of solutions to
\eqref{approximate cusp equation} with respect to the norm $C^{1,\alpha}_{loc}$ and
\(\sup_\epsilon\lVert\varphi^\epsilon-\widetilde{\varphi}_0\rVert_{C_{\triangle}}<\infty\).
\end{defn}

\begin{defn}\label{c11bw}
A path \(\{\varphi_t\}_{t\in[0,1]}\) is a cusp geodesic if it is the
limit, as \(\epsilon\to0\), of solutions to
\eqref{approximate cusp equation} with respect to the norm $C^{1,\alpha}_{loc}$ and
\(\sup_\epsilon\lVert\varphi^\epsilon-\widetilde{\varphi}_0\rVert_{C_{\mathbf w}^{1,1}}<\infty\).
\end{defn}

\begin{thm}\label{weakgeodesicclosedness}
Assume that \(\varphi_0\) and \(\varphi_1\) are cusp potentials whose
associated K\"ahler cusp metrics have bounded \(C^4\)-norm with respect
to the background cusp metric. Then there exists a unique
\(C_{\mathbf w}^{1,1}\) cusp geodesic joining \(\varphi_0\) and
\(\varphi_1\). Moreover, the approximating solutions satisfy
\begin{equation*}
\lVert\varphi^\tau-\widetilde{\varphi}_0\rVert_{C_{\mathbf w}^{1,1}}
\leq C,
\end{equation*}
where \(C\) is independent of \(\tau\).
\end{thm}

\begin{rem}
The constant depends on a lower bound for the bisectional curvature of
\(\omega_1\) and on bounds for the first derivatives of \(\omega_1\).
\end{rem}

\begin{q}
Can two \(C^{1,1}\) cusp potentials be joined by a generalized cusp
geodesic?
\end{q}

 \subsection{Estimates}\label{l infty} 
 We omit the upper index $\epsilon$ in this subsection.
\begin{lem}\label{low}(Lower bound of $\Psi$)
The following estimate holds
$$
\Psi-\Psi_1\geq 0\text{ in }\mathfrak X\; .
$$
\end{lem}
\begin{proof}
Taking logarithms of \eqref{approximate cusp equation}, we have the equation for $\tilde\Psi=\Psi-\Psi_1$:
\begin{align}\label{log equ}
\log\frac{\det(g_{1,i\bar j}+\tilde\Psi_{i\bar j})}
{\det(g_{1,i\bar j})}=\log\epsilon\leq 0,
\end{align}
where \((g_{1,i\bar j})\) is the local coefficient matrix of
\(\omega_1\). We denote
\(g_{a,i\bar j}=g_{1,i\bar j}+a\tilde\Psi_{i\bar j}\). Then
\begin{align}
\left(\int_0^1g_a^{i\bar j}\,da\right)
\tilde\Psi_{i\bar j}\leq0\quad\text{on }\mathfrak M,
\qquad
\tilde\Psi=0\quad\text{on }\partial\mathfrak X.
\end{align}
The coefficients have Poincar\'e singularities. Thus, the conclusion
follows from the maximum principle for the boundary value problem
\cite[Lemma 2.6]{A}.
\end{proof}


\begin{lem}\label{upp}(Upper bound of $\Psi$)
The following estimate holds
\begin{align*}
\Psi-\Psi_1\leq h \text{ in }\mathfrak X .
\end{align*}
where
$h$ satisfies the linear equation 
\begin{equation*}
  \left\{
   \begin{array}{ll}
\Delta_1 h =n+1 & \text{ in }\mathfrak{M}\; ,\\
h=0 & \text{ on }\partial\mathfrak{X} \;,
   \end{array}
  \right.
\end{equation*}
and \(\Delta_1\) is the Laplacian with respect to \(\omega_1\).
\end{lem}
\begin{proof}
Taking the trace of
\(\omega_1+\sqrt{-1}\partial\bar\partial\tilde\Psi\geq0\)
with respect to \(\omega_1\),
we have
$$
\Delta_1 \tilde\Psi\leq n+1 \; .
$$
The existence of $h$ is provided by \cite[Lemma 2.7]{A}. Then we
apply the maximum principle for the boundary value problem \cite[Lemma 2.6]{A} to obtain the required bound.
\end{proof}

\begin{lem}[Boundary gradient estimate]\label{boundary gradient estimate}
\begin{equation}
	\begin{split}
	\partial_X\Psi= 
	\partial_X\Psi_{1} &,\quad 
	0\leq \partial_t(\Psi-\Psi_1)\leq \partial_th \text{ for } t=0, \\
	\partial_X\Psi= 
	\partial_X\Psi_{1} &,\quad 
	0\geq \partial_t(\Psi-\Psi_1)\geq \partial_th \text{ for } t=1.
	\end{split}
\end{equation}
\end{lem}
\begin{proof}
Along the directions tangent to the boundary, the gradient estimates are obtained from the boundary value $\Psi_1$. Along the $t$ direction, we have
\begin{align*}
0\leq \partial_t(\Psi-\Psi_1)\leq \partial_th
\quad\text{on }t=0.
\end{align*}
The analogous reversed inequalities hold on \(t=1\). Thus,
\(\lvert\partial_t\Psi\rvert\) on the boundary is bounded in terms of
\(\sup_{\partial\mathfrak X}\lvert\partial_t\Psi_1\rvert\) and
\(\sup_{\partial\mathfrak X}\lvert\partial_t h\rvert\).
\end{proof}

\begin{lem}[$t$ gradient estimate]\label{vdot estimate}
	\begin{align*}
	-C\leq\partial_t(\Psi-\Psi_1)\leq C
	\quad\text{on }\mathfrak X.
	\end{align*}
	Here \(C\) is a fixed positive constant depending on the choice of
	the convex function \(c(t)\) in the construction of \(\Psi_1\).
\end{lem}
\begin{proof}
We compute 
\begin{align*}
\partial_t\Psi=\dot\varphi_t+1-2t,
\qquad
\partial_t\Psi_1=\varphi_1-\varphi_0+c'(t).
\end{align*} 
Set
\(u=\Psi+t(t-1)-\Psi_1+c(t)=\varphi_t-\widetilde{\varphi}_0\).
From \eqref{geo ma} and \eqref{approximate cusp equation}, we have
\begin{align*}
&\partial_t u
=\partial_t(\Psi-\Psi_1)+2t-1+c'(t)
=\dot\varphi_t-(\varphi_1-\varphi_0),\\
&\partial_t^2u=\ddot\varphi_t\geq0,
\qquad t\in[0,1].
\end{align*}
Thus, \(\partial_tu\) is controlled by its boundary values:
\begin{align}\label{pt u bound}
\left.\partial_tu\right|_{t=0}
\leq\partial_tu
\leq\left.\partial_tu\right|_{t=1}.
\end{align}
Using Lemma~\ref{boundary gradient estimate}, we have
\begin{align*}
&\left.\partial_tu\right|_{t=0}
=\left.\partial_t(\Psi-\Psi_1)\right|_{t=0}+c'(0)-1
\geq c'(0)-1,\\
&\left.\partial_tu\right|_{t=1}
=\left.\partial_t(\Psi-\Psi_1)\right|_{t=1}+c'(1)+1
\leq c'(1)+1.
\end{align*}
Returning to \eqref{pt u bound}, we get
\begin{equation*}
\begin{split}
c'(0)
&\leq\left.\partial_tu\right|_{t=0}+1
\leq\partial_tu+1\\
&=\partial_t(\Psi-\Psi_1)+2t+c'(t)
\leq\left.\partial_tu\right|_{t=1}+1
\leq c'(1)+2.
\end{split}
\end{equation*}
That gives
\begin{align*}
&\partial_t(\Psi-\Psi_1)
\geq c'(0)-2t-c'(t)
\geq c'(0)-2-c'(t),\\
&\partial_t(\Psi-\Psi_1)
\leq2-2t+c'(1)-c'(t)
\leq2+c'(1)-c'(t).
\end{align*}
The bound for \(C\) follows from the choice of \(c\); see
Remark~\ref{choose c}.
 \end{proof}
 
The following result is the form of
\cite[Theorem~2.1 and Corollary~2.2]{A} used below.

\begin{lem}\label{eps geodesic}
Let
\(\varphi_0,\varphi_1\in\widetilde{\mathcal{PM}}_{\Omega}\).
For every sufficiently small \(\epsilon>0\), there exists an
\(\epsilon\)-geodesic joining \(\varphi_0\) to \(\varphi_1\).
Moreover, there is a constant \(C>0\), independent of \(\epsilon\), such
that on \(M\times[0,1]\),
\begin{equation}\label{epsilon geodesic estimates}
\begin{split}
&\lvert\varphi^\epsilon-\widetilde{\varphi}_0\rvert
+\lvert d_X\varphi^\epsilon\rvert_\omega
+\lvert\ddot{\varphi}^\epsilon\rvert 
+\lvert d_X\dot{\varphi}^\epsilon\rvert_\omega
+\lvert d_Xd_X^c\varphi^\epsilon\rvert_\omega
\le C.
\end{split}
\end{equation}
The difference
\(\varphi^\epsilon-\widetilde{\varphi}_0\) is smooth on
\(M\times[0,1]\).
\end{lem}

Letting \(\epsilon\to0\) in the monotone family constructed in
\cite{A} gives the weak geodesic and the estimates needed later.

\begin{lem}\label{C11 geodesic}
For any
\(\varphi_0,\varphi_1\in\widetilde{\mathcal{PM}}_{\Omega}\), there is a
continuous geodesic \(\{\varphi_t\}_{t\in[0,1]}\) joining them such that
\[
\bigl(\pi^\ast\omega+dd^c\Phi\bigr)^{n+1}=0
\]
on \(\mathfrak M\) in the sense of currents. The bounds
\eqref{epsilon geodesic estimates} and Lemma~\ref{vdot estimate} hold for \(\varphi\) in the weak
sense, with the same type of constant. Furthermore, for every
\(\alpha\in(0,1)\) and every compact set
\(K\subset M\times[0,1]\),
\[
\lim_{\epsilon\to0}
\lVert\varphi^\epsilon-\varphi\rVert_{C^{1,\alpha}(K)}=0.
\]
We refer to this path as the Poincar\'e type \(C^{1,1}\) geodesic.
\end{lem}

\begin{proof}
From Lemma~\ref{vdot estimate}, we get
\[
\begin{split}
\lvert\dot{\varphi^\epsilon}(x,t)\rvert
&\le\lvert\varphi_1(x)-\varphi_0(x)\rvert+C.
\end{split}
\]
Then
\(\varphi_1-\varphi_0=O(\mathbf{u})\) by the definition of
\(\widetilde{\mathcal{PM}}_{\Omega}\). Combining the above estimates with the Lemma \ref{eps geodesic}, we can get that for any compact set $K \subset M \times [0,1]$, the $C^{1,1}$ norm of $\varphi^{\epsilon}$ is uniformly bounded on $K$. As a result, we can conclude the proof of the Lemma.
\end{proof}
We note that when we say $C^{1,1}$ norm, it means the complex Hessian bound. 

\subsection{Energy functionals}\label{subsection: Energy functionals}

We next record the energy functionals and their first variations on
\(\widetilde{\mathcal{PM}}_{\Omega}\). Define
\begin{equation}\label{mathcal E}
\mathcal{E}(\varphi)
\coloneqq
\int_M\varphi\sum_{j=0}^{n}
\omega_\varphi^{\,n-j}\wedge\omega^j.
\end{equation}
If \(T\) is a closed real \((1,1)\)-form satisfying
\(\lvert T\rvert_\omega\le C\), set
\[
\mathcal{E}^T(\varphi)
\coloneqq
\int_M\varphi\sum_{j=0}^{n-1}
\omega_\varphi^{\,n-j-1}\wedge\omega^j\wedge T.
\]
Let \(\mu_0=\omega^n\). For a positive measure \(\mu\) that is absolutely
continuous with respect to \(\mu_0\), define its relative entropy by
\[
H_{\mu_0}(\mu)
\coloneqq
\int_M\log\left(\frac{d\mu}{d\mu_0}\right)d\mu,
\]
whenever the integral is finite.

The mean scalar curvature is
\[
\bar R
\coloneqq
\frac{\int_M R_\varphi\,\omega_\varphi^n}
{\int_M\omega_\varphi^n}
=-2\pi n\,
\frac{c_1(K_X[D])\cdot[\omega]^{n-1}}{[\omega]^n}.
\]
With the convention \eqref{operator conventions}, the \(K\)-energy is
\begin{equation}\label{decom kenergy}
\mathcal{M}(\varphi)
\coloneqq
\frac{\bar R}{n+1}\mathcal{E}(\varphi)
-\mathcal{E}^{\operatorname{Ric}(\omega)}(\varphi)
+H_{\omega^n}(\omega_\varphi^n).
\end{equation}
The growth and quasi-isometry conditions in the definition of
\(\widetilde{\mathcal{PM}}_{\Omega}\) ensure that these integrals are
finite.

We shall use the following version of Gaffney's Stokes theorem
\cite{MR0062490}.

\begin{lem}\label{G-S}
Let \((N,g)\) be a complete oriented Riemannian manifold of real
dimension \(m\), with \(g\) of class \(C^2\). If
\(\Theta\in C^1(\Lambda^{m-1}T^\ast N)\) and both
\(\lvert\Theta\rvert_g\) and \(\lvert d\Theta\rvert_g\) are integrable,
then
\[
\int_N d\Theta=0.
\]
\end{lem}

\begin{lem}\label{functional first order}
Let
\(\varphi\in\widetilde{\mathcal{PM}}_{\Omega}\), and let
\(v=O(\mathbf{u})\) have bounded derivatives of every positive order
with respect to a Poincar\'e type metric. Then
\[
d\mathcal{E}|_\varphi(v)
=(n+1)\int_M v\,\omega_\varphi^n
\]
and
\[
d\mathcal{E}^T|_\varphi(v)
=n\int_M v\,\omega_\varphi^{n-1}\wedge T.
\]
Moreover,
\[
dH_{\omega^n}(\omega_\varphi^n)(v)
=
\int_M v\left[
-R_\varphi
+\operatorname{tr}_{\omega_\varphi}\operatorname{Ric}(\omega)
\right]\omega_\varphi^n.
\]
\end{lem}

\begin{proof}
Differentiating \(\mathcal E(\varphi+tv)\) at \(t=0\) gives
\begin{equation}\label{Eutv}
\begin{split}
d\mathcal E|_\varphi(v)
&=\sum_{j=0}^n\int_M
v\,\omega_\varphi^{n-j}\wedge\omega^j\\
&\quad
+\sum_{j=0}^{n-1}(n-j)\int_M
\varphi\,dd^c v\wedge
\omega_\varphi^{n-j-1}\wedge\omega^j.
\end{split}
\end{equation}
The stated growth assumptions make the relevant forms integrable, so
Lemma~\ref{G-S}, or equivalently the cutoff functions constructed in
Lemma~\ref{snc-cutoff}, justifies integration by parts:
\[
\int_M
\varphi\,dd^c v\wedge
\omega_\varphi^{n-j-1}\wedge\omega^j
=
\int_M
v\,dd^c\varphi\wedge
\omega_\varphi^{n-j-1}\wedge\omega^j.
\]
Since \(dd^c\varphi=\omega_\varphi-\omega\), the resulting sum
telescopes:
\[
\begin{split}
d\mathcal E|_\varphi(v)
&=\sum_{j=0}^n\int_M
v\,\omega_\varphi^{n-j}\wedge\omega^j\\
&\quad+\sum_{j=0}^{n-1}(n-j)\int_M
v\left(
\omega_\varphi^{n-j}\wedge\omega^j
-\omega_\varphi^{n-j-1}\wedge\omega^{j+1}
\right)\\
&=(n+1)\int_M v\,\omega_\varphi^n.
\end{split}
\]
The same calculation, with the closed form \(T\) included in every
term, yields
\[
d\mathcal E^T|_\varphi(v)
=n\int_M v\,\omega_\varphi^{n-1}\wedge T.
\]

Finally, set
\[
F_\varphi
\coloneqq
\log\left(\frac{\omega_\varphi^n}{\omega^n}\right).
\]
By \eqref{basic variation conventions},
\[
\left.\frac{d}{dt}\right|_{t=0}F_{\varphi+tv}
=-\Delta_\varphi v,
\qquad
\left.\frac{d}{dt}\right|_{t=0}\omega_{\varphi+tv}^n
=-(\Delta_\varphi v)\omega_\varphi^n.
\]
Using Lemma~\ref{G-S} and the self-adjointness of
\(\Delta_\varphi\) on the functions under consideration, we obtain
\[
\begin{split}
dH_{\omega^n}(\omega_\varphi^n)(v)
&=-\int_M(1+F_\varphi)\Delta_\varphi v\,\omega_\varphi^n\\
&=-\int_M v\Delta_\varphi F_\varphi\,\omega_\varphi^n\\
&=\int_M v\left[
-R_\varphi
+\operatorname{tr}_{\omega_\varphi}\operatorname{Ric}(\omega)
\right]\omega_\varphi^n,
\end{split}
\]
because
\[
\Delta_\varphi F_\varphi
=
R_\varphi
-\operatorname{tr}_{\omega_\varphi}\operatorname{Ric}(\omega).
\]
\end{proof}

\begin{cor}\label{cri mabuchi}
Under the hypotheses of Lemma~\ref{functional first order},
\[
d\mathcal{M}|_\varphi(v)
=
\int_Mv(\bar R-R_\varphi)\,\omega_\varphi^n.
\]
In particular, every Poincar\'e type cscK metric is a critical point of
\(\mathcal M\) on \(\widetilde{\mathcal{PM}}_{\Omega}\).
\end{cor}

\begin{proof}
Combine Lemma~\ref{functional first order} with
\[
n\,\operatorname{Ric}(\omega)\wedge\omega_\varphi^{n-1}
=
\operatorname{tr}_{\omega_\varphi}\operatorname{Ric}(\omega)\,
\omega_\varphi^n.
\]
\end{proof}
\subsection{Holomorphic vector fields}\label{subsection: Holomorphic vector fields}

We collect the notation for holomorphic vector fields and automorphism
groups used below.

\begin{defn}\label{definition: Holomorphic vector fields}
\begin{enumerate}
\item
\(\mathbf{h}_{\parallelsum}^D\) is the Lie algebra of holomorphic vector
fields on \(X\) tangent to every stratum of \(D\).

\item
\(\mathbf{h}_{\parallelsum,0}^D\) is the ideal in
\(\mathbf{h}_{\parallelsum}^D\) consisting of vector fields that vanish
along \(D\).

\item
With respect to the metric \(\omega\), set
\[
\mathbf{h}_{\parallelsum,\mathbb C}^D
\coloneqq
\left\{
V\in\mathbf{h}_{\parallelsum}^D:
V=\nabla_\omega^{1,0}f
\text{ for some complex-valued function }f
\right\}.
\]

\item
\(\mathbf{a}_{\parallelsum}^D(M)\) is the Lie subalgebra of
\(\mathbf{h}_{\parallelsum}^D\) consisting of holomorphic vector fields
that are parallel on \(M\) with respect to the Levi--Civita connection
of \(\omega\).

\item
Define
\[
\mathbf{h}_{\parallelsum,\mathbb R}^D
\coloneqq
\left\{
V\in\mathbf{h}_{\parallelsum}^D:
V=\nabla_\omega^{1,0}f
\text{ for some real-valued function }f
\right\}.
\]
We also write
\[
\overline{\mathbf{h}_{\parallelsum,\mathbb R}^D}
\coloneqq
\left\{
f\in C^\infty_{\mathbb R}(M):
\nabla_\omega^{1,0}f\in
\mathbf{h}_{\parallelsum,\mathbb R}^D
\right\}
\]
and
\[
\overline{\mathbf{h}_{\parallelsum,\mathbb R,0}^D}
\coloneqq
\left\{
f\in C^\infty_{\mathbb R}(M):
\nabla_\omega^{1,0}f\in
\mathbf{h}_{\parallelsum,0}^D
\right\}.
\]

\item
\(\mathbf{h}^D\) denotes the Lie algebra of holomorphic vector fields on
\(D\).

\item
\(\operatorname{Aut}_0(X)\) and \(\operatorname{Aut}_0(D)\) denote the
identity components of the biholomorphism groups of \(X\) and \(D\),
respectively.

\item
\(\operatorname{Aut}_0^D(X)\) is the identity component of the group of
biholomorphisms of \(X\) preserving \(D\), and
\[
\operatorname{Aut}_{0,\mathrm{id}}^D(X)
\coloneqq
\left\{
g\in\operatorname{Aut}_0^D(X):
\left.g\right|_D=\operatorname{id}_D
\right\}.
\]

\item
\(\operatorname{Iso}_0^D(X,\omega)\) is the subgroup of
\(\operatorname{Aut}_0^D(X)\) preserving \(\omega\).

\item
We use the notation
\[
\bar s=\bar R
=-2\pi n\,
\frac{c_1(K_X[D])\cdot[\omega]^{n-1}}{[\omega]^n},
\]
and, for each component \(D_j\),
\[
\bar s_{D_j}
=-2\pi(n-1)\,
\frac{
c_1(D_j)\cdot c_1(K_X[D])\cdot[\omega]^{n-2}
}{
c_1(D_j)\cdot[\omega]^{n-1}
}.
\]

\item
Let
\[
\varphi_0,\varphi_1\in\widetilde{\mathcal{PM}}_{\Omega},
\]
and let \(\varphi_t\) be the Poincar\'e type \(C^{1,1}\) geodesic
supplied by Lemma~\ref{C11 geodesic}. Put
\[
\omega_t=\omega+dd^c\varphi_t,\qquad
b_t
\coloneqq
\frac{\int_M\dot\varphi_t\,\omega_t^n}
{\int_M\omega_t^n}.
\]
The squared Mabuchi distance between the associated metrics is
\[
d(\omega_0,\omega_1)^2
\coloneqq
\int_0^1\int_M
\lvert\dot\varphi_t-b_t\rvert^2\,\omega_t^n\,dt.
\]
\end{enumerate}
\end{defn}
\section{Convexity of the K-energy}\label{section: Convexity of the K-energy}

We essentially follow the approach of the proof by Berman and Berndtsson \cite{BB}. Some modifications need to be made according to the Introduction.  Let us sketch the proof of the convexity of the $K$-energy. As in the Introduction, we can see a geodesic as a $S^1$ invariant function on $\mathfrak X=X\times R$, where $R$ is a cylinder. In this section, we will denote $d_R$ and $d_R^c$ as the differential operators on $R$. We denote $d_X$ and $d_X^c$ as the differential operators on $X$. We denote $d$ and $d^c$ as the differential operators on $\mathfrak X$. Firstly, for any $C^{1,1}$ Poincar\'e type geodesic $u_t$, we show that $$d_{R}d_{R}^c \mathcal{M}(u_{t})\ge 0$$ in the weak sense. Secondly, we prove that the $K$-energy is continuous along $u_t$. Then we can use the fact that the geodesic is $S^1$-symmetric to prove that the $K$-energy is in fact convex along $u_t$.

\subsection{Subharmonicity of the K-energy}
First we calculate the Hessian of some functionals along Poincar\'e type $C^{1,1}$ geodesics. 

\begin{lem}\label{ddce}
Let $U(x,t)=u_t(x)$ be a Poincar\'e type $C^{1,1}$ geodesic. Then we have that:
\begin{equation}
d_{R}d^c_{R}\mathcal{E}^{Ric_{\omega}}(u_{t})=\int_X(\pi^* \omega+dd^c U)^n \wedge \pi^*Ric_{\omega}, 
\end{equation}
and
\begin{equation}\label{linear E}
    d_{R}d^c_{R}\mathcal{E}(u_{t})=\int_X (\pi^* \omega+dd^c U)^{n+1}.
\end{equation}
\end{lem}
\begin{proof}
Let $u_{t}^{\epsilon}=U^{\epsilon}(x,t)$ be the $\epsilon$-geodesic in the Lemma \ref{eps geodesic}. By Lemma \ref{eps geodesic}, we have that $u_t^{\epsilon}-[(1-t)u_0+ tu_1]\in C^{\infty}$. This means that there exist constants $C(k,\epsilon)$ such that 
\begin{equation*}
    |\nabla_{\omega}^k\{u_t^{\epsilon}-[(1-t)u_0+ tu_1]\}|_{\omega} \le C(k,\epsilon)
\end{equation*}
Since all the derivatives of $u_0$ and $u_1$ are bounded with respect to $\omega$, we have that 
\begin{equation}\label{epsilon high estimate}
    |\nabla_{\omega}^k u_t^{\epsilon}|_{\omega} \le C_1(k,\epsilon)
\end{equation}
for $k \ge 1$ for some constants $C_1(k,\epsilon)$.

For any function $v\in C_0^{\infty}(R)$, we have that:
\begin{equation}\label{e 3.3}
\begin{split}
& d_R d_R^c \mathcal{E}^{Ric_{\omega}} (u^{\epsilon}_{t})(v)
= \int_{\mathfrak X} u^{\epsilon}_t \cdot \sum_{j=0}^{n-1}\omega_{u^{\epsilon}_t}^{n-j-1} \wedge \omega^j \wedge Ric_{\omega}  \wedge  d_Rd_R^c v \\
& =\frac{-1}{2} \int_{\mathfrak X} (\Delta_R v) \cdot u^{\epsilon}_t \cdot \sum_{j=0}^{n-1}\omega_{u^{\epsilon}_t}^{n-j-1} \wedge \omega^j \wedge Ric_{\omega}  \wedge \sqrt{-1}dw \wedge d\bar w \\
& =\frac{1}{4}\int_{\mathfrak X} v \frac{d^2}{dt^2}( u_t^{\epsilon}  \sum_{j=0}^{n-1}\omega_{u^{\epsilon}_t}^{n-j-1} \wedge \omega^j \wedge Ric_{\omega})  \wedge \sqrt{-1}dw \wedge d\bar w.
\end{split}
\end{equation}
The first line of the above formula follows because $d_R d_R^c \mathcal{E}^{Ric_{\omega}} (u^{\epsilon}_{t})$ is defined in the distribution sense. 

Since $u_t^{\epsilon}$ is smooth, we can calculate that:
\begin{align*}
&\frac{d}{dt} ( u_t^{\epsilon}  \sum_{j=0}^{n-1}\omega_{u^{\epsilon}_t}^{n-j-1} \wedge \omega^j \wedge Ric_{\omega}) \\
&= \dot{u}_t^{\epsilon} \sum_{j=0}^{n-1}\omega_{u^{\epsilon}_t}^{n-j-1} \wedge \omega^j \wedge Ric_{\omega} + u_t^{\epsilon} \sum_{j=0}^{n-2} (n-j-1) d_Xd_X^c \dot{u}_t^{\epsilon} \wedge \omega_{u_t^{\epsilon}}^{n-j-2} \wedge \omega^j \wedge Ric_{\omega},
\end{align*}
and
\begin{equation}\label{e 3.4.1}
\begin{split}
& \frac{d^2}{dt^2} ( u_t^{\epsilon}  \sum_{j=0}^{n-1}\omega_{u^{\epsilon}_t}^{n-j-1} \wedge \omega^j \wedge Ric_{\omega})= \ddot{u}_t^{\epsilon}\sum_{j=0}^{n-1} \omega_{u_t^{\epsilon}}^{n-j-1}\wedge \omega^j \wedge Ric_{\omega} \\
&+  \dot{u}_t^{\epsilon} \sum_{j=0}^{n-2} (n-j-1) d_Xd_X^c \dot{u}_t^{\epsilon} \wedge \omega_{u_t^{\epsilon}}^{n-j-2} \wedge \omega^j \wedge Ric_{\omega}\\
& + \dot{u}_t^{\epsilon} \sum_{j=0}^{n-2} (n-j-1) d_Xd_X^c \dot{u}_t^{\epsilon} \wedge \omega_{u_t^{\epsilon}}^{n-j-2} \wedge \omega^j \wedge Ric_{\omega}\\
& + u_t^{\epsilon}\sum_{j=0}^{n-2} (n-j-1) d_X d_X^c \ddot{u_t^{\epsilon}} \wedge \omega_{u_t^{\epsilon}}^{n-j-2}\wedge \omega^j \wedge Ric_{\omega} \\
&+ u_t^{\epsilon} \sum_{j=0}^{n-3} (n-j-1)(n-j-2) d_X d^c_X \dot{u}^{\epsilon}_t  \wedge d_X d^c_X \dot{u}^{\epsilon}_t \wedge \omega_{u_t^{\epsilon}}^{n-j-3} \wedge \omega^j \wedge Ric_{\omega}.
\end{split}
\end{equation}

Plugging in (\ref{e 3.4.1}) into (\ref{e 3.3}), we can get 
\begin{equation*}
\begin{split}
    & 4d_R d_R^c \mathcal{E}^{Ric_{\omega}} (u^{\epsilon}_{t})(v) \\
    & =\int_{\mathfrak X}v\cdot \bigg[ \ddot{u}_t^{\epsilon}\sum_{j=0}^{n-1} \omega_{u_t^{\epsilon}}^{n-j-1}\wedge \omega^j \wedge Ric_{\omega} \\
&\quad + 2 \dot{u}_t^{\epsilon} \sum_{j=0}^{n-2} (n-j-1) d_X d_X^c \dot{u}_t^{\epsilon} \wedge \omega_{u_t^{\epsilon}}^{n-j-2} \wedge \omega^j \wedge Ric_{\omega} \\
&\quad + u_t^{\epsilon}\sum_{j=0}^{n-2} (n-j-1) d_X d_X^c \ddot{u_t^{\epsilon}} \wedge \omega_{u_t^{\epsilon}}^{n-j-2}\wedge \omega^j \wedge Ric_{\omega} \\
&\quad+ u_t^{\epsilon} \sum_{j=0}^{n-3} (n-j-1)(n-j-2) d_X d^c_X \dot{u}^{\epsilon}_t  \wedge d_X d^c_X \dot{u}^{\epsilon}_t \wedge \omega_{u_t^{\epsilon}}^{n-j-3} \wedge \omega^j \wedge Ric_{\omega}\bigg]\\
&\quad \wedge \sqrt{-1} dw \wedge d \bar w.
\end{split}
\end{equation*}
Using (\ref{epsilon high estimate}), Lemma \ref{eps geodesic}  and Lemma \ref{vdot estimate}, all the terms in the above integral are integrable. 

Lemma \ref{G-S} implies that the third line in (\ref{e 3.3}) is valid. Then we can use Lemma \ref{G-S} to integrate by parts,
\begin{equation}\label{e 3.6}
\begin{split}
    &4 d_R d_R^c \mathcal{E}^{Ric_{\omega}} (u^{\epsilon}_{t})(v)\\
    &=\int_{\mathfrak X}v\cdot \bigg[ \ddot{u}_t^{\epsilon}\sum_{j=0}^{n-1} \omega_{u_t^{\epsilon}}^{n-j-1}\wedge \omega^j \wedge Ric_{\omega} \\
    & -  2 \sum_{j=0}^{n-2} (n-j-1)d_X \dot{u}_t^{\epsilon} \wedge  d_X^c \dot{u}_t^{\epsilon} \wedge \omega_{u_t^{\epsilon}}^{n-j-2} \wedge \omega^j \wedge Ric_{\omega} \\
    & + \sum_{j=0}^{n-2} (n-j-1)\ddot{u_t^{\epsilon}} d_X d_X^c u_t^{\epsilon}   \wedge \omega_{u_t^{\epsilon}}^{n-j-2}\wedge \omega^j \wedge Ric_{\omega} \\
    & - \sum_{j=0}^{n-3} (n-j-1)(n-j-2) d_X d_X^c u_t^{\epsilon} \wedge  d_X \dot{u}^{\epsilon}_t \wedge d_X^c \dot{u}^{\epsilon}_t \wedge \omega_{u_t^{\epsilon}}^{n-j-3}\wedge \omega^j \wedge Ric_{\omega}\bigg]\\
    &\wedge \sqrt{-1}dw \wedge d\bar w \\
    &=\int_{\mathfrak X}v\cdot \big[I+II+III+IV\big]\wedge \sqrt{-1}dw \wedge d\bar w .
\end{split}
\end{equation}

We also compute
\begin{equation}\label{e 3.7}
    \begin{split}
     &IV\\&=-\sum_{j=0}^{n-3} (n-j-1)(n-j-2) d_X d_X^c u_t^{\epsilon} \wedge  d_X \dot{u}^{\epsilon}_t \wedge d_X^c \dot{u}^{\epsilon}_t \wedge \omega_{u_t^{\epsilon}}^{n-j-3}\wedge \omega^j \wedge Ric_{\omega}\\
        &= - \sum_{j=0}^{n-3} (n-j-1)(n-j-2)(\omega_{u_t^{\epsilon}}- \omega) \wedge d_X \dot{u}^{\epsilon}_t \wedge d_X^c \dot{u}^{\epsilon}_t \wedge \omega_{u_t^{\epsilon}}^{n-j-3} \wedge \omega^j \wedge Ric_{\omega} \\
        & =-\sum_{j=0}^{n-3} (n-j-1)(n-j-2)\omega_{u_t^{\epsilon}}^{n-j-2} \wedge d_X \dot{u}^{\epsilon}_t \wedge d_X^c \dot{u}_t^{\epsilon} \wedge \omega^j \wedge Ric_{\omega} \\
        &\quad  +\sum_{j=1}^{n-2} (n-j)(n-j-1) \omega_{u_t^{\epsilon}}^{n-j-2} \wedge d_X \dot{u}_t^{\epsilon} \wedge d_X^c \dot{u}_t^{\epsilon} \wedge \omega^j \wedge Ric_{\omega}\\
        & =- (n-1)(n-2)\omega_{u_t^{\epsilon}}^{n-2} \wedge d_X \dot{u}_t^{\epsilon} \wedge d_X^c \dot{u}_t^{\epsilon} \wedge Ric_{\omega} \\
        & \quad + \sum_{j=1}^{n-3}2 (n-j-1)\omega_{u_t^{\epsilon}}^{n-j-2} \wedge d_X \dot{u}_t^{\epsilon} \wedge d_X^c \dot{u}_t^{\epsilon} \wedge \omega^j \wedge Ric_{\omega} \\
        &\quad  +2 d_X \dot{u}_t^{\epsilon} \wedge d_X^c \dot{u}_t^{\epsilon} \wedge \omega^{n-2} \wedge Ric_{\omega}.
    \end{split}
\end{equation}

We can simplify the above formula as follows:
\begin{equation}\label{e 3.71}
    \begin{split}
    &IV+II\\
    &= - \sum_{j=0}^{n-3} (n-j-1)(n-j-2) d_X d_X^c u_t^{\epsilon} \wedge  d_X \dot{u}^{\epsilon}_t \wedge d_X^c \dot{u}^{\epsilon}_t \wedge \omega_{u_t^{\epsilon}}^{n-j-3}\wedge \omega^j \wedge Ric_{\omega}\\
    &\quad  -2 \sum_{j=0}^{n-2}(n-j-1)d_X \dot{u}_t^{\epsilon} \wedge d_X^c \dot{u}^{\epsilon}_t \wedge \omega_{u_t^{\epsilon}}^{n-j-2} \wedge \omega^j \wedge Ric_{\omega}\\
        & =- (n-1)(n-2)\omega_{u_t^{\epsilon}}^{n-2} \wedge d_X \dot{u}_t^{\epsilon} \wedge d_X^c \dot{u}_t^{\epsilon} \wedge Ric_{\omega} \\
        & \quad + \sum_{j=1}^{n-3}2 (n-j-1)\omega_{u_t^{\epsilon}}^{n-j-2} \wedge d_X \dot{u}_t^{\epsilon} \wedge d_X^c \dot{u}_t^{\epsilon} \wedge \omega^j \wedge Ric_{\omega} \\
        & \quad +2 d_X \dot{u}_t^{\epsilon} \wedge d_X^c \dot{u}_t^{\epsilon} \wedge \omega^{n-2} \wedge Ric_{\omega}\\
        &\quad -2 \sum_{j=0}^{n-2}(n-j-1)d_X \dot{u}_t^{\epsilon} \wedge d_X^c \dot{u}^{\epsilon}_t \wedge \omega_{u_t^{\epsilon}}^{n-j-2} \wedge \omega^j \wedge Ric_{\omega}\\
      &= -n(n-1) \omega_{u_t^{\epsilon}}^{n-2} \wedge d_X \dot{u}^{\epsilon}_t \wedge d_X^c \dot{u}_t^{\epsilon} \wedge Ric_{\omega}.
    \end{split}
\end{equation}
We can also calculate that:
\begin{equation}\label{e 3.8}
\begin{split}
    I+III&=\ddot{u}_t^{\epsilon} \sum_{j=0}^{n-1}\omega_{u_t^{\epsilon}}^{n-j-1} \wedge \omega^j \wedge Ric_{\omega} \\
   &\quad +  \ddot{u}_t^{\epsilon} \sum_{j=0}^{n-2} (n-j-1)(\omega_{u_t^{\epsilon}}- \omega)\wedge \omega_{u_t^{\epsilon}}^{n-j-2} \wedge \omega^j \wedge Ric_{\omega}\\
    & =\ddot{u}_t^{\epsilon} \sum_{j=0}^{n-1}\omega_{u_t^{\epsilon}}^{n-j-1} \wedge \omega^j \wedge Ric_{\omega} +\ddot{u}_t^{\epsilon} \sum_{j=0}^{n-1} (n-j-1)\omega_{u_t^{\epsilon}}^{n-j-1}\wedge \omega^j \wedge Ric_{\omega} \\
    &\quad  -\ddot{u}_t^{\epsilon} \sum_{j=0}^{n-2}(n-j-1)\omega_{u_t^{\epsilon}}^{n-j-2}\wedge \omega^{j+1}\wedge Ric_{\omega}\\
    & =\ddot{u}_t^{\epsilon} \sum_{j=0}^{n-1}\omega_{u_t^{\epsilon}}^{n-j-1} \wedge \omega^j \wedge Ric_{\omega} +\ddot{u}_t^{\epsilon} \sum_{j=0}^{n-1} (n-j-1)\omega_{u_t^{\epsilon}}^{n-j-1}\wedge \omega^j \wedge Ric_{\omega} \\
    & \quad  -\ddot{u}_t^{\epsilon} \sum_{j=1}^{n-1}(n-j)\omega_{u_t^{\epsilon}}^{n-j-1}\wedge \omega^{j}\wedge Ric_{\omega}\\
    &= \ddot{u}_t^{\epsilon}\sum_{j=0}^{n-1} \omega_{u_t^{\epsilon}}^{n-j-1}\wedge \omega^j \wedge Ric_{\omega} +(n-1)\ddot{u}_t^{\epsilon} \omega_{u_t^{\epsilon}}^{n-1}\wedge Ric_{\omega} \\
    &\quad  - \ddot{u}_t^{\epsilon} \sum_{j=1}^{n-1} \omega_{u_t^{\epsilon}}^{n-j-1}\wedge \omega^j \wedge Ric_{\omega} \\
    & = n \ddot{u}_t^{\epsilon} \omega_{u_t^{\epsilon}}^{n-1}\wedge Ric_{\omega}.
\end{split}
\end{equation}

Therefore, we plug (\ref{e 3.71}) and (\ref{e 3.8}) into (\ref{e 3.6}) to get that:
\begin{equation}\label{ddceric eps}
d_R d_R^c \mathcal{E}^{Ric_{\omega}}(u_t^{\epsilon})(v) = \int_{\mathfrak X} v(\pi^* \omega + dd^c U^{\epsilon} )^{n} \wedge \pi^* Ric_{\omega}.
\end{equation}

Now we want to let $\epsilon$ go to zero. We first show that $\mathcal{E}^{Ric_{\omega}}(u_t^{\epsilon})$ converges uniformly to $\mathcal{E}^{Ric_{\omega}}(u_t)$ if we see them as functions on $R$. In fact, we can calculate that:
\begin{equation*}
\begin{split}
& \mathcal{E}^{Ric_{\omega}}(u_t^{\epsilon}) - \mathcal{E}^{Ric_{\omega}}(u_t) \\
&=\int_X \sum_{j=0}^{n-1} (u^{\epsilon}-u) \omega_{u^{\epsilon}}^{n-j-1} \wedge \omega^j \wedge Ric_{\omega}  +\int_X u \sum_{j=0}^{n-1} (\omega_{u^{\epsilon}}^{n-j-1} -\omega_u^{n-j-1}) \wedge \omega^j\wedge Ric_{\omega} \\
& =\int_{\mathcal{N}_A}\sum_{j=0}^{n-1} (u^{\epsilon}-u) \omega_{u^{\epsilon}}^{n-j-1} \wedge \omega^j \wedge Ric_{\omega}  + \int_{X \setminus \mathcal{N}_A}\sum_{j=0}^{n-1} (u^{\epsilon}-u)\omega_{u^{\epsilon}}^{n-j-1} \wedge \omega^j \wedge Ric_{\omega}  \\
&\quad +\int_X u \sum_{j=0}^{n-1} (\omega_{u^{\epsilon}}^{n-j-1} -\omega_u^{n-j-1}) \wedge \omega^j \wedge Ric_{\omega}.
\end{split}
\end{equation*}

Recall that $\mathcal{N}_A$ is a neighbourhood of $D$ such that it shrinks to $D$ as $A \rightarrow \infty$. Here $A$ will be determined later on. By the Lemma \ref{eps geodesic} and the Lemma \ref{C11 geodesic}, we have that there exists a uniform $C$ independent of $\epsilon$ such that $|u^{\epsilon}|+|u|\le C \cdot \mathbf{u}$ and $\omega_{u^{\epsilon}}\le C \cdot \omega$. Since by the definition of Poincar\'e type K\"ahler metrics, we have that $|Ric_{\omega}|_{\omega}\le C$ for some constant $C$. Then we have that:
\begin{equation*}
|\int_{\mathcal{N}_A}\sum_{j=0}^{n-1} (u^{\epsilon}-u) \omega_{u^{\epsilon}}^{n-j-1} \wedge \omega^j \wedge Ric_{\omega}  | \le C \int_{\mathcal{N}_A} \mathbf{u}\, \omega^n. 
\end{equation*}
Fix $\delta>0$. We can use the fact that $\mathbf{u}\in L^1(\omega^n)$ to get that there exists an $A$ such that $$C\int_{\mathcal{N}_A} \mathbf{u}\, \omega^n\le 1/3 \delta.$$ We fix this $A$.

Using the Lemma \ref{C11 geodesic} again we can get that $U^{\epsilon}$ converges in $C^{1,\alpha}_{loc}(M \times R)$ to $U$ for any $\alpha \in (0,1)$. In particular, $U^{\epsilon}$ converge in $C^{1,\alpha}((X\setminus N_A) \times R)$ to $U$. This implies that for $\epsilon$ small enough:
\begin{equation*}
    |\int_{X\setminus \mathcal{N}_A}\sum_{j=0}^{n-1} (u^{\epsilon}-u) \, \omega_{u^{\epsilon}}^{n-j-1} \wedge \omega^j \wedge Ric_{\omega}  | \le  C ||u-u^{\epsilon}||_{L^{\infty}(X\setminus \mathcal{N}_A)} \int_{X\setminus \mathcal{N}_A} \omega^n \le 1/3 \delta.
\end{equation*}

For the third term $\int_X u \sum_{j=0}^{n-1} (\omega_{u^{\epsilon}}^{n-j-1} -\omega_u^{n-j-1}) \wedge \omega^j \wedge Ric_{\omega}$, we can estimate as follows:
\begin{equation*}
\begin{split}
    &\int_X u \cdot\sum_{j=0}^{n-1} (\omega_{u^{\epsilon}}^{n-j-1} -\omega_u^{n-j-1}) \wedge \omega^j \wedge Ric_{\omega}\\
    &= \int_X u \cdot\sum_{j=0}^{n-1} d_Xd_X^c (u^{\epsilon}-u )\cdot \sum_{l=0}^{n-j-2} \omega_{u^{\epsilon}}^l \wedge \omega_u^{n-j-l-2} \wedge \omega^j \wedge Ric_{\omega} \\
    & = \int_X (u^{\epsilon}-u) \cdot d_X d_X^c u \cdot \sum_{j=0}^{n-1} \sum_{l=0}^{n-j-2} \omega_{u^{\epsilon}}^l \wedge \omega_u^{n-j-l-2} \wedge \omega^j  \wedge Ric_{\omega}.
\end{split}
\end{equation*}
The absolute value of this term can be proved to be smaller than $1/3 \delta$ in the same way as before, if we let $\epsilon$ be small enough. 

In conclusion, $\mathcal{E}^{\operatorname{Ric}(\omega)}(u_t^{\epsilon_k})-\mathcal{E}^{\operatorname{Ric}(\omega)}(u_t)$ converges uniformly to zero. Since the associated currents converge when the potentials converge uniformly, we obtain
\begin{equation}\label{ddceric converge}
    d_R d_R^c \mathcal{E}^{Ric_{\omega}}(u_t^{\epsilon_k}) \rightarrow d_R d_R^c \mathcal{E}^{Ric_{\omega}}(u_t)
\end{equation}
in the sense of distributions, meaning that for any $v\in C_0^{\infty}(R)$, we have that:
\begin{align*}
&\lim_{\epsilon_k \rightarrow 0} d_R d_R^c \mathcal{E}^{Ric_{\omega}}(u_t^{\epsilon_k}) (v)= \lim_{\epsilon_k \rightarrow 0}  \int_R  ( \mathcal{E}^{Ric_{\omega}} (u_t^{\epsilon_k})d_R d_R^c v )\, dvol_R\\
&= \int_R   (\mathcal{E}^{Ric_{\omega}} (u_t)d_R d_R^c v) \, dvol_R = d_R d_R^c \mathcal{E}^{Ric_{\omega}}(u_t) (v).
\end{align*}
Next, we calculate that 
\begin{equation*}
\begin{split}
& \int_{\mathfrak X} v\cdot(\pi^* \omega + dd^c U^{\epsilon} )^{n} \wedge \pi^* Ric_{\omega}- \int_{\mathfrak X} v\cdot(\pi^* \omega + dd^c U )^{n} \wedge \pi^* Ric_{\omega}    \\
& =\int_{\mathfrak X} v\cdot dd^c (U^{\epsilon}-U) \wedge \sum_{l=0}^{n-1} (\pi^* \omega + dd^c U)^l \wedge (\pi^* \omega +dd^c U_{\epsilon})^{n-1-l} \wedge \pi^* Ric_{\omega}\\
& =\int_{\mathfrak X} (U^{\epsilon}-U)\cdot dd^c v  \wedge \sum_{l=0}^{n-1} (\pi^* \omega + dd^c U)^l \wedge (\pi^* \omega +dd^c U_{\epsilon})^{n-1-l} \wedge \pi^* Ric_{\omega}.
\end{split}
\end{equation*}
The integral in the last line above goes to zero with $\epsilon=\epsilon_k$. This can be proved in the same way as before by considering the integral over $\mathcal{N}_A \times R$ and $(X \setminus \mathcal{N}_A) \times R$. 

Combining this fact with (\ref{ddceric converge}) and (\ref{ddceric eps}), we have that:
\begin{equation*}
    d_R d_R^c \mathcal{E}^{Ric_{\omega}}(u_t)(v) = \int_{\mathfrak X} v\cdot(\pi^* \omega + dd^c U )^{n} \wedge \pi^* Ric_{\omega}.
\end{equation*}
This finishes the proof of the first part of the Lemma. The second part can be proved in a very similar way.
\end{proof}

\begin{lem}\label{subharmonic}
Suppose that we have a function:
     \begin{equation*}
        f^{\Psi}(t) \triangleq \frac{\bar R}{n+1}\mathcal{E}(u_{t})-\mathcal{E}^{\operatorname{Ric}(\omega)}(u_t)+\int_X \psi_t\,\omega_{u_t}^n.
    \end{equation*}
    Here $$\Psi=\Psi(t,x)=\psi_t(x)= O(\mathbf{u})$$ is a function defined on $M\times R$ which is independent of the imaginary part of the coordinate $w$ of $R$. Then we have that 
\begin{equation*}
d_Rd_R^c f^{\Psi}(\tau)=\int_X \left[dd^c ( \psi_t (\pi^* \omega +dd^c U)^n) -(\pi^* \omega+dd^c U)^n \wedge \pi^* \operatorname{Ric}(\omega)\right].
\end{equation*}
\end{lem}
\begin{proof}
First we claim that there exists a sequence of smooth functions $\Psi_{j}$ such that $\supp \Psi_j \subset \subset M$ that converge almost everywhere to $\Psi$ on $M\times R$. We also require that $$|\Psi_{j}|\le C \, \mathbf{u}$$ for some uniform constant $C$. In fact, we can define $\Psi_j$ as the regularisation of $\Psi \chi_{ (X\setminus \mathcal{N}_j) \times R}$, where $\chi_{ (X\setminus \mathcal{N}_j) \times R}$ denotes the characteristic function of $(X\setminus \mathcal{N}_j) \times R$. Since 
$$\supp (\Psi \chi_{ (X\setminus \mathcal{N}_j) \times R}) \subset \subset  M\times R,$$
 the regularisation of $\Psi \chi_{ (X\setminus \mathcal{N}_j)\times R}$ can still be compactly supported in $M\times R$.  Since we assume that $$|\Psi| \le C_0 \,\mathbf{u},$$ the regularisation of  $\Psi \chi_{ (X\setminus \mathcal{N}_j)\times R}$ can be assumed to be bounded by $(C_0+1)\mathbf{u}$ from above and $-(C_0+1)\mathbf{u}$ from below.  So the claim is proved. 

Denote $\psi_j$ as the restriction of $\Psi_j$ on each fibre $M \times \{w\}$. For any $v\in C^{\infty}(R)$ which is compactly supported in the interior of $R$, we can calculate that:
\begin{equation}\label{e 3.2.1}
    d_R d_R^c \int_X \psi_{j}\,\omega_{u_t}^n (v) \triangleq \int_R  \int_X \psi_{j}\,\omega_{u_t}^n d_R d^c_R v = \int_{\mathfrak X} \psi_j \, (\pi^* \omega + dd^c U)^n \wedge dd^c v.
\end{equation}
The last equality holds because $\Psi_j$ is smooth and compactly supported in $M\times R$, so the calculation is the same as on a closed manifold with a smooth K\"ahler metric. Since $|\psi_j|\le C\mathbf{u}$, the dominated convergence theorem gives
\begin{equation}\label{e 3.2.2}
    \lim_{j \rightarrow \infty} \int_{\mathfrak X} \psi_j \,(\pi^* \omega + dd^c U)^n \wedge dd^c v = \int_{\mathfrak X} \psi \,(\pi^* \omega + dd^c U)^n \wedge dd^c v.
\end{equation}
We can also use the dominated convergence theorem to get that:
\begin{equation}\label{e 3.2.3}
\begin{split}
   \lim_{j \rightarrow \infty} d_R d_R^c (\int_X \psi_j \,\omega_{u_t}^n) (v) &=\lim_{j \rightarrow \infty} \int_{\mathfrak X} \psi_j \,\omega_{u_t}^n \wedge d_R d_R^c v\\
   & = \int_{\mathfrak X} \psi \,\omega_{u_t}^n \wedge d_R d_R^c v =d_R d_R^c (\int_X \psi_t \,\omega_{u_t}^n) (v).
   \end{split}
\end{equation}

Then we can combine (\ref{e 3.2.1}), (\ref{e 3.2.2}) and (\ref{e 3.2.3}) to get that:
\begin{equation*}
    d_R d_R^c \int_X \psi_t\, \omega_{u_t}^n (v) = \int_{\mathfrak X} \psi  \, (\pi^* \omega +dd^c U)^n \wedge dd^c v.
\end{equation*}
Combining this formula with the Lemma \ref{ddce}, we conclude the proof of this lemma.
\end{proof}

Note that $(\omega+dd^c u_t)$ can be degenerate. In order to remove the degeneracy, we use the following modification of $\Psi_A$:
\begin{equation}\label{def psia}
    \Psi_A =\max\{\log\frac{(\omega+dd^c u_t)^n}{\omega^n}, \chi-A\},
\end{equation}
where $\chi$ is a properly chosen continuous function on $M\times R$ which satisfies that $\chi=O(\mathbf{u})$. Then we have the following lemma:
\begin{lem}\label{posi}
   Let $k_0$ be a constant such that $\pi^* Ric_{\omega}-k_0 \pi^* \omega \le 0$. Let $\chi$ be $-k_0 U$. Then we have that 
    \begin{equation}\label{positivity}
         dd^c [ \psi_A (\pi^* \omega +dd^c U)^n] -(\pi^* \omega+dd^c U)^n \wedge \pi^* Ric_{\omega} \ge 0 \text{ on }\mathfrak X.
    \end{equation}
\end{lem}
\begin{proof}
First we claim that $\chi-A$ satisfies (\ref{positivity}). In fact, we have that 
    \begin{align*}
    dd^c(\chi-A) &=dd^c \chi =-k_0 dd^c U =-k_0(\pi^* \omega + dd^c U)+k_0 \pi^* \omega \\
    &\ge -k_0(\pi^* \omega + dd^c U) +\pi^* Ric_{\omega}.
    \end{align*}
    Then the claim follows by wedging the above formula with $(\pi^* \omega+ dd^c U)^n$ and the fact that $U$ is a geodesic.

Note that in the proof of \cite[Theorem 3.3]{BB}, Berman and Berndtsson use the local Bergman approximation to show that (\ref{positivity}) holds which is a purely local method. Since we are dealing with Poincar\'e type K\"ahler metrics which is smooth away from $D$, we can use the same method as Berman and Berndtsson to show that (\ref{positivity}) holds on $\mathfrak M=M \times R$. 

Now we want to show that (\ref{positivity}) holds on the entire manifold $X$. That is, for any smooth nonnegative function $v$ on $\mathfrak X$, we want to show that 
\begin{equation*}
\int_{\mathfrak X} v \cdot dd^c [ \psi_A (\pi^* \omega +dd^c U)^n] -v\cdot (\pi^* \omega+dd^c U)^n \wedge \pi^* Ric_{\omega} \ge 0.
\end{equation*}
Since $M$ with a Poincar\'e type K\"ahler metric is complete and noncompact, for any $m$ we can define a cut-off function $\eta_{m}$ as in the Lemma \ref{snc-cutoff}. Then, we have that 
\begin{equation*}
\begin{split}
&\int v  \cdot  dd^c [ \psi_A (\pi^* \omega +dd^c U)^n]-v\cdot  (\pi^* \omega+dd^c U)^n \wedge \pi^* Ric_{\omega}\\
& = \int dd^c v \wedge \psi_A (\pi^* \omega +dd^c U)^n  - v \cdot (\pi^* \omega+dd^c U)^n \wedge \pi^* Ric_{\omega} \\
& = \int dd^c [(1-\eta_{m})v +\eta_{m} v] \wedge \psi_A(\pi^* \omega +dd^c U)^n \\
& \quad - [(1-\eta_{m})v +\eta_{m} v  ] \cdot (\pi^* \omega+dd^c U)^n \wedge \pi^* Ric_{\omega}.
\end{split}
\end{equation*}

Since (\ref{positivity}) holds on $(X \setminus D) \times R$, we have that
\begin{equation*}
 \int dd^c [(1-\eta_{m})v] \wedge \psi_A(\pi^* \omega +dd^c U)^n - (1-\eta_{m})v \cdot  (\pi^* \omega+dd^c U)^n \wedge \pi^* Ric_{\omega} \ge 0.
\end{equation*}
We also have that:
\begin{align*}
&\int \psi_A  \cdot  (\pi^* \omega +dd^c U)^n \wedge dd^c (\eta_{m}v) \\
&= \int \psi_A \cdot (\pi^* \omega +dd^c U)^n \wedge [ dd^c (\eta_{m})v + \eta_{m} dd^c v + d\eta_{m} \wedge d^c v  +d v \wedge d^c \eta_{m} ].
\end{align*}

Since $|\psi_A| \leq C_A + C u$ and $(1+u)$ is integrable against the Poincar\'e volume, the measure $|\psi_A| (\pi^* \omega +dd^c U)^n$ has finite mass. Hence every term containing $d \eta_m$ or $dd^c \eta_m$ is bounded by $C/m$ times the mass of the annulus ${m<T<2m}$, and tends to zero.


Since $\eta_m$ decreases to zero, the second term also tends to zero. Indeed, set 
$$\Sigma=\{x\in \mathfrak X: \psi_A \cdot (\pi^* \omega+dd^c U)^n \wedge dd^c v (x) \text{ is positive}\}.$$ The monotone convergence theorem gives
\begin{equation*}
    \lim_{m \rightarrow \infty}\int_{\Sigma} \eta_{m} \psi_A  \cdot  (\pi^* \omega +dd^c U)^n \wedge dd^c v  =0
\end{equation*}
and
\begin{equation*}
    \lim_{m \rightarrow \infty}\int_{X\setminus \Sigma} \eta_{m} \psi_A  \cdot  (\pi^* \omega +dd^c U)^n \wedge dd^c v  =0.
\end{equation*}

Similarly, we have that $$\int \eta_{m} v \cdot  (\pi^* \omega +dd^c U)^n \wedge \pi^* Ric_{\omega}$$ goes to zero  as $m$ goes to $\infty$. Combining the above arguments, we have that 
\begin{equation*}
\int v  \cdot  dd^c [ \psi_A(\pi^* \omega +dd^c U)^n] -v \cdot (\pi^* \omega+dd^c U)^n \wedge \pi^* Ric_{\omega} \ge 0.
\end{equation*}
This finishes the proof of this lemma.
\end{proof}

\begin{thm} \label{subharmonicity}
Let $u_{t}$ be a Poincar\'e type $C^{1,1}$ geodesic. Then the Mabuchi functional $\mathcal{M}(u_{t})$ is weakly subharmonic when it is seen as a function on $R$.
\end{thm}
\begin{proof}
Using the Lemma \ref{subharmonic} and Lemma \ref{posi}, we have that $f^{\Psi_A}$ is subharmonic. Using the monotone convergence theorem, we have that $f^{\Psi_A}$ decrease to $f^{\Psi}$ as $A \rightarrow \infty$. Because the subharmonicity is preserved under the decreasing sequence, we have that $f^{\Psi}$ is also subharmonic on $R$. This concludes the proof of this theorem.
\end{proof}

\subsection{Continuity of the K-energy}

Next, we need to show that the $K$-energy is continuous.

\begin{lem}\label{lem cont e}
$\mathcal{E}$ and $\mathcal{E}^{Ric_{\omega}}$ are continuous  along Poincar\'e type $C^{1,1}$ geodesics.
\end{lem}
\begin{proof}
Let $u_t$ be a Poincar\'e type $C^{1,1}$ geodesic.  First, we show that $\mathcal{E}(u_t)$ is continuous. We can use the Lemma \ref{functional first order} to calculate that:
\begin{equation}\label{difference E}
    \mathcal{E}(u_{t+s})-\mathcal{E}(u_{t})= \int_t^{t+s}d \mathcal{E}|_{u_{\lambda}}(\dot{u}_{\lambda})d\lambda =\int_t^{t+s} 
    (n+1)\int_M \dot{u}_{\lambda} \omega_{u_{\lambda}}^n.
\end{equation}

Using the Lemma \ref{C11 geodesic} and the Lemma \ref{vdot estimate}, we have that $$(n+1)\Big|\int_M \dot{u}_{\lambda} \omega_{u_{\lambda}}^n\Big|$$ is uniformly bounded with respect to $\lambda$. Then we can use (\ref{difference E}) to get that:
\begin{equation*}
    \lim_{s \rightarrow 0} \mathcal{E}(u_{t+s})-\mathcal{E}(u_{t}) =0.
\end{equation*}

Next we talk about the continuity of $\mathcal{E}^{Ric_{\omega}}$. Using Lemma \ref{functional first order} we have that:
\begin{equation}\label{diff ERic}
    \mathcal{E}^{Ric_{\omega}}(u_{t+s})- \mathcal{E}^{Ric_{\omega}}(u_{t})=\int_t^{t+s}d\mathcal{E}^{Ric_{\omega}}\Big|_{u_{\lambda}}(\dot{u_{\lambda}}) d\lambda\\
    = \int_t^{t+s} n\int_M \dot{u_{\lambda}}\omega_{u_{\lambda}}^{n-1} \wedge \pi^* Ric_{\omega}.
\end{equation}
Again using the Lemma \ref{C11 geodesic} and the Lemma \ref{vdot estimate}, we obtain that $$n\Big|\int_M \dot{u}_{\lambda}\omega_{u_{\lambda}}^{n-1} \wedge \pi^* Ric_{\omega}\Big|$$ is uniformly bounded with respect to $\lambda$. Then we can let $s$ go to zero in (\ref{diff ERic}) and get:
\begin{equation*}
    \lim_{s \rightarrow 0} \mathcal{E}^{Ric_{\omega}}(u_{t+s})- \mathcal{E}^{Ric_{\omega}}(u_{t})=0.
\end{equation*}
\end{proof}

\begin{prop}\label{H upp semicotinuous}
The entropy term $H_{\omega^n}$ is upper-semicontinuous along Poincar\'e type $C^{1,1}$ geodesic.
\end{prop}

\begin{proof}
Let $U_k=\{p\in X\setminus D:d(p,D)\le 2^{-k}\}$, where $d$ is the distance induced by a smooth K\"ahler metric on $X$. For every $k$, cover $U_k\setminus U_{k+1}$ by finitely many balls of radius $2^{-k-2}$, and cover $X\setminus U_1$ by finitely many balls whose closures are disjoint from $D$. Denote all these balls by $\{V_i\}_{i=1}^{\infty}$. Then $\{V_j\}_{j=1}^{\infty}$ is a locally finite cover of $X\setminus D$ and
    \begin{equation*}
    X\setminus D = \cup_{j=1}^{\infty} V_j.
    \end{equation*}
     Let $\xi_j$ be a partition of unity subordinate to $\{V_j\}_{j=1}^{\infty}$. 
     
  Step 1.    We want to prove that
\begin{equation*}
        \widetilde{H}_j \triangleq \int_M \xi_j \,\log(\frac{\omega_{u_t}^n}{\omega^n}) \,\omega_{u_t}^n
    \end{equation*}
    is upper semicontinuous. The support of each partition-of-unity function $\xi_j$ has positive distance from $D$, so $\omega^n$ is smooth there. Thus the preceding argument is the same as the smooth argument of Berman--Berndtsson \cite{BB}; we include the details for the reader's convenience.
     
     As before, we use $\int_M \log (\frac{e^{\psi_{A}}\omega^n}{\omega^n})\omega_{u_t}^n$ to do approximation.  Define $\kappa_{\epsilon}(s)=s+\epsilon e^s$ and
    \begin{equation*}
        H_j=\int_M \xi_j \,\kappa_{\epsilon}\,(\log(\frac{e^{\psi_{At}}\omega^n}{\omega^n}))\,\omega_{u_t}^n.
    \end{equation*}
We can define $H_j^{(k)}$ in a similar way, replacing $\Psi_A$ by its kth approximation by local Bergman measure on $V_j$. Note that $V_j$ has positive distance from $D$. So such approximation is the same as the smooth case and has nothing to do with Poincar\'e type singularity. According to \cite{BB}, we have that:
    \begin{equation*}
        dd^c H_j^{(k)}\ge -C_{\epsilon,j} \sqrt{-1} dw\wedge d\bar w .
    \end{equation*}

    Since the local Bergman kernel is continuous, we have that $H_j^{(k)}$ is continuous. Since $H_j^{(k)}$ is $S^1$ invariant, we have that $H_j^{(k)}+C_{\epsilon,j}t^2$ is also convex. 
    
 Since 
    \begin{equation*}
        \lim_{k \rightarrow \infty} H_j^{(k)}=H_j,
    \end{equation*}
    we have that $H_j+C_{\epsilon,j}t^2$ is also convex and therefore upper semicontinuous. Hence $H_j$ is upper semicontinuous.
    
    Recall the definition of $\Psi_A$:
    \begin{equation*}
        \Psi_A =\max\Big\{\log \frac{(\omega+dd^c u_t)^n}{\omega^n}, \chi-A\Big\},
    \end{equation*}
    we have that $\Psi_A$ is decreasing with respect to $A$ and $$\lim_{A\rightarrow \infty}\Psi_A = \log \frac{(\omega+dd^c u_t)^n}{\omega^n}.$$ As a result, $H_j$ is decreasing with respect to $A$ and the limit is 
    \begin{equation*}
        H_j^* \triangleq \int_M \xi_j \,\kappa_{\epsilon}\, (\log(\frac{\omega_{u_t}^n}{\omega^n}))\,\omega_{u_t}^n
    \end{equation*}
    Then $H_j^*$ is upper semicontinuous because this property is preserved under decreasing limits. Since $H_j^*$ decreases as $\epsilon$ tends to zero, its limit $\widetilde{H}_j$ is also upper semicontinuous.

     Step 2. We show that the entropy term $H_{\omega^n}$ is also upper semicontinuous. We begin with a heuristic argument.
     
     Suppose that $\frac{\omega_{u_t}^n}{\omega^n}\le 1$, then we have that $\log  \frac{\omega_{u_t}^n}{\omega^n}$ is non-positive. Then, by the monotone convergence theorem, entropy term can be written as
     \begin{equation*}
     \begin{split}
         \int_M \log(\frac{\omega_{u_t}^n}{\omega^n})\omega_{u_t}^n 
         &= \int_M \lim_{l \rightarrow \infty}\sum_{j=1}^l  \xi_j \log(\frac{\omega_{u_t}^n}{\omega^n})\omega_{u_t}^n\\
& = \lim_{l\rightarrow \infty}\sum_{j=1}^l \int_M  \xi_j \log(\frac{\omega_{u_t}^n}{\omega^n})\omega_{u_t}^n \\
         & = \lim_{l \rightarrow \infty}\sum_{j=1}^l \widetilde{H}_j.
         \end{split}
     \end{equation*}
Under the assumption $\omega_{u_t}^n/\omega^n\le1$, each $\widetilde{H}_j$ is nonpositive. Therefore, $\sum_{j=1}^l\widetilde{H}_j$ decreases with $l$. Each partial sum is upper semicontinuous, and so is its decreasing limit $$\int_M \log(\frac{\omega_{u_t}^n}{\omega^n})\omega_{u_t}^n.$$ 
     
     We now prove upper semicontinuity of the entropy. Although $\omega_{u_t}^n/\omega^n$ need not be at most one, Lemma \ref{C11 geodesic} gives a uniform constant $C$ such that $\omega_{u_t}^n/\omega^n\le C$. Hence
     \begin{equation*}
     \begin{split}
 &\int_M \log(\frac{\omega_{u_t}^n}{\omega^n})\omega_{u_t}^n = \int_M \log(\frac{\omega_{u_t}^n}{C \omega^n})\omega_{u_t}^n + \int_M \log C \omega_{u_t}^n
     \end{split}         
     \end{equation*}
   Since $\omega_{u_t}^n/(C\omega^n)\le1$, the preceding argument shows that $\int_M \log(\omega_{u_t}^n/(C\omega^n))\omega_{u_t}^n$ is upper semicontinuous. The term $\int_M\log(C)\omega_{u_t}^n$ is continuous by Lemma \ref{lem cont e}. Thus the entropy is upper semicontinuous.
\end{proof}

\begin{lem}\label{entropy lower con}
    Let $u_t$ be a Poincar\'e type geodesic. Then we have that $f_s \triangleq \frac{\omega_{u_{t+s}}^n}{\omega^n}$ converge to  $ f \triangleq \frac{\omega_{u_{t}}^n}{\omega^n}$ weakly in $L^1(\omega^n)$ sense as $s \rightarrow 0$.
\end{lem}
\begin{proof}
    First, since $u_{t+s}-u_t=\int_t^{t+s}\dot{u}_\lambda d\lambda$, for $v\in C^2(X)$, we can calculate that:
     \begin{equation*}
     \begin{split}
  \int_M \omega_{u_{t+s}}^n v -\int_M \omega_{u_{t}}^n v 
  &= \int_M v \cdot dd^c(u_{t+s}-u_t)\cdot \sum_{i=0}^{n-1} \omega_{u_{t+s}}^i \wedge \omega_{u_t}^{n-1-i}\\
    &=\int_M (u_{t+s}-u_t)\cdot  dd^c v \cdot \sum_{i=0}^{n-1} \omega_{u_{t+s}}^i \wedge \omega_{u_t}^{n-1-i}\\
& =\int_M \int_t^{t+s}\dot{u}_\lambda d\lambda\cdot dd^c v \cdot \sum_{i=0}^{n-1} \omega_{u_{t+s}}^i \wedge \omega_{u_t}^{n-1-i}.
\end{split}
\end{equation*}
After taking absolute values, we get
     \begin{equation*}
     \begin{split}
     \Big| \int_M \omega_{u_{t+s}}^n v -\int_M \omega_{u_{t}}^n v\Big| 
& \le \int_M \Big|\int_t^{t+s}\dot{u}_\lambda d\lambda\Big| \, \Big|dd^c v\Big|_{\omega} \omega \wedge  \sum_{i=0}^{n-1} \omega_{u_{t+s}}^i \wedge \omega_{u_t}^{n-1-i} \\
& \le \int_M C|s|\, \mathbf{u} \,\Big|dd^c v\Big|_{\omega} \omega \wedge  \sum_{i=0}^{n-1} \omega_{u_{t+s}}^i \wedge \omega_{u_t}^{n-1-i}\\
&\le C|s|\int_M \mathbf{u}\,\Big |dd^c v\Big|_{\omega} \omega^n.
\end{split}
\end{equation*}
In the above calculation we use the Lemma \ref{C11 geodesic} and the Lemma \ref{vdot estimate}.

Then we have that 
\begin{equation}\label{e 3.5.1}
    \lim_{s\rightarrow 0} \int_X f_s v \omega^n = \int_X f v \omega^n, \text{ for any } v\in C^2(X). 
\end{equation}

For $v\in L^{\infty}(X)$, we can find a sequence of functions $v_k\in C^2(X)$ such that $v_k$ are uniformly bounded and converge to $v$ almost everywhere. We  calculate that:
\begin{equation*}
    \Big|\int_X \omega_{u_{t+s}}^n (v_k-v )\Big|+ \Big|\int_X \omega_{u_{t}}^n (v_k-v )\Big|\le C \int_X \Big|v_k-v\Big| \omega^n.
\end{equation*}
Here we again use the $C^{1,1}$ estimate for the geodesic:  $\omega_{u_t}\le C \omega$ and $\omega_{u_{t+s}}\le C \omega$. Then we can use the dominated convergence theorem to get that  
\begin{equation}\label{e 3.5.2}
    \lim_{\epsilon \rightarrow 0} \Big|\int_X \omega_{u_{t+s}}^n (v_k-v )\Big|+ \Big|\int_X \omega_{u_{t}}^n (v_k-v )\Big| =0
\end{equation}
uniformly with respect to $s$. Combining (\ref{e 3.5.1}) and (\ref{e 3.5.2}), we get that 
\begin{equation*}
    \lim_{s\rightarrow 0} \int_X f_s v \omega^n = \int_X f v \omega^n 
\end{equation*}
for any $v\in L^{\infty}(X)$. This shows that $f_s$ converge to $f$ weakly in $L^1(\omega^n)$.
\end{proof}

\begin{lem}\label{entropy lower con2}
$H_{\omega^n}$ is lower semicontinuous.
\end{lem}
\begin{proof}
Let $f$ and $f_s$ be defined as in the Lemma \ref{entropy lower con}. First, we consider the special case that $f$ and $f_s$ have uniformly positive lower bound, i.e. $f\ge \delta >0$ and $f_s \ge \delta >0$ for some constant $\delta>0$. Denote $\mathcal{F}(t)= t \log t$. Then we can put $$F_s(t)=\mathcal{F}(tf_s +(1-t)f) =\mathcal{F}(at+b),$$ where $a=f_s-f$ and $b=f$, then 
\begin{equation*}
    F_s'(t) =a(\log u_t +1),
\end{equation*}
where $u_t =tf_s +(1-t)f$, and 
\begin{equation*}
    F_s'' (t) =\frac{a^2}{u_t} \ge \frac{a^2}{C},
\end{equation*}
for some constant $C$ such that $f_s$ and $f$ are smaller than $C$. Hence
\begin{equation}\label{entropy converge}
\begin{split}
      \int_X (f_s \log f_s - f \log f)\omega^n &=\int_X (\int_0^1 \int_0^t F_s''(\lambda) d\lambda\, dt +\int_0^1 F_s'(0)dt) \omega^n \\
      & \ge \frac{1}{2C}\int_X (f_s -f)^2 \omega^n +\int_X F_s'(0) \omega^n \\
      &\ge  \int_X(f_s -f) (\log f +1)\omega^n.
\end{split}
\end{equation}
Here, we use the lower bound
\[\int_0^1 \int_0^t F_s''(\lambda) d\lambda\, dt =\int_0^1 \int_\lambda^1 F_s''(\lambda) dt \, d\lambda
=\int_0^1 (1-\lambda) F_s''(\lambda) d\lambda
\geq \int_0^1 (1-\lambda)  d\lambda \cdot \frac{a^2}{C}=\frac{a^2}{2C}.\]

Then we can use that $f_s$ converges weakly to $f$ in $L^1(\omega)$ to get that 
\begin{equation*}
    \lim_{s\rightarrow 0} \int_X(f_s -f) (\log f +1)\omega^n =0
\end{equation*}
So we have proved the Lemma for such $f$ and $f_s$. 

For general nonnegative $f$, without loss of generality we can assume that $$|f|+|f_s|\le \frac{e^{-1}}{2}$$ by dividing them by a constant. Then we have that $f_s+\delta \le e^{-1}$ for small $\delta$. We also have that $$f_s \log f_s \ge (f_s+\delta)\log (f_s +\delta),$$ since $\mathcal{F}(t)$ is decreasing on $(0,e^{-1})$. It is easy to show that $$\int f \log f - (f+\delta) \log (f+\delta)=o(1)$$ with respect to $\delta.$ 

So we see that:
\begin{equation*}
\begin{split}
        \int_X (f_s \log f_s -f \log f)\omega^n &\ge \int_X [(f_s +\delta)  \log (f_s+\delta) -(f+\delta) \log (f+\delta)] \omega^n +o(1) \\
       & \ge \int_X (f_s +\delta - f-\delta) ( \log (f+\delta)+1) \omega^n +o(1)
\end{split}
\end{equation*}
Here in the second line, we use (\ref{entropy converge}) with $f$ (resp. $f_s$) replaced by $f+\delta$ (resp. $f_s+\delta$).

Using the Lemma \ref{entropy lower con}, we have that
\begin{equation*}
    \lim_{s\rightarrow 0} \int_X (f_s -f) (\log (f+\delta) +1) \omega^n =0.
\end{equation*}
Therefore, we obtain that $$\liminf_{s\rightarrow 0}  \int_X (f_s \log f_s -f \log f)\omega^n \ge 0.$$
\end{proof}

Combining the Proposition \ref{H upp semicotinuous}, Lemma \ref{lem cont e},   Lemma \ref{entropy lower con2}, we can get:
\begin{prop}\label{continuity}
    The $K$-energy is continuous. 
\end{prop}

Combining the above results, we have shown that:
\begin{proof}
   of the Theorem \ref{thm convex}. This theorem follows from the Theorem \ref{subharmonicity} and the Proposition \ref{continuity}.
\end{proof}

\section{Solvability of the Lichnerowicz operator}\label{section: Solvability of the Lichnerowicz operator}
The main Proposition in this section is:
\begin{prop}\label{solve l operator}
Suppose that $D$ is a divisor with simple normal crossings. Suppose that $\omega$ is a Poincar\'e type cscK metric. Suppose that $f\in W_0^{0,2}$ such that $\int_{M\setminus D}f u \omega^n=0$ for any $u\in \overline{\mathbf{h}^D_{\parallelsum,\mathbb{R}}}$ . Then we can find a function $v\in W_0^{4,2}$ such that $Lv =f$, where $L$ is the Lichnerowicz operator of $\omega$.
\end{prop}

Then we can prove the Proposition \ref{solve l operator}:
\begin{proof}
of the Proposition \ref{solve l operator}. Using the Proposition \ref{prop e 5.1}, we have that (\ref{estimate for closed range}) holds for $\delta=0$. Then, we can use the Lemma \ref{lem 4.3} to get that $dim N(L|_{W_{0}^{m,2}})< \infty$ and  $R(L|_{W_{0}^{m,2}})$ are closed. Then, we can use the Lemma \ref{lem 4.2} to get that for any 
$f\in W_0^{0,2}$ such that $\int_{M\setminus D}f v \omega^n=0$ for any $v\in N(L|_{W_{0}^{4,2}})$, there exists $u \in W^{4,2}_0(M)$ such that $Lu=f$. Using the Lemma \ref{ker a}, we have that $$N(L|_{W_{0}^{4,2}})=\overline{\mathbf{h}^D_{\parallelsum,\mathbb{R}}}.$$ This concludes the proof of this Proposition.
\end{proof}

\subsection{Kernel and range}
To begin with, we need the following lemma which characterises the image of operators with closed range (See \cite[Theorem 2.19]{B}).
\begin{lem}\label{lem 4.2}
Let $A: D(A) \subset E \rightarrow F$ be an unbounded linear operator that is densely defined and closed. The following properties are equivalent:
\begin{enumerate}
\item  $R(A)$ is closed,
\item $R(A^*)$ is closed,
\item $R(A)=N(A^*)^{\bot}$,
\item $R(A^*)=N(A)^{\bot}$.
\end{enumerate}
\end{lem}
In the above lemma, we denote the range of $A$ as $R(A)$ and the kernel of $A$ as $N(A)$.
In order to use the above lemma, we set $$E=F=L_0^2\text{ and }A=L \text{ and } D(A)=W_0^{4,2}.$$ Denote $A^*$ as the adjoint of $A$.  The Lichnerowicz operator $L$ of a cscK metric is a real operator and satisfies: For any $v,u\in W_0^{4,2}$,
 $$\int_M v \cdot L u  \cdot \omega^n=\int_M L v  \cdot  u  \cdot \omega^n.$$ Then with the Lemma \ref{lem 4.12}, we have that $$D(A^*)=W_{0}^{4,2} \text{ and }A^*=L|_{W_{0}^{4,2}}.$$

Now we want to show that $R(L)$ is closed. We need to use the lemma below:
\begin{lem}\label{lem 4.3}
Suppose that $L$ satisfies the following formula for any $v \in W_{\delta}^{m,2}$ :
\begin{equation}\label{estimate for closed range}
||v||_{W_{\delta}^{m,2}(M \setminus D)}\le C(||L v||_{W_{\delta}^{m-4,2}(M\setminus D)}+||v||_{L^2(K)})
\end{equation}
for some compact set $K \subset \subset M $ and some constant $\delta$ and $m \ge 4$. Then we have that $dim N(L|_{W_{\delta}^{m,2}})< \infty$ and  $R(L|_{W_{\delta}^{m,2}})$ is closed.
\end{lem}
\begin{proof}
For simplicity, we denote $L|_{W_{\delta}^{m,2}}$ as  $L$ in this proof.  Firstly, we want to show that $dim N(L)< \infty$ by contradiction. Suppose that $$dim N(L)=\infty$$ Then we can find a sequence of functions $v_k \in N(L)$ such that $||v_k||_{W_{\delta}^{m,2}}=1$ and for any $k\neq l$,
\begin{equation*}
<v_k,v_l>_{W_{\delta}^{m,2}}=0.
\end{equation*}
Then we can extract a subsequence of $v_k$ (still denoted as $v_k$) converging to some function $v_{\infty}$ in $L^2_{loc}(K)$ sense. Then we can use (\ref{estimate for closed range}) to get that $v_k$ converge in the $W_{\delta}^{m,2}$ sense. This contradicts the fact that $||v_k||_{W_{\delta}^{m,2}}=1$ and for any $k\neq l$,
\begin{equation*}
<v_k,v_l>_{W_{\delta}^{m,2}}=0.
\end{equation*}
This concludes the proof for the first part of the proposition.

Secondly, we want to show that $R(L)$ is closed. We can replace $W_{\delta}^{m,2}$ by  $W_{\delta}^{m,2}/ N(L)$. Here $W_{\delta}^{m,2}/ N(L)$ is understood as the maximal subspace of $W_{\delta}^{m,2}$ which is perpendicular to $N(L)$. This is well-defined since $N(L)$ is finite-dimensional. Suppose that we have a sequence of functions $v_k \in W_{\delta}^{m,2}/ N(L)$ such that $f_k \triangleq L v_k$ converge to some function $f$ in $W_{\delta}^{m-4,2}$. We want to show that there exists $v \in W_{\delta}^{m,2}/ N(L)$ such that $$L(v)=f.$$

 First, we consider the following case: Suppose that we have that $||v_k||_{W_{\delta}^{m,2}}$ are uniformly bounded. Then we can get that $v_k$ (up to a subsequence) converge in $L_{loc}^2(M)$ to a function $v$. We can use (\ref{estimate for closed range}) to show that $v_k$ actually converges to $v$ in $W_{\delta}^{m,2}$. This implies that $Lv=f$. We are done.

Then we consider the second case: Suppose that $||v_k||_{W_{\delta}^{m,2}}$ is not bounded. We can assume that (by taking a subsequence) $$||v _k||_{W_{\delta}^{m,2}}\rightarrow +\infty.$$ We also define
\begin{equation*}
u_k =\frac{v_k}{||v_k||_{W_{\delta}^{m,2}}}
\end{equation*}
Since $||u_k||_{W_{\delta}^{m,2}}=1$ are uniformly bounded, we can take a subsequence (still denoted as $u_k$) such that $u_k$ converge to some function $u_{\infty}$ in $L^2(K)$. Note that $$||L u_k||_{W_{\delta}^{m-4,2}} \rightarrow 0.$$ 

Next we use (\ref{estimate for closed range}) again to get that $u_k$ converge to $u_{\infty}$ in $W_{\delta}^{m,2}/N(L)$ and $$L u_{\infty}=0.$$ Since $$(W_{\delta}^{m,2}/N(L))\cap N(L)=0,$$ we have that $u_{\infty}=0$. 

This contradicts the fact that  $u_k$ converges to $u_{\infty}$ in $W_{\delta}^{m,2}$ and $||u_k||_{W_{\delta}^{m,2}}=1$.
\end{proof}

\subsection{When does (\ref{estimate for closed range}) hold}

In this section, we want to prove the following Proposition:
\begin{prop}\label{prop e 5.1}
(\ref{estimate for closed range}) holds for $\delta=0$.
\end{prop}
Note that the above Proposition was proved by Auvray in \cite{A3} in the case that $D$ is smooth. Since we are dealing with the case that $D$ has simple normal crossings, we use a totally different approach.

Before we prove the Proposition \ref{prop e 5.1}, we first prove the following Lemma:

\begin{lem}\label{lem:hardy}
Let $0<r_0<1$, write $s=-\log\lvert z\rvert$, and let
$q\in C_c^\infty(\{0<\lvert z\rvert<r_0\};\bC)$. For $p\in\{0,1\}$,
\begin{equation}\label{eq:hardy}
\int_{0<\lvert z\rvert<r_0}
\frac{\lvert q\rvert^2}{\lvert z\rvert^2s^{2p+2}}\,dA
\le\frac{16}{4p+1}
\int_{0<\lvert z\rvert<r_0}
\frac{\lvert \partial_{\bar z}q\rvert^2}{s^{2p}}\,dA,
\end{equation}
where $dA$ is Euclidean area measure.
\end{lem}

\begin{proof}
Put $s_0=-\log r_0$ and expand
\[
q(s,\theta)=\sum_{m\in\mathbb Z}a_m(s)e^{im\theta}.
\]
For each Fourier mode,
\[
\partial_{\bar z}\bigl(a_m(s)e^{im\theta}\bigr)
=-\frac{e^{i(m+1)\theta}}{2\lvert z\rvert}
\bigl(a_m'(s)+ma_m(s)\bigr).
\]
For $a\in C_c^\infty((s_0,\infty);\bC)$, integration by parts gives
\begin{align*}
\int_{s_0}^{\infty}\lvert a'+ma\rvert^2s^{-2p}\,ds
&=\int_{s_0}^{\infty}\lvert a'\rvert^2s^{-2p}\,ds\\
&\quad+\int_{s_0}^{\infty}
\bigl(m^2s^2+2pms\bigr)\lvert a\rvert^2s^{-2p-2}\,ds.
\end{align*}
The weighted Hardy inequality
\begin{equation}\label{eq:radial-hardy}
\int_{s_0}^{\infty}\lvert a'\rvert^2s^{-2p}\,ds
\ge\frac{(2p+1)^2}{4}
\int_{s_0}^{\infty}\lvert a\rvert^2s^{-2p-2}\,ds
\end{equation}
follows from
\[
(2p+1)\int_{s_0}^{\infty}\lvert a\rvert^2s^{-2p-2}\,ds
=2\operatorname{Re}\int_{s_0}^{\infty}a'\bar a\,s^{-2p-1}\,ds
\]
and Cauchy--Schwarz. Since
\[
m^2s^2+2pms=(ms+p)^2-p^2\ge-p^2,
\]
we obtain
\[
\int_{s_0}^{\infty}\lvert a'+ma\rvert^2s^{-2p}\,ds
\ge\left(p+\frac14\right)
\int_{s_0}^{\infty}\lvert a\rvert^2s^{-2p-2}\,ds.
\]
Summing over the Fourier modes and using the factor $1/4$ in the
formula for $\lvert \partial_{\bar z}q\rvert^2$ proves
\eqref{eq:hardy}.
\end{proof}

In a cusp coordinate chart $\big( U, (z) \big)$, where $U=\{z: |z_i|\le \epsilon \text{ for any }i\}$ and $D\cap U= U \cap ( \cup_{i=1}^k \{z_i=0\})$, we use the holomorphic frame
\[
e_i=z_i\frac{\partial}{\partial z_i}\quad(1\le i\le k),
\qquad
e_i=\frac{\partial}{\partial z_i}\quad(k<i\le n).
\]
For the model metric \eqref{standard cusp}, this frame is orthogonal and
\[
\lvert e_i\rvert_{\om_0}^2 \asymp  (-\log |z_i|^2)^{-2}\quad(i\le k),
\qquad \lvert e_i\rvert_{\om_0}^2 \asymp 1 \quad(i>k).
\]
Write $V=\sum_i a_i e_i$. Since the frame is holomorphic,
\[
\bar\partial V=\sum_i\bar\partial a_i\otimes e_i.
\]
For a fixed cusp direction $j\le k$, set
\[
\bar\partial_j V=d\bar z_j\otimes
\sum_i\frac{\partial a_i}{\partial\bar z_j}e_i.
\]
We compute
$$|V|_{\omega_0}^2 = \sum_{i=1}^n |a_i|^2|e_i|_{\omega_0}^2, \quad |\bar\partial_jV|_{\omega_0}^2 = (\omega_0)^{j\bar j} \sum_{i=1}^n |\partial_{\bar z_j}a_i|^2|e_i|_{\omega_0}^2. $$

\begin{lem}\label{lem:local-vector}
For every smooth vector field $V$ compactly supported in $U \setminus D$, we have that:
\begin{equation}\label{eq:local-vector}
\lVert V\rVert_{L^2(U)}^2
\le C\lVert \bar\partial_jV \rVert_{L^2(U)}^2,\quad \forall 1\le j \le k.
\end{equation}
\end{lem}
\begin{proof}

For a fixed normal direction \(j\le k\), apply Lemma 5.5 separately to the coefficients:
\begin{equation*}
\begin{split}
\lVert V\rVert_{L^2(U)}^2 & = \int_U \sum_{i=1}^k |a_i e_i|^2  \omega_0^n +  \int_U \sum_{i=k+1}^n |a_i e_i|^2  \omega_0^n \\
& \le C \int_U \sum_{i=1}^k |a_i|^2 (-\log |z_i|^2)^{-2}  \omega_0^n + C \int_U \sum_{i=k+1}^n |a_i |^2  \omega_0^n \\
& \le C \int_U \sum_{i=1}^k |\bar \partial_j a_i |^2 |z_j|^2 (-\log |z_j|^2)^2 (-\log |z_i|^2)^{-2} \omega_0^n \\
&+ C \int_U \sum_{i=k+1}^n |\bar \partial_j a_i |^2  |z_j|^2 (-\log |z_j|^2)^2 \omega_0^n \\
& \le C \int_U \sum_{i=1}^n |\bar \partial_j (a_i e_i) d \bar z^j |^2 \omega_0^n
\end{split}
\end{equation*}
In the third line above, we use the Lemma \ref{lem:hardy} in the $z_j$ direction.
\end{proof}

\begin{prop}\label{prop:vector}
There exist a compact set $K\Subset M$ and $C>0$ such that
\begin{equation}\label{eq:vector}
\|V\|_{L^2(M)}^2\le C\left(
\|\bar \partial  V\|_{L^2(M)}^2+\|V\|_{L^2(K)}^2\right)
\end{equation}
for every $V\in C_c^\infty(M,T^{1,0}M)$.
\end{prop}

\begin{proof}
Choose a finite family of  polydiscs with cusp coordinates and interior coordinate
charts, together with real smooth functions $\chi_\alpha$ on $X$
subordinate to this cover, such that
\[
\sum_\alpha\chi_\alpha^2=1,
\qquad \supp\chi_\alpha\Subset U_\alpha.
\]
Let $I_\alpha$ be the set of components represented by coordinate
hyperplanes in $U_\alpha$. The cover can be chosen so that, for $j$
not in $I_\alpha$, $\mathbf{u}_j$ is bounded on $\supp\chi_\alpha$.

Fix $\zeta\in C^\infty(\R,[0,1])$ with $\zeta=0$ on $(-\infty,1]$
and $\zeta=1$ on $[2,\infty)$, and put
\[
\zeta_{j,R}=\zeta(\mathbf{u}_j/R).
\]
Then
\begin{equation}\label{eq:zeta-gradient}
|d\zeta_{j,R}|_g\le C/R.
\end{equation}
For $R$ sufficiently large, every chart contributing to
$\chi_\alpha\zeta_{j,R}$ has $j\in I_\alpha$. 

We define
$ O_R=\bigcup_{j=1}^N\{\mathbf u_j>R\} = \{\max_j\mathbf u_j>R\}.  $
Then
$$ M\setminus O_{2R} = \{\mathbf u_j\le2R\text{ for every }j\} \Subset M. $$

 For any $z\in O_{2R}$, there exists $i$ such that $\mathbf{u}_i (z) \ge 2R$. This implies that $$\zeta_{i,R}(z)=1.$$ Then, we can get that:
\begin{equation*}
\begin{split}
\|V\|^2_{L^2(O_{2R})} &\le C \sum_{\alpha} \sum_{j \in I_{\alpha}} \int_{O_{R}\cap U_{\alpha}} |\chi_\alpha^2 \zeta_{j,R} V|^2 \omega_0^n \le C \sum_{\alpha} \sum_{j \in I_{\alpha}} \int_{O_{R}\cap U_{\alpha}} | \bar \partial_j (\chi_\alpha^2 \zeta_{j,R} V)d \bar z^j|^2 \omega_0^n \\
& \le C \sum_{\alpha} \sum_{j \in I_{\alpha}} \int_{O_{R}\cap U_{\alpha}} \Big( | \bar \partial_j (\chi_\alpha^2) \zeta_{j,R} V d \bar z^j|^2 + |\chi_{\alpha}^2  \bar \partial_j  \zeta_{j,R}   V d \bar z^j|^2 + |\chi_{\alpha}^2 \zeta_{j,R} \bar \partial_j V d \bar z^j|^2 \Big) \omega_0^n \\
& \le C \sum_{\alpha} \sum_{j \in I_{\alpha}} \int_{O_{R}\cap U_{\alpha}} \Big( | \zeta_{j,R} V |^2 |z_j|^2 (\log^2 |z_j|^2) + |\chi_{\alpha}^2   V|^2\frac{1}{R^2} + |\chi_{\alpha}^2 \zeta_{j,R} \bar \partial_j V d \bar z^j|^2 \Big) \omega_0^n \\
& \le \frac{C}{R^2}  \int_{O_{R}} |V|^2 \omega_0^n + C \int_{O_{R}} |\bar \partial V|^2 \omega_0^n.
\end{split}
\end{equation*}
In the first line of the above formula, we use the Lemma \ref{lem:local-vector}. In the third line, we use the fact that $\chi_{\alpha}$ is smooth on $X$ and the estimate (\ref{eq:zeta-gradient}).
Thus, we can get that:
\begin{equation*}
\|V\|^2_{L^2(M)}= \|V\|^2_{L^2(O_{2R})} + \|V\|^2_{L^2(M\setminus O_{2R})}  \le \frac{C}{R^2}  \int_{O_{R}} |V|^2 \omega_0^n + C \int_{O_{R}} |\bar \partial V|^2 \omega_0^n + \|V\|^2_{L^2(M\setminus O_{2R})}.
\end{equation*}
Then, we can let $R$ be big enough such that $\frac{C}{R^2}\le \frac{1}{2}$ to get that:
 \begin{equation*}
\|V\|^2_{L^2(M)} \le  C \int_{M} |\bar \partial V|^2 \omega_0^n +C \|V\|^2_{L^2(M\setminus O_{2R})}.
\end{equation*}
This concludes the proof of this Proposition, with the compact set $K$ to be $M \setminus O_{2R}$.
\end{proof}

Similarly, we can get that:
\begin{lem}\label{lem u com est}
There exists a compact set $K \Subset M$ suchthat for any $u \in C_c^{\infty}(M)$, we have that:
\begin{equation*}
\|u\|_{L^2(M)}^2 \le C (\|du\|_{L^2(M)}^2 + \|u\|_{L^2(K)}^2).
\end{equation*}
\end{lem}
\begin{proof}
We use the same notations as in the proof of the Proposition \ref{prop:vector}. Then, we have that:
\begin{equation*}
\begin{split}
&\|u\|_{L^2(O_{2R})}^2 \\
& \le C \sum_{\alpha} \sum_{j \in I_{\alpha}} \int_{O_{R}\cap U_{\alpha}}|\chi_{\alpha} u \zeta_{j,R} |^2 \omega_0^n \le C \sum_{\alpha} \sum_{j \in I_{\alpha}} \int_{O_{R}\cap U_{\alpha}}| \bar \partial_j (\chi_{\alpha} u \zeta_{j,R} ) |^2 |z_j|^2 \log^2 |z_j|^2  \omega_0^n\\
& \le C \sum_{\alpha} \sum_{j \in I_{\alpha}} \int_{O_{R}\cap U_{\alpha}} \Big( | \bar \partial_j (\chi_{\alpha}) u \zeta_{j,R}  |^2+ | \bar \partial_j (u) \chi_{\alpha} \zeta_{j,R}  |^2 + | \bar \partial_j (\zeta_{j,R}) \chi_{\alpha}  u |^2 \Big) |z_j|^2 \log^2 |z_j|^2  \omega_0^n \\
& \le  C \sum_{\alpha} \sum_{j \in I_{\alpha}}  \int_{O_{R}\cap U_{\alpha}} (\frac{1}{R^2}|u|^2 + |\bar \partial u|^2)\omega_0^n \\
& \le \frac{C}{R^2} \int_{O_{R}} |u|^2 \omega_0^n + C \int_{O_{R}} |\bar \partial u|^2 \omega_0^n.
\end{split}
\end{equation*}
In the first line of the above formula, we use the Lemma \ref{lem:hardy}. In the second line, we use the fact that $\chi_{\alpha}$ is  smooth on $X$ and (\ref{eq:zeta-gradient}).
Then, we can get that:
\begin{equation*}
\|u\|_{L^2(M)}^2  = \|u\|_{L^2(O_{2R})}^2 + \|u\|_{L^2(M\setminus O_{2R})}^2 \le  \frac{C}{R^2} \int_{O_{R}} |u|^2 \omega_0^n + C \int_{O_{R}} |\bar \partial u|^2 \omega_0^n + \|u\|_{L^2(M\setminus O_{2R})}^2.
\end{equation*}
Then, we can let $R$ be big such that $\frac{C}{R^2}\le \frac{1}{2}$ to get:
\begin{equation*}
\|u\|_{L^2(M)}^2  \le C \int_M |\bar \partial u|^2 \omega_0^n + C \|u\|_{L^2(M \setminus O_{2R})}^2. 
\end{equation*}
\end{proof}

\begin{lem}\label{lem 4.12}
Suppose that $u\in W_{\delta}^{0,2}$ and $L u \in W_{\delta}^{k,2}$ for some $\delta$. Then we have that $u\in W_{\delta}^{k+4,2}$ and
\begin{equation*}
||u||_{W^{k+4,2}_{\delta}}\le C(||Lu||_{W^{k,2}_{\delta}} +||u||_{W^{0,2}_{\delta}}).
\end{equation*}
\end{lem}
\begin{proof}
This lemma can be proved in the same way as \cite[Lemma 1.12]{A}. We just sketch the proof here. We can use a covering of $M$ using the quasi-coordinate mentioned in Section 3. In each coordinate, the Poincar\'e type K\"ahler metric is quasi-isometric to the Euclidean metric. The operator $Lu$ is understood distributionally. Then we can use the standard $L^p$ estimate for $L$ in each quasi-coordinate. Then we patch them together to prove the lemma. 
\end{proof}

Finally, we can prove the Proposition \ref{prop e 5.1}:

\begin{proof}
of the Proposition \ref{prop e 5.1}. First, we assume that $u\in C_c^{\infty}(M)$. Define $V= \nabla_{\omega}^{1,0}u$. Applying the Proposition \ref{prop:vector} to $V$, we can get that:
\begin{equation*}
||du||^2_{L^2(M)}= 2 ||V||^2_{L^2(M)} \le C (||\mathcal{D} u||^2_{L^2(M)} + ||du||^2_{L^2(K)}).
\end{equation*}
Combine this formula with Lemma \ref{lem u com est} and denote $K_1$ as the compact set which is the union of $K$ and the compact set appearing in Lemma \ref{lem u com est}. Then, we can get that:
\begin{equation}\label{e u Du}
||u||^2_{L^2(M)}\le C(||\mathcal D u||^2_{L^2(M)} + ||u||^2_{W^{1,2}(K_1)}).
\end{equation}
Recall that $L= \mathcal{D}^* \mathcal{D}$. Then, we can get that:
\begin{equation*}
||\mathcal{D}u||_{L^2(M)}^2= \int_M u L u \omega^n \le ||u||_{L^2(M)} ||Lu||_{L^2(M)}\le \epsilon ||u||^2_{L^2(M)} + C_{\epsilon} ||Lu ||^2_{L^2(M)},
\end{equation*}
for any $\epsilon>0$ and a constant $C_{\epsilon}$ depending on $\epsilon$.
Combining this formula with (\ref{e u Du}) and let $\epsilon$ be small enough, we can get that:
\begin{equation*}
||u||_{L^2(M)} \le C (||Lu||_{L^2(M)} + ||u||_{W^{1,2}(K_1)}).
\end{equation*}
Plugging the above formula into the formula given by the Lemma \ref{lem 4.12} with $\delta=0$, we can get that:
\begin{equation}\label{e uw42 com}
||u||_{W_0^{4,2}(M)}\le C(||Lu||_{L^2(M)} + ||u||_{W^{1,2}(K_1)}).
\end{equation}
Choose
$K_1\Subset\operatorname{int}K_2\Subset M. $
Then the $L^p$ interior estimate gives
\begin{equation}\label{e uw12 com}
 ||u||_{W^{1,2}(K_1)} \le C(||u||_{L^2(K_2)}+ ||Lu ||_{L^2(K_2)}).
\end{equation}
Combining (\ref{e uw42 com}) with (\ref{e uw12 com}), we obtain \eqref{estimate for closed range} for smooth $u$.

For every $m\geq4$, add the bootstrap
$$\|u\|_{W^{m,2}}
\leq C(\|Lu\|_{W^{m-4,2}}+\|u\|_{L^{2}(K)}),$$
using Lemma \ref{lem 4.12} and the established $m=4$ estimate.

For \(u\in W_0^{m,2}\),
we use the compact-support smoothing and the estimates \eqref{eq:constructive-cutoff-four} of the cutoff function
to see $ C_c^\infty(M)$ is dense in $W_0^{4,2}$
and
$ L:W_0^{4,2}\rightarrow L_0^2 $
is continuous. Thus the estimate passes to the limit $u$
and we can conclude the proof of this Proposition.
\end{proof}

\subsection{Kernel and holomorphic vector fields}
For a real-valued function $u$, set
\[
\mathcal{D}_\omega u\coloneqq\bar\partial(\nabla_\omega^{1,0}u),
\qquad \mathcal{D}\coloneqq\mathcal{D}_\omega.
\]
We use the standard Hermitian tensor norm, so that $L=\mathcal{D}^*\mathcal{D}$ and $\mathcal{D}u=0$ precisely when $\nabla_\omega^{1,0}u$ is holomorphic.
 Recall that the  $$\overline{\mathbf{h}^D_{\parallelsum, \mathbb{R}}}=\{f\in C_{\mathbb{R}}^{\infty}(X \setminus D): \nabla^{1,0}f \in \mathbf{h}_{\parallelsum}^D\}.$$  We record the following Lemma:

\begin{lem}\label{ker a}
Suppose that $\omega$ is a Poincar\'e type cscK metric, then
\begin{equation*}
N(L|_{W_{0}^{k,2}})=\overline{\mathbf{h}^D_{\parallelsum,\mathbb{R}}},\quad k\geq 4.
\end{equation*}
\end{lem}
\begin{proof}
The formula $N(L |_{W_{0}^{k,2}}) \subset \overline{\mathbf{h}^D_{\parallelsum,\mathbb{R}}}$ can be shown by using the Laurent expansion of holomorphic functions near the divisor.  Indeed, for any $u\in N(L|_{W_{0}^{k,2}})$, we have that 
\begin{equation*}
0=\int_M L u \cdot u\cdot \omega^n = \int_M \mathcal{D}^* \mathcal{D}u \cdot u \omega^n = \int_M |\mathcal{D}u|^2 \omega^n .
\end{equation*}
This implies that $\mathcal{D}u=0$ which means that $v \triangleq \nabla^{1,0}u$ is a holomorphic vector field on $M$. We will prove that $v$ can be extended holomorphically across $D$.

Since $u\in W_0^{k,2}$, we can get that $|v| \in L^2(\omega^n)$. In an arbitrary cusp coordinate chart $(U,(z))$ where $D=\{(z): z_1 \cdots z_k =0\}$, we denote $v=v^i \frac{\partial}{\partial z^i}$. Since $\omega$ is equivalent to the standard Poincar\'e type K\"ahler metric (\ref{standard cusp}), we have that:
\begin{equation}\label{e 4.5.1}
\int_U |v|_{\omega_0}^2 \omega_0^n \le C \int_M |v|_{\omega}^2 \omega^n < +\infty
\end{equation}
Then we have that:
\begin{equation}\label{e 4.5.2}
\begin{split}
&\int_U |v|_{\omega_0}^2 \omega_0^n \ge C \int_U |v^i|^2 (\omega_0)_{i \bar i} \omega_0^n\\
& =C \int_U |v^i|^2 \frac{1}{|z_i|^2 log^2 (|z_i|^2)}  \Pi_{l=1}^k \frac{1}{|z_l |^2 log^2 (|z_l|^2)} dVol_E \\
&= C\int_U |v^i|^2 \frac{n!}{|z_i|^4 log^4(|z_i|^2)}   \Pi_{l=1}^{i-1} \frac{1}{|z_l |^2 log^2 (|z_l|^2)}   \Pi_{l=i+1}^{k} \frac{1}{|z_l |^2 log^2 (|z_l|^2)}  dVol_E.
\end{split}
\end{equation}
Here $dVol_E=(\sqrt{-1}dz^1 \wedge d \bar z^1) \wedge \cdots \wedge (\sqrt{-1}dz^n \wedge d\bar z ^n)$. 

Note that we have a Laurent series of $v^i$(see \cite[Proposition 1.4]{R})
 \begin{equation*}
 v^i =\sum_{\mu \in \mathbb{Z}^{k-1}} \sum_{\nu \in \mathbb{N}^{n-k}} \sum_{j\in \mathbb{Z}}C_{\mu \nu j} z'^{\mu} w^{\nu} z_i^j.
 \end{equation*}
Here $z'=(z_1,\ldots,z_{i-1}, z_{i+1}, \ldots, z_k )$, $w=(z_{k+1},\ldots, z_n)$.

 Let $\epsilon>0$ be a constant such that $$U_{\epsilon} \triangleq \{z: |z_i|\le \epsilon \text{ for any }i\} \subset U.$$ We also denote $$U'_{\epsilon} \triangleq \{(z_1,\ldots, z_{i-1}, z_{i+1},\ldots z_k ): |z_l|\le \epsilon \text{ for any }l\neq i \text{ and } l \le k\},$$
 and
 $$U_{\epsilon}''=\{(z_{k+1},\ldots, z_n): |z_i|\le \epsilon, \text{ for any }i \ge k+1\}$$
   Then we have that:
\begin{equation}\label{e 4.5.3}
\begin{split}
&\int_{U_{\epsilon}} |v^i|^2 \frac{n!}{|z_i|^4 log^4(|z_i|^2)}   \Pi_{l=1}^{i-1} \frac{1}{|z_l |^2 log^2 (|z_l|^2)}   \Pi_{l=i+1}^{k} \frac{1}{|z_l |^2 log^2 (|z_l|^2)}  dVol_E \\
=& \sum_{\mu \in \mathbb{Z}^{k-1}} \sum_{\nu \in \mathbb{N}^{n-k}} \sum_{j\in \mathbb{Z}} |C_{\mu \nu j}|^2 (\int_{U'_{\epsilon}} |z'|^{2\mu} \Pi_{l=1}^{i-1} \frac{1}{|z_l |^2 log^2 (|z_l|^2)}   \Pi_{l=i+1}^{k} \frac{1}{|z_l |^2 log^2 (|z_l|^2)}  dVol_E(z')) \cdot \\
& (\int_{U_{\epsilon}''} |w|^{2\nu} dVol_E (w) ) \cdot  ( \int_{B_{\epsilon}(0)} \frac{|z_i|^{2j-4}}{log^4(|z_i|^2)}  \sqrt{-1} dz^i \wedge d\bar z^i).
\end{split}
\end{equation}

Combining (\ref{e 4.5.1}), (\ref{e 4.5.2}) and (\ref{e 4.5.3}), we have that $C_{\mu \nu j} \neq 0$ implies that $\mu \in \mathbb{N}^{k-1}$ and $j \ge 1$. This proves that  $v^i$ can be extended holomorphically to $D$ and vanishes on $\{z_i=0\}$.  

For \(i>k\), the integrability condition \eqref{e 4.5.2} is
$$ \int_U |v^i|^2 \prod_{\ell=1}^k \frac{1}{|z_\ell|^2\log^2|z_\ell|^2}\,dVol_{\mathrm E}<\infty. $$
A Laurent monomial with a negative exponent in any normal variable is not integrable. Thus \(v^i\) extends holomorphically across every component.

This concludes the proof of $$N(L |_{W_{0}^{k,2}}) \subset \overline{\mathbf{h}^D_{\parallelsum,\mathbb{R}}}.$$

The formula $\overline{\mathbf{h}^D_{\parallelsum,\mathbb{R}}} \subset N(L|_{W_{0}^{k,2}}) $ can be shown as follows:
For any $f\in \overline{\mathbf{h}^D_{\parallelsum,\mathbb{R}}}$, we have that $v\triangleq \nabla^{1,0}_{\omega}f$
 is a holomorphic vector field on $X$.  First, we claim that:
 \begin{equation}\label{nabla diff metric 1}
 v =\nabla_{\omega_X}^{1,0}[f-v(u)],
 \end{equation}
 where $u$ is the potential such that $\omega =\omega_X +dd^c u$. Indeed, we can calculate that:
 \begin{equation*}
v^i [(g_X)_{i\bar l}+u_{i\bar l}]= v^i g_{i\bar l}=g^{i\bar j}f_{\bar j}g_{i\bar l}=f_{\bar l}.
 \end{equation*}
 Multiply both sides of the above formula by $g_X^{k\bar l}$ and sum over $l$. We obtain
 \begin{equation*}
 g^{k\bar j} f_{\bar j} \frac{\partial}{\partial_z^k}+\nabla_{\omega_X}^{1,0} (v(u))=\nabla_{\omega_X}^{1,0}f.
 \end{equation*}
 Here we use the fact that $v$ is holomorphic. This concludes the proof of the claim.  
 
 By the definition of the Poincar\'e type metric, we have that the derivatives of $u$ of any order are bounded with respect to a given Poincar\'e type metric. As a result, we have that $$v(u)\in C_0^{k,\alpha}\text{ for any }k.$$ Note that $v(u)$ can be a complex-valued function. Here we mean that both the real part and the imaginary part of $v(u)$ belong to $C_0^{k,\alpha}$. 
 
 Set \(F=f-v(u)\). Equation (\ref{nabla diff metric 1}) gives \(F_{\bar j}=v^i(g_X)_{i\bar j}\). Near any point of \(D\), choose a smooth local K\"ahler potential \(\psi\) for \(\omega_X\). Since \(v\) is holomorphic, \(\bar\partial(v(\psi))=\bar\partial F\), so \(H=F-v(\psi)\) is holomorphic off \(D\). Moreover, \(f\) is real-valued and \(v(u)\) is bounded, hence \(\operatorname{Im}H=-\operatorname{Im}[v(u)+v(\psi)]\) is bounded. The bounded holomorphic functions \(e^{\sqrt{-1}H}\) and \(e^{-\sqrt{-1}H}\) extend across \(D\), and their extensions have product \(1\). A local holomorphic logarithm therefore shows that \(H\) extends holomorphically across \(D\). Consequently, \(F\) extends smoothly to \(X\). Thus, under quasi-coordinates, \(F\in C_0^{k,\alpha}\) and \(f=F+v(u)\in C_0^{k,\alpha}\). Since \(\omega\) has finite volume, \(f\in W_0^{k,2}\); the identity \(\mathcal D_\omega f=0\) now gives \(Lf=0\).
 
 
 This concludes the proof of the Lemma. 
 \end{proof}  
 
 \subsection{Approximation}
 We fix a Poincar\'e type cscK metric $\omega$ and use the real Hilbert space
$L^2_0(M,\omega)=W^{0,2}_0(M,\omega)$. 
Recall that in Proposition~\ref{solve l operator} and Lemma~\ref{ker a}, we already have
\begin{equation}\label{eq:constructive-inputs}
\begin{gathered}
R\bigl(L|_{W_0^{4,2}}\bigr)
=N\bigl(L|_{W_0^{4,2}}\bigr)^\perp,
\quad
N\bigl(L|_{W_0^{4,2}}\bigr)
=\overline{\mathbf{h}_{\mathbin{\!/\mkern-5mu/\!},\mathbb{R}}^D}.
\end{gathered}
\end{equation}

\begin{prop}\label{kernel approximation}
Let $f\in L_0^2(M,\omega)$ satisfy
\begin{equation}\label{eq:constructive-orthogonality}
\int_M f\cdot h\cdot\omega^n=0
\quad\text{for every }h\in
\overline{\mathbf{h}_{\mathbin{\!/\mkern-5mu/\!},\mathbb{R}}^D}.
\end{equation}
Then there is a uniquely normalised real solution $u\in W_0^{4,2}$ of $Lu=f$, with $u$ orthogonal to the kernel. There are real $v_m\in C_c^\infty(M)$ such that
\begin{equation}\label{eq:constructive-conclusion}
v_m\longrightarrow u\quad\text{in }W_0^{4,2},
\qquad
Lv_m\longrightarrow f\quad\text{in }L_0^2.
\end{equation}

In particular, if $f$ is smooth on $M$, the cutoff function $\eta_m$ from Lemma~\ref{snc-cutoff} is used to construct the approximations
\begin{equation}\label{eq:constructive-vm}
v_m=(1-\eta_m)u.
\end{equation}
Furthermore,
\begin{equation}\label{eq:constructive-tail-estimate}
\begin{split}
\lVert Lv_m-f\rVert_{L_0^2}
\le{}&\lVert f\rVert_{L^2(\{T>m\},\omega^n)}
+\frac{C}{m}\sum_{r=0}^3
\lVert\nabla^r u\rVert_{L^2(\{m\le T\le2m\},\omega^n)}.
\end{split}
\end{equation}
The constants depend on the fixed metric and the cutoff function, but not on $m$.
\end{prop}
\begin{proof}
\noindent\textbf{Step 1:}
By \eqref{eq:constructive-inputs} and \eqref{eq:constructive-orthogonality}, there exists a solution such that
\begin{equation}\label{eq:constructive-solution}
Lu=f,\qquad
u\in W_0^{4,2},\qquad
\int_Mu\cdot h\cdot\omega^n=0
,\quad \forall h\in N(L|_{W_0^{4,2}}).
\end{equation}
Since constants belong to the kernel, \eqref{eq:constructive-solution} also implies $\int_Mu\,\omega^n=0$. 

The operator $L:W_0^{4,2}\to L_0^2$ is bounded by 
\begin{equation}\label{eq:S5-L-continuity}
\lVert L u\rVert_2\le C\lVert u\rVert_{W_0^{4,2}}.
\end{equation} 

Due to Proposition \ref{prop e 5.1} with $m=4$,
\begin{equation}\label{eq:compact-remainder}
\|u\|_{W_0^{4,2}}
\le C_0\left(
\|f\|_{L_0^2}+\|u\|_{L^2(K,\omega^n)}
\right),
\end{equation}
where $K\Subset M$ is fixed. 
In fact, we further have
\begin{equation}\label{eq:constructive-inverse-bound}
\lVert u\rVert_{W_0^{4,2}}\le C\lVert f\rVert_{L_0^2}.
\end{equation}

That follows from the following claim
\begin{equation}\label{eq:coercive-estimate}
\|u\|_{W_0^{4,2}}\le C\|Lu\|_{L_0^2}
\quad\text{for every }
 u\in W_0^{4,2}\cap N\bigl(L|_{W_0^{4,2}}\bigr)^\perp.
\end{equation}
Suppose that \eqref{eq:coercive-estimate} is false. Then for every positive integer $j$ there is a nonzero function $u_j\in W_0^{4,2}$, orthogonal to $N(L|_{W_0^{4,2}})$, such that
\[
\|u_j\|_{W_0^{4,2}}>j\|Lu_j\|_{L_0^2}.
\]
Set the normalisation $v_j=u_j/\|u_j\|_{W_0^{4,2}}$. Then
\begin{equation}\label{eq:contradiction-sequence}
\|v_j\|_{W_0^{4,2}}=1,
\qquad
\|Lv_j\|_{L_0^2}<\frac1j,
\qquad
v_j\perp N\bigl(L|_{W_0^{4,2}}\bigr).
\end{equation}

Choose a relatively compact smooth domain $U$ with $K\subset U\Subset M$. The sequence $v_j$ is bounded in $W^{4,2}(U)$. Then by compactness in $L^2(K,\omega^n)$, after passing to a subsequence, we have
\begin{equation}\label{eq:local-cauchy}
\|v_j-v_k\|_{L^2(K,\omega^n)}\longrightarrow0
\qquad\text{as }j,k\longrightarrow\infty.
\end{equation}

Applying \eqref{eq:compact-remainder} to $v_j-v_k$ and using \eqref{eq:contradiction-sequence} and \eqref{eq:local-cauchy}, we have
\begin{align}
\|v_j-v_k\|_{W_0^{4,2}}
&\le C_0\left(
\|L(v_j-v_k)\|_{L_0^2}
+\|v_j-v_k\|_{L^2(K,\omega^n)}
\right)\rightarrow 0.
\label{eq:global-cauchy}
\end{align}
Hence there exists a limit $v_\infty\in W_0^{4,2}$ with
\begin{equation}\label{eq:strong-limit}
v_j\longrightarrow v_\infty\quad\text{in }W_0^{4,2},
\qquad
\|v_\infty\|_{W_0^{4,2}}=1.
\end{equation}

By \eqref{eq:S5-L-continuity} and \eqref{eq:strong-limit},
\[
Lv_j\longrightarrow Lv_\infty\quad\text{in }L_0^2.
\]
On the other hand, the second term in \eqref{eq:contradiction-sequence} gives $Lv_j\to0$, so $Lv_\infty=0$ and
\[
v_\infty\in N\bigl(L|_{W_0^{4,2}}\bigr).
\]
The third term in \eqref{eq:contradiction-sequence} implies that for every $h\in N(L|_{W_0^{4,2}})$,
\[
\int_Mv_\infty h\,\omega^n
=\lim_{j\to\infty}\int_Mv_jh\,\omega^n=0.
\]
In particular, we take $h=v_\infty$. Then $\|v_\infty\|_{L_0^2}^2=0$, hence $v_\infty=0$. This contradicts $\|v_\infty\|_{W_0^{4,2}}=1$ in \eqref{eq:strong-limit}. The contradiction proves \eqref{eq:coercive-estimate}. For a normalized solution of $Lu=f$, it gives \eqref{eq:constructive-inverse-bound}.


\noindent\textbf{Step 2:}
We assume that $f$ is smooth. Using \eqref{eq:constructive-vm} and $Lu=f$, we get
\begin{equation}\label{eq:constructive-residual-identity}
\begin{split}
Lv_m-f
&=L((1-\eta_m)u)-f
=(1-\eta_m)Lu+[L,1-\eta_m]u-f\\
&=-\eta_m f+[L,1-\eta_m]u.
\end{split}
\end{equation}
In quasi-coordinate estimates, we have the bound for the lower order term
\begin{equation}\label{eq:constructive-pointwise-commutator}
\bigl|[L,1-\eta_m]u\bigr|
\le C\sum_{r=1}^4|\nabla^r\eta_m|_\omega
\sum_{j=0}^{4-r}|\nabla^ju|_\omega.
\end{equation}
The order of the derivative of $u$ is at most three.

Inserting \eqref{eq:constructive-cutoff-four} and \eqref{eq:constructive-cutoff-support} into \eqref{eq:constructive-pointwise-commutator}, we obtain its bound
\begin{equation}\label{eq:constructive-commutator-L2}
\lVert[L,1-\eta_m]u\rVert_{L_0^2}
\le\frac{C}{m}\sum_{r=0}^3
\lVert\nabla^ru\rVert_{L^2(\{m\le T\le2m\},\omega^n)}
\le\frac{C}{m}\lVert u\rVert_{W_0^{3,2}}.
\end{equation}

For the first term in \eqref{eq:constructive-residual-identity}, 
since $0\le\eta_m\le1$ and $\eta_m=0$ on $\{T\le m\}$,
\[
\lVert\eta_mf\rVert_{L_0^2}
\le\lVert f\rVert_{L^2(\{T>m\},\omega^n)}\rightarrow0,
\]
which follows from $f\in L_0^2$ and properness of $T$. Inserting it and \eqref{eq:constructive-commutator-L2} into \eqref{eq:constructive-residual-identity}, we prove \eqref{eq:constructive-tail-estimate}:
\begin{equation}\label{eq:constructive-tail-estimate 1}
\begin{split}
&\lVert Lv_m-f\rVert_{L_0^2}
\le\lVert f\rVert_{L^2(\{T>m\},\omega^n)}
+\frac{C}{m}\sum_{r=0}^3
\lVert\nabla^r u\rVert_{L^2(\{m\le T\le2m\},\omega^n)}\\
&\le\lVert f\rVert_{L^2(\{T>m\},\omega^n)}
+\frac{C}{m}\lVert u\rVert_{W_0^{3,2}}
\le\lVert f\rVert_{L^2(\{T>m\},\omega^n)}
+\frac{C}{m}\lVert f\rVert_{L_0^2}.
\end{split}
\end{equation}
As a result, $Lv_m\to f$ in $L_0^2$.

Since $v_m-u=-\eta_mu$, \eqref{eq:constructive-cutoff-four} gives
\begin{equation}\label{eq:constructive-domain-convergence}
\begin{split}
\lVert v_m-u\rVert_{W_0^{4,2}}
\le{}& C\lVert u\rVert_{W_0^{4,2}(\{T>m\})}
+\frac{C}{m}\sum_{r=0}^3
\lVert\nabla^ru\rVert_{L^2(\{m\le T\le2m\},\omega^n)}
\rightarrow0.
\end{split}
\end{equation}
Hence by \eqref{eq:constructive-inverse-bound} the approximating sequence is bounded in $W_0^{4,2}$:
\begin{equation}\label{eq:constructive-uniform-Sobolev}
\sup_{m\ge1}\lVert v_m\rVert_{W_0^{4,2}}
\le C\lVert u\rVert_{W_0^{4,2}}
\le C\lVert f\rVert_{L_0^2}.
\end{equation}

\noindent\textbf{Step 3:}
For $f \in L_0^2$, $v_m=(1-\eta_m)u$ is a compactly supported $W^{4,2}$ function. For each fixed $m$, choose a relatively compact neighborhood of its support inside $\{T<2m+1\}$. We apply the finite partition of unity in that neighborhood and the  local mollification to obtain smooth approximations of $v_m$, denoted $\tilde v_m$, such that
\begin{equation}\label{eq:constructive-smoothing}
\lVert \tilde v_m-v_m\rVert_{W_0^{4,2}}
\le \frac{1}{m(1+C)},
\end{equation}
where $C$ is the bound in \eqref{eq:S5-L-continuity}. 
Then adding \eqref{eq:constructive-tail-estimate 1} and \eqref{eq:constructive-smoothing}, we get
\begin{equation}\label{eq:constructive-nonsmooth-data}
\begin{split}
\lVert L\tilde v_m-f\rVert_{L_0^2}
\le{}&\lVert f\rVert_{L^2(\{T>m\},\omega^n)}
+\frac{C}{m}\lVert u\rVert_{W_0^{3,2}}+\frac1m
\rightarrow0.
\end{split}
\end{equation}
Together with \eqref{eq:constructive-domain-convergence} and \eqref{eq:constructive-smoothing}, this proves \eqref{eq:constructive-conclusion}. 
\end{proof}

\begin{rem}
In particular, on each smooth branch of $D$, we have
\[
\operatorname{Vol}_\omega(\{T>m\})\le Ce^{-m}
\]
for all sufficiently large $m$. Consequently, if $f\in L^\infty$, inserting it into \eqref{eq:constructive-tail-estimate 1}, we have
\begin{equation}\label{eq:constructive-bounded-source-rate}
\lVert Lv_m-f\rVert_{L_0^2}
\le Ce^{-m/2}\lVert f\rVert_{L^\infty}
+\frac{C}{m}\lVert f\rVert_{L_0^2}.
\end{equation}
\end{rem}

\begin{rem}
For every $h\in\overline{\mathbf{h}_{\mathbin{\!/\mkern-5mu/\!},\mathbb{R}}^D}$, we have
\begin{equation}\label{eq:constructive-residual-orthogonal}
\int_M(Lv_m-f)h\,\omega^n
=\int_Mv_mLh\,\omega^n-\int_Mfh\,\omega^n=0.
\end{equation}
In particular, $\int_M(Lv_m-f)\omega^n=0$. 
\end{rem}


\section{Reductivity of holomorphic vector fields}\label{section: Reductivity of holomorphic vector fields}
In this section, we prove the analogue of a proposition established by Calabi in \cite{C} in the smooth case. Auvray proved the Hodge decomposition of holomorphic vector fields parallel to $D$ in \cite[Proposition 1.2]{A4}.
\begin{prop}\label{reductivity}
Let $\omega$ be a Poincar\'e type extremal K\"ahler metric. Then $\omega$ determines a unique semidirect-sum splitting of the Lie algebra $\mathbf{h}^D_{\parallelsum}$:
\begin{equation*}
\mathbf{h}^D_{\parallelsum}=\mathbf{a}^D_{\parallelsum}(M) \oplus \mathbf{h}^{D}_{\parallelsum,\mathbb{C}},
\end{equation*}
where $\mathbf{a}^D_{\parallelsum}(M)$ and $\mathbf{h}^{D}_{\parallelsum,\mathbb{C}}$ are defined in the Section 3.
\end{prop}

We need the following Lemma:
\begin{lem}\label{lem 5.2}
Let $\omega$ and $\omega_0$ be two K\"ahler metrics with $\omega =\omega_0 +dd^c \varphi$. Let $\alpha=\alpha_{\bar k} d \bar z^k$ be a $(0,1)$ form. Suppose that $v=g_0^{i \bar k}\alpha_{\bar k} \frac{\partial}{\partial z^i}$ is a holomorphic vector field. Then we have that:
 \begin{equation*}
g_0^{i \bar k}\alpha_{\bar k} = g^{i \bar k}\alpha_{\bar k}+ g^{i \bar k}(v(\varphi))_{\bar k}.
 \end{equation*}
\end{lem}
\begin{proof}
We can calculate that:
 \begin{equation*}
v^i (g_{i\bar l}-\varphi_{i\bar l})= v^i (g_0)_{i\bar l}=g_0^{i\bar j}\alpha_{\bar j}(g_0)_{i\bar l}=\alpha_{\bar l}.
 \end{equation*}
 Multiply both sides of the above formula by $g^{k\bar l}$ and sum over $l$. Using the fact that $v$ is holomorphic, we obtain
  \begin{equation*}
v^i- g^{i\bar k} (v(\varphi))_{\bar k}= g^{i \bar l} \alpha_{\bar l}.
 \end{equation*}
 This concludes the proof of the lemma.
\end{proof}

\begin{proof}
(of the Proposition \ref{reductivity})
Let $\omega_0$ be a K\"ahler form on $X$. Then by the definition of the Poincar\'e type K\"ahler metric, there exists a function $\varphi$ such that
\begin{equation*}
    \omega= \omega_0 +dd^c \varphi.
\end{equation*}
For any $Z\in \mathbf{h}_{\parallelsum}^D$, denote $$\xi_{\omega_0}^Z=Z^i (g_0)_{i\bar j}d \bar z^j$$ as the dual 1-form of $Z$ with respect to $\omega_0$ which is $\bar \partial $-closed. We can use the traditional Hodge decomposition to get that:
\begin{equation*}
  Z^i (g_0)_{i\bar j}d \bar z^j =(\xi_{harm}^Z)_{\bar \beta} d \bar z^{\beta}+ \psi_{\bar \beta} d \bar z^{\beta}.
\end{equation*}
Here $\psi$ is a smooth function on $X$ and $\xi_{harm}^Z$ is a harmonic $(0,1)$-form which implies that $\xi_{harm}^Z$ is conjugate holomorphic.

Define $\uparrow_{\omega_0}$ by $\uparrow_{\omega_0} (\alpha_{\bar i} d \bar z^i)= \alpha_{\bar i}g_0^{j \bar i} \frac{\partial}{\partial z^j}$. Then we apply $\uparrow_{\omega_0}$ (using $\omega_0$) to both sides to get that:
\begin{equation*}
Z= \uparrow_{\omega_0} \xi_{harm}^Z + \nabla^{1,0}_{\omega_0}\psi.
\end{equation*}
We further use the Lemma \ref{lem 5.2} to get that:
\begin{equation*}
\uparrow_{\omega_0} \xi_{harm}^Z + \nabla^{1,0}_{\omega_0}\psi= \uparrow_{\omega} \xi_{harm}^Z +\nabla^{1,0}_{\omega}(\psi +Z(\varphi)).
\end{equation*}
Here $\varphi$ is a function such that $\omega = \omega_0 + dd^c \varphi$.
Denote $F=\psi+Z(\varphi)$. We get that:
\begin{equation}\label{decomposition Z}
   Z=\uparrow_{\omega} \xi_{harm}^Z + \uparrow_{\omega}  \bar \partial F.
\end{equation}
Using the fact that $Z$ is holomorphic, we have that:
\begin{equation*}
    0=\bar \partial Z =\bar \partial \uparrow_{\omega} \xi_{harm}^Z +\bar \partial \uparrow_{\omega} \bar \partial F.
\end{equation*}

Recall that $\mathcal{D}=\mathcal{D}_\omega=\bar\partial\nabla_\omega^{1,0}$ and $L=\mathcal{D}^*\mathcal{D}$. Applying $\mathcal{D}^*$ to the above formula gives
\begin{equation*}
\begin{split}
        LF & =\mathcal{D}^* \mathcal{D} F =-\mathcal{D}^* \bar \partial (\uparrow_{\omega} \xi_{harm}^Z)=-g^{\alpha \bar \beta}(\xi_{harm}^Z)_{\bar \beta}^{\,\,\,,\gamma}{}_{\gamma \alpha}\\
&=-g^{\alpha \bar \beta} (\xi_{harm}^Z)_{\bar \beta }^{\,\,\,,\gamma}{}_{\alpha \gamma } =-g^{\alpha \bar \beta}((\xi_{harm}^Z)_{\bar \beta, \alpha}{}^{\gamma}{}_{\gamma}+((\xi_{harm}^Z)_{\bar \tau}R^{\bar \tau}_{\,\,\,\bar \beta}{}^{\gamma}{}_{\alpha})_{,\gamma}).
\end{split}
\end{equation*}
Since $\xi_{harm}^{Z}$ is conjugate-holomorphic, we have that $(\xi_{harm}^Z)_{\bar \beta, \alpha}=0$. Thus we arrive at
\begin{equation}\label{DF}
   L F =-((\xi_{harm}^Z)_{\bar \tau}R^{\bar \tau \gamma})_{,\gamma}=-(\xi_{harm}^Z)_{\bar \tau} R^{, \bar \tau}.
\end{equation}

Since $\xi_{harm}^Z$ is bounded with respect to $\omega_0$, it is also bounded with respect to the Poincar\'e metric $\omega$. Quasi-coordinates show that all covariant derivatives of the metric coefficients are bounded with respect to a fixed background Poincar\'e type K\"ahler metric; hence the covariant derivatives of its scalar curvature are bounded as well. Thus
\begin{equation*}
    |(\xi_{harm}^Z)_{\bar \tau}R^{,\bar \tau}| \le |\xi_{harm}^Z|_{\omega}|\nabla_{\omega}R|_{\omega}\le C.
\end{equation*}

Since $(\xi_{harm}^Z)_{\bar \tau}R^{,\bar \tau}$ is the inner product of a conjugate-holomorphic $(0,1)$-form and a conjugate-holomorphic tangent-bundle section, it is a conjugate-holomorphic function on $X\setminus D$.

We have shown that this function is bounded, so it can be extended to $X$ and as a result, it is a constant $C$. Then we can use  (\ref{DF}) to get that
\begin{equation*}
    L F=C.
\end{equation*}
Note that $F$ is bounded together with all covariant derivatives with respect to a Poincar\'e type K\"ahler metric. The coefficients of $L$ and their derivatives of any order are also bounded with respect to a Poincar\'e type K\"ahler metric. This implies that we can use the Lemma \ref{G-S} to do integration by parts without a boundary term to get that:
\begin{equation*}
    \int_M C \omega^n =\int_M L F \omega^n =0.
\end{equation*}
This implies that $C$ is zero. So we have that $LF=0$. Consequently,
$$\int_M \bar F\cdot L F\cdot\om^n=0.$$

Since $L=\mathcal{D}^*\mathcal{D}$, we get $\mathcal{D}F=0$, so the complex gradient $\nabla^{1,0}F$ is holomorphic. Returning to \eqref{decomposition Z}, we conclude that $\uparrow_{\omega}\xi_{harm}^Z$ is holomorphic. Since $\xi_{harm}^Z$ is harmonic, it is also antiholomorphic; hence $\uparrow_{\omega}\xi_{harm}^Z$ is parallel. This completes the proof of Calabi's decomposition of holomorphic vector fields.
\end{proof}

\section{Uniqueness of the Poincar\'e type cscK metrics}
\label{section :Uniqueness of the Poincare type cscK metrics}
In this section, we prove Theorem \ref{thm:main}.

\begin{defn}
We call $d\mu$ a Poincar\'e type volume form if it is a smooth positive volume form on $M$ such that, for some constant $C$, $C^{-1}\omega^n \le d\mu \le C\omega^n$, and all derivatives of $d\mu/\omega^n$ are bounded with respect to $\omega$. We also assume that $\int_X d\mu=\int_X\omega^n$.
\end{defn}

\subsection{Twisted functional}
Recall the definition of $\mathcal E$ in \eqref{mathcal E}. We define 
\begin{align*}
\mathcal{F}_{\mu}(\varphi)=I_{\mu}(\varphi)-\frac{1}{n+1}\mathcal{E}(\varphi),\quad I_{\mu}(\varphi)=\int_X \varphi d\mu.
\end{align*}
According to Lemma \ref{functional first order}, we have
    \begin{equation*}
    d \mathcal{F}_{\mu}|_\varphi (v)=\int_X v 
        [d\mu- \omega_\varphi^n]. 
    \end{equation*}
Furthermore,
    \begin{equation*}
    d^2 \mathcal{F}_{\mu}|_\varphi (v,w)= \int_X v\Delta_\varphi w\,\omega_\varphi^n
    =\int_X\langle \partial v,\partial w\rangle_{\omega_\varphi}\,\omega_\varphi^n,
    \end{equation*}
which is positive semidefinite and is positive definite modulo constants.

We write $\overline v_\mu=V^{-1}\int_Mv\,d\mu$, for an integrable function $v$.

Then we have the following lemma:
\begin{lem}\label{lem 6.2}
Let $\mu$ be a Poincar\'e type volume form. $I_{\mu}$ is strictly convex along $C^{1,1}$-geodesics in the sense that if $u_t$ is a Poincar\'e type $C^{1,1}$-geodesic , then
\begin{equation*}
(\frac{d}{dt}\Big|^-_{t=1} - \frac{d}{dt}\Big|^+_{t=0}) I_{\mu}(u_t)   \ge c\int_0^1\int_M|\dot u_t-a_t|^2\,d\mu\,dt,
\end{equation*}
where $a_t$ is the average of $\dot{u}_t$ with respect to $\mu$ and $c>0$ depends only on $\mu$, $\omega$, and the uniform $C^{1,1}$ bounds for the geodesic approximation.
Consequently,
\begin{equation}\label{eq:endpoint-variance}
\left\|u_1-u_0-\overline{(u_1-u_0)}_\mu\right\|_{L^2(\mu)}^2
\le c^{-1}\bigl[(I_\mu\circ u)'(1^-)-(I_\mu\circ u)'(0^+)\bigr].
\end{equation}
\end{lem}
\begin{proof}
Let $u_t^{\epsilon}$ be the $\epsilon$-geodesic approximating $u_t$ given by the Lemma \ref{eps geodesic}. Note that $u_t^{\epsilon}$ is a subgeodesic, meaning that:
\begin{equation*}
\ddot{u_t}^{\epsilon} \ge |\bar \partial \dot{u_t}^{\epsilon}|_{\omega_{u_t^{\epsilon}}}^2.
\end{equation*}
Denote $f^{\epsilon}(t)= I_{\mu}(u_t^{\epsilon})$. We can calculate that:
\begin{equation*}
(f^{\epsilon})''(t) = \int_X \ddot{u_t^{\epsilon}}\,d\mu \ge \int_X |\bar \partial \dot{u}^{\epsilon}_t|_{\omega_{u_t^{\epsilon}}}^2 d\mu \ge \frac{1}{C} \int_X  |\bar \partial \dot{u}^{\epsilon}_t|_{\omega}^2 d\mu.
\end{equation*}
Here we use the $C^{1,1}$ estimate for the $\epsilon$-geodesic: $\omega_{u_t^{\epsilon}}\le C\omega$. Using the Poincar\'e inequality from Lemma \ref{poincare ine} and its remark below, we have
\begin{equation*}
\int_X  |\bar \partial \dot{u}^{\epsilon}_t|_{\omega}^2 d\mu \ge C \int_X |\dot{u}^{\epsilon}_t -a^{\epsilon}_t|^2 d\mu,
\end{equation*}
where $a^{\epsilon}_t$ is the average of $\dot{u}^{\epsilon}_t$ with respect to $\mu$ and $C$ depends on $\mu$, $\omega$ and $X$.
Hence we have that:
\begin{equation*}
(f^{\epsilon})''(t) \ge C\int_X |\dot{u}^{\epsilon}_t- a_t^{\epsilon}|^2 d\mu.
\end{equation*}
Integrate the above formula from $0$ to $1$, we have that:
\begin{equation}\label{f diff}
f'^{\epsilon}(1)-f'^{\epsilon}(0) \ge C \int_0^1 dt  \int_X |\dot{u}^{\epsilon}_t- a^{\epsilon}_t|^2 d\mu.
\end{equation}
Note that $f'(t)= \int_X \dot{u} d\mu$ and $f'^{\epsilon}(t)= \int_X \dot{u}^{\epsilon} d\mu$. Since $u^{\epsilon}$ converge to $u$ in $C^{1,\alpha}_{loc}((X\setminus D) \times [0,1])$ and $|\dot{u^{\epsilon}}|+|\dot{u}|\le C \mathbf{u}$ for some uniform constant $C$ independent of $\epsilon$, we have that 
\begin{equation*}
\lim_{\epsilon \rightarrow 0} f'^{\epsilon} (t)= f'(t)
\end{equation*}
and
\begin{equation*}
\lim_{\epsilon \rightarrow 0}\int_0^1 dt  \int_X |\dot{u}^{\epsilon}_t- a^{\epsilon}_t|^2 d\mu = \int_0^1 dt  \int_X |\dot{u}_t- a_t|^2 d\mu.
\end{equation*}
Letting $\epsilon$ tend to zero in \eqref{f diff}, we obtain
\begin{equation*}
f'(1)-f'(0) \ge C \int_0^1 dt  \int_X |\dot{u}_t- a_t|^2 d\mu .
\end{equation*}

Applying Cauchy-Schwarz and integrating in $M$ to
\[
u_1-u_0-\overline{(u_1-u_0)}_\mu
=\int_0^1(\dot u_t-a(t))\,dt,
\]
we prove \eqref{eq:endpoint-variance}.
\end{proof}

Auvray proved the following Poincar\'e inequality for Poincar\'e type K\"ahler metrics (cf. \cite[Lemma 1.11]{A}):
\begin{lem}\label{poincare ine}
Assume $M=X\setminus D$ is equipped with a Poincar\'e type K\"ahler metric $\omega$. Then there exists a constant $C_P>0$ such that for all $v\in W^{1,2}(M,\omega)$, we have 
\begin{equation*}
\int_{X}|v-a|^2  \omega^n \le C_P \int_{X}|dv|_{\omega}^2 \omega^n,
\end{equation*}
where $a=\big(\int_Mv\omega^n\big)/\big(\int_M\omega^n\big)$.\end{lem}
\begin{rem}
The Poincar\'e inequality also holds when replacing $\om^n$ by a Poincar\'e volume form $d\mu$: there is a constant $C_\mu$ such that
\begin{equation}\label{eq:mu-poincare}
\int_X|v-\overline v_\mu|^2\,d\mu
\le C_\mu\int_X|\partial v|_\omega^2\,d\mu
\end{equation}
for $v\in W^{1,2}$. The constant depends only on $C_P$ and the upper and lower bounds on $\frac{d\mu}{\omega^n}$.

\end{rem}
\subsection{K-energy}

We also need the following Lemma:
\begin{lem}\label{lem 6.9}
Given two Poincar\'e type K\"ahler potentials, let $u_t$ be the corresponding Poincar\'e type $C^{1,1}$ geodesic. Then
\begin{equation*}
\lim_{ t \rightarrow 0^+} \frac{\mathcal{M}(u_t)-\mathcal{M}(u_0)}{t} \ge \int_X (-R_{\omega_{u_0}}+\bar R) \frac{d u_t}{dt}\Big|_{t^+=0} \omega_{u_0}^n
\end{equation*}
\end{lem}
\begin{proof}
The proof is similar to the proof of \cite[Lemma 3.5]{BB}.  We only need to deal with the entropy part because other parts in the decomposition of the $K$-energy are differentiable along the geodesic using a calculation similar to the Lemma \ref{functional first order}. We want to show that:
\begin{equation*}
\begin{split}
&\lim_{t \rightarrow 0^+} (\frac{1}{t})(H_{\omega^n}(\omega_{u_t}^n)- H_{\omega^n}(\omega_{u_0}^n)) \ge \\
& -n \int_X \frac{d u_t}{dt}\Big|_{t=0} \operatorname{Ric}(\omega_{u_0}) \wedge \omega_{u_0}^{n-1} + n \int_X \frac{d u_t}{dt}\Big|_{t=0} \operatorname{Ric}(\omega) \wedge \omega_{u_0}^{n-1}.
\end{split}
\end{equation*}
We use the convexity of entropy on the affine space of finite measures having the same total mass:
\begin{equation*}
H_{\omega^n}(\nu_1) - H_{\omega^n}(\nu_0) \ge \frac{d}{ds}\Big|_{s=0} H_{\omega^n} (\nu_s),
\end{equation*}
with $\nu_s= s\nu_1+(1-s)\nu_0$. Since $\log (\frac{\nu_s}{\omega^n}) \nu_s$ is convex in $s$, we can use the monotone convergence to get that:
\begin{equation*}
\frac{d}{ds}\Big|_{s=0} H_{\omega^n}(\nu_s) = \int_X \left(1+\log \frac{d\nu_0}{\omega^n}\right) (d\nu_1-d\nu_0)
=\int_X \log \frac{d\nu_0}{\omega^n}\,(d\nu_1-d\nu_0),
\end{equation*}
where the last equality uses $\nu_0(X)=\nu_1(X)$.
Taking $\nu_1=\omega_{u_t}^n$ and $\nu_0=\omega_{u_0}^n$, we obtain
\begin{equation*}
\begin{split}
 \frac{1}{t}(H_{\omega^n}(\omega_{u_t}^n) - H_{\omega^n}(\omega_{u_0}^n))
& \ge \int_X \log(\frac{\omega_{u_0}^n}{\mu_0})\frac{1}{t}(\omega_{u_t}^n -\omega_{u_0}^n) \\
& = \int_X \log(\frac{\omega_{u_0}^n}{\mu_0}) \frac{1}{t} dd^c (u_t- u_0) \wedge \sum_{i=0}^{n-1}\omega_{u_t}^i \wedge \omega_{u_0}^{n-1-i}\\
& = \int_X dd^c \log(\frac{\omega_{u_0}^n}{\omega^n}) \frac{u_t -u_0}{t} \wedge \sum_{i=0}^{n-1}\omega_{u_t}^i \wedge \omega_{u_0}^{n-1-i} \\
&= \int_X (- Ric_{\omega_{u_0}} + Ric_{\omega}) \frac{u_t -u_0}{t} \wedge \sum_{i=0}^{n-1}\omega_{u_t}^i \wedge \omega_{u_0}^{n-1-i}.
\end{split}
\end{equation*}
Letting $t\to0^+$ in the above formula gives
\begin{equation*}
\lim_{t \rightarrow 0+}\frac{1}{t}(H_{\omega^n}(\omega_{u_t}^n) - H_{\omega^n}(\omega_{u_0}^n)) \ge n\int_X \dot u_0\,\big(-\operatorname{Ric}(\omega_{u_0})+\operatorname{Ric}(\omega)\big)\wedge\omega_{u_0}^{n-1}.
\end{equation*}
This proves the entropy estimate. Combining this estimate with the first variations of $\frac{\bar R}{n+1}\mathcal E$ and $-\mathcal E^{\operatorname{Ric}(\omega)}$ in \eqref{decom kenergy} gives the stated inequality.
\end{proof}

\subsection{Difference between Poincar\'e type cscK metrics}

\begin{lem}\label{lem uvbound}
Let $\omega_u$ and $\omega_v$ be two Poincar\'e type cscK metrics in the same cohomology class. Then there exists a constant $C$ depending on $\omega$, $u$ and $v$ such that $|u-v|\le C$.
\end{lem}
\begin{proof}
It was pointed out in the section 3.1 of \cite{A3} that \cite{A2} proved that if a metric of
Poincar\'e type $\omega_\varphi=\omega+dd^c\varphi$ has constant
scalar curvature, then
\[
\varphi-\sum_{j=1}^{k}a_j \mathbf{u}_j \in C^\infty(X\setminus D),
\]
that is, the remainder and all its derivatives are
bounded with respect to the model
Poincar\'e type metric $\omega$, where the $a_j<A$ and given by
\[
\frac{2}{A-a_j}=\bar{s}_{D_j}-\bar{s}.
\]
Here $\bar{s}$ (respectively $\bar{s}_{D_j}$) denotes the mean
scalar curvature of Poincar\'e type metrics of class $\omega_X$
(respectively $\omega_{D_j}$) on $X\setminus D$
(respectively on $D_j\setminus\sum_{\ell\neq j}D_\ell$)---recall
the formulae
\[
\bar{s}
=
-2\pi n\,
\frac{c_1(K_X[D])\cdot[\omega_X]^{n-1}}
     {[\omega_X]^n}
\qquad\text{and}\qquad
\bar{s}_{D_j}
=
-2\pi (n-1)\,
\frac{c_1(D_j)\cdot c_1(K_X[D])\cdot[\omega_X]^{n-2}}
     {c_1(D_j)\cdot[\omega_X]^{n-1}},
\]
typical of the Poincar\'e setting. As a result, different Poincar\'e type cscK metrics have the same asymptotic along the divisor and $u-v$ is bounded.
\end{proof}

Then we can show the following lemma:
\begin{lem}\label{lem 6.5}
Let $\omega_u$ and $\omega_v$ be two Poincar\'e type cscK metrics.  Denote $\mu=\omega_u^n$ and $\nu= \omega_v^n$.  Then there exists a constant $C$ depending on $\omega$, $u$ and $v$ such that for any K\"ahler potential $\phi$, 
\begin{equation*}
|I_{\mu}(\phi)-I_{\nu}(\phi)|\le C.
\end{equation*}
\end{lem}
\begin{proof}
By Lemma \ref{lem uvbound}, $|u-v|$ is bounded. We can calculate that:
\begin{equation*}
I_{\mu}(\phi)-I_{\nu}(\phi)=\int_X \phi(\omega_u^n-\omega_v^n) =\int_X \phi\cdot [dd^c(u-v)\wedge \sum \omega_u^{n-k-1}\wedge \omega_v^{k}].
\end{equation*}
Then we can use the Lemma \ref{G-S} to get that:
\begin{equation}\label{diff i mu}
\begin{split}
I_{\mu}(\phi)-I_{\nu}(\phi)&=\int_X (u-v) dd^c \phi \wedge \sum \omega_u^{n-k-1}\wedge \omega_v^{k}\\
&=\int_X (u-v) (\omega_{\phi}-\omega) \wedge \sum \omega_u^{n-k-1}\wedge \omega_v^{k}\\
&\le \|u-v\|_\infty \int_M(\omega_\phi+\omega)\wedge \sum_{k=0}^{n-1}\omega_u^{n-1-k}\wedge\omega_v^k\\ &=2nV\|u-v\|_\infty. 
\end{split}
\end{equation}
Here all forms represent the same fixed class, and the last identity follows from the integrations by parts already used in Section \ref{subsection: Energy functionals}.
Since we have shown that $|u-v|$ is bounded, we have that the above formula is bounded independent of $\phi$.
\end{proof}

\subsection{Linearised equation}
Consider the twisted $K$-energy
\begin{equation*}
\mathcal{M}_s \triangleq \mathcal{M} +s \mathcal{F}_{\mu}.
\end{equation*}

We recall that
\begin{equation}\label{eq:functionals}
\begin{split}
&\mathcal{F}_\mu(\vphi)=I_\mu(\vphi)-\frac{1}{n+1}E(\vphi),\\
&E(\vphi)=\int_M\vphi\sum_{k=0}^n\omega_\vphi^{n-k}\wedge\omega^k,\quad
I_\mu(\vphi)=\int_M\vphi\,d\mu.
\end{split}
\end{equation}
where $\mu(M)=V$. The first variations are
\begin{equation}\label{eq:functional-variations}
\begin{split}
dE|_\vphi(v)=(n+1)\int_Mv\,\omega_\vphi^n,\quad
d\mathcal{F}_\mu|_\vphi(v)=\int_Mv\,(d\mu-\omega_\vphi^n).
\end{split}
\end{equation}

Due to the Corollary \ref{cri mabuchi}, we have
    \begin{equation}\label{dMs}
        d\mathcal{M}_s|_\varphi(v) =- \int_X v\cdot \bigg\{ R_\varphi-\bar R -s  
        \bigg[\frac{d\mu}{\omega_\varphi^n}-1\bigg] \bigg\}   \omega_\varphi^n.
    \end{equation}

According to Mabuchi \cite{M} and Donaldson \cite{D}, if $\omega_\varphi$ is cscK, then the Hessian of $\mathcal{M}$ at $\varphi$ is
\begin{equation}\label{Hessian M}
d^2\mathcal M\vert_{\varphi}(v,w)= Re \int_X (\mathcal{D}_\varphi^* \mathcal{D}_\varphi v)\cdot w \cdot \omega_\varphi^n.
\end{equation}
The conclusion is still valid for any Poincar\'e type cscK metric $\vphi$ and $v,w\in C^{4,\alpha}_0$, by recalling the cutoff functions and Gaffney's theorem, as shown in the Section \ref{subsection: Energy functionals}.

\subsection{Gauge fixing}
\subsubsection{The tangent space}
Note that for any $Y\in \mathbf{h}^D_{\parallelsum}$, the holomorphic transformation induced by $Y$ sends a Poincar\'e type cscK metric to another Poincar\'e type cscK metric. We can define $Iso_0^D(X,\omega)$ as the set of elements in $Aut_0^D(X)$ that fix $\omega$. Define the orbit of \(\omega\):
\begin{equation*}
    \mathcal{O}_\om=Aut_0^D(X)/Iso_0^D(X,\omega)
\end{equation*}

For \(\omega'\in\mathcal O_\omega\), let $\Gamma \subset \mathbf{h}^D_{\parallelsum}$ be the vector space of real part $Y=Re Z$ of holomorphic vector fields $Z\in \mathbf{h}^D_{\parallelsum}$ such that the Lie derivative $\mathcal{L}_Y \omega' =0$. Then the tangent space of $\mathcal{O}_\om$ can be written as 
\begin{equation*}
T_{\om'}{\mathcal{O}_\om}=  \{Y=Re Z\vert Z\in \mathbf{h}^D_{\parallelsum}\}/ \Gamma.
\end{equation*}

From Proposition \ref{reductivity}, a holomorphic \((1,0)\)-vector field $Z$ has the decomposition
$$ Z=P+\nabla_\omega^{1,0}F, $$
where \(P\) is parallel and
$ F=\psi+Z(\varphi) $.
Since $\om$ is cscK, $L$ is real. Both $\nabla_\omega^{1,0}Re F$ and $\nabla_\omega^{1,0}Im F$ are holomorphic. Then
$$ Re Z = ReP+Re\nabla_\omega^{1,0}Re F+JRe\nabla_\omega^{1,0}Im F. $$
Both the first and third terms preserve \(\omega\). Also, 
$$\mathcal L_{Y_{Re F}}\omega=dd^cRe F. $$
Therefore, every tangent direction to the orbit is represented by a real holomorphy potential.
Conversely, every real holomorphy potential \(h\) generates such an orbit direction and the normalisation makes the representative unique. Thus, we have proved an equivalent characterisation of the tangent space.
\begin{lem}\label{tangent space of orbit}
The tangent space $T_{\omega'}\mathcal O_\omega$ is naturally identified with
$$ \mathcal{H}_{\omega'}
=\left\{h\in\overline{\mathbf{h}^D_{\parallelsum,\mathbb{R}}}(\omega'):\int_Mh\,(\omega')^n=0\right\}. $$
\end{lem}

\subsubsection{Flow generated by the holomorphic gradient is geodesic}
Fix $ h\in\mathcal H_\omega$,
by Lemma \ref{ker a}, the holomorphic gradient $Y_h=Re \nabla_\omega^{1,0}h$ extends to \(X\) and is tangent to \(D\). Denote its real flow by \(f_s\). Then
\begin{align}\label{generate geodesic} \tilde u_t=\int_0^tf_s^*h\,ds
\end{align}
satisfies the smooth geodesic equation
$$ \ddot  {\tilde u}_t-|\partial\dot  {\tilde u}_t|_{\omega_{ \tilde u_t}}^2=0 $$
and normalisation condition $\mathcal E( \tilde u_t)=0$, via the computation
\[\frac{1}{n+1}\frac{d}{dt}\mathcal{E}( \tilde u_t)
= \int_M\dot  {\tilde u}_t\,\omega_{ \tilde u_t}^n = \int_Mf_t^*(h\omega^n) = \int_Mh\,\omega^n=0. \]

The smooth, positive-definite, \(G\)-invariant Riemannian metric defined on the orbit is
$$ \mathcal G_{\omega'}(h,k) = \int_Mhk\,(\omega')^n. $$
Consider the energy of a normalised path in the orbit
$$ A(u) = \frac12\int_a^b\int_M\dot u_t^2\,\omega_{u_t}^n\,dt,\quad \mathcal{E}(u_t)=0.$$
Applying the bounded cusp derivatives, finite volume, the cutoffs and Gaffney argument we get its variation
$$ \delta A = -\int_a^b\int_M v_t \left( \ddot u_t-|\partial\dot u_t|_{\omega_t}^2 \right)\omega_t^n\,dt. $$
Thus $\tilde u_t$ are geodesics for the finite-dimensional orbit metric. 

Every initial tangent vector at every point of the orbit generates a geodesic $\tilde u_t$ for all real time. Hence the finite-dimensional orbit metric is geodesically complete. Since the orbit is connected, Hopf-Rinow implies that every orbit point can be joined to the base point by one of these geodesics. That gives a surjective map
$$ \Theta:\mathcal H_\omega\longrightarrow\mathcal O_\omega, \quad \Theta(h)=(f_1)^*\omega.$$

\subsubsection{A minimizing gauge}
\begin{lem}\label{gauge}
Let $\omega$ be a Poincar\'e type cscK metric. For any Poincar\'e type cscK metric $\omega_1$ in the same class as \(\omega\), we can define $\mu= \omega_1^n$. Then $\mathcal{F}_{\mu}$ has a minimum and hence a critical point, $\omega_u$, on $\mathcal{O}_{\omega}$. 
At this metric,
\begin{equation}\label{eq:gauge-orthogonality}
\int_M h\,(d\mu-\omega_u^n)=0,
\quad h\in\overline{\mathbf{h}^D_{\parallelsum,\mathbb{R}}}(\omega_u).
\end{equation}
\end{lem}
\begin{proof}
If $\mu_0=\omega^n$, then $0$ is a critical point of $\mathcal{F}_{\mu_0}$ and $\mathcal{F}_{\mu_0}(0)=0$. Since $I_{\mu_0}$ is strictly convex along each ray using Lemma \ref{lem 6.2} and $\mathcal{E}$ is affine along each ray, we have that $\mathcal{F}_{\mu_0}$  is strictly convex along each ray. 
Actually, along the geodesic $\tilde u_t$ generated by \(h\neq0\) in \eqref{generate geodesic},
we have
$$ \frac{d^2}{dt^2}\mathcal{F}_{\mu_0}({\tilde u}_t) = \int_M|\partial\dot {\tilde u}_t|_{\omega_t}^2\,\omega^n>0.$$
Since \(\mathcal{F}_{\mu_0}'(0)=0\), this implies
$$ \mathcal{F}_{\mu_0}(\Theta(h))>0. $$
By compactness of its unit sphere and continuity, we get 
$$ c:=\min_{\|h\|=1}\mathcal{F}_{\mu_0}(\Theta(h))>0. $$
For \(r\ge1\), convexity along the same geodesic gives
$$ \mathcal{F}_{\mu_0}(\Theta(rh))\ge cr,\qquad \|h\|=1.$$
Then each sublevel set of \(\mathcal{F}_{\mu_0}\) is contained in the continuous image of a compact ball in \(\mathcal H_\omega\). Since the sublevel set is closed in the Hausdorff orbit, it is compact.
So, we have that $\mathcal{F}_{\mu_0}$ is proper on $\mathcal{O}_{\omega}$, if $\mu_0=\omega^n$.

Then we can use the Lemma \ref{lem 6.5} to show that $$ |\mathcal F_\mu-\mathcal F_{\omega^n}|\le C, $$ and furthermore $\mathcal{F}_{\mu}$ is also proper. Then we can find a minimum $\omega_u$ of $\mathcal{F}_{\mu}$ over $\mathcal{O}_{\omega}$. Note that in that case the connected orbit is a point, the conclusion
is immediate.

At a minimising metric $\omega_u$, the differential vanishes on all orbit tangents, which are precisely $ \mathcal H_{\omega_u} $ by Lemma~\ref{tangent space of orbit}. Therefore, \eqref{eq:functional-variations} proves the orthogonality \eqref{eq:gauge-orthogonality} for these normalised potentials on the whole space $\overline{\mathbf{h}^D_{\parallelsum,\mathbb{R}}}(\omega_u)$. Also, that orthogonality to normalised holomorphy potentials extends to constant $h$ because the two measures have equal total mass.
\end{proof}

Now we are ready to prove the main theorems in this section:
\subsection{Proof of Theorem \ref{thm:main}}

In this subsection, we prove Theorem \ref{thm:main}.

\subsubsection{Gauge fixing and the nondecaying source}
We apply Lemma~\ref{gauge} with 
\[\mu=\omega_2^n,\quad\om=\om_1.\] 
Then we find $g\in Aut_0^D(X)$ such that $g^* \omega_1=\omega_{1u}$ is the minimiser of $\mathcal{F}_{\mu}$ over $\mathcal{O}_{\omega_1}$. That is $\mathcal{F}_{\mu}$ has a minimum $\omega_{1u}$, on $\mathcal{O}_{\omega_1}$ and
\begin{equation}\label{eq:gauge-orthogonality 1}
\int_M h\,(d\mu-\omega_{1u}^n)=0,
\quad h\in\overline{\mathbf{h}^D_{\parallelsum,\mathbb{R}}}(\omega_{1u}).
\end{equation}

Write $\om=\omega_{1u}$ and write
\begin{equation*}
\omega_2=\omega+dd^c b,\qquad
 f=1-\frac{d\mu}{\omega^n}.
\end{equation*}
From Lemma~\ref{lem uvbound}, we have $b\in C_0^\infty(M)$.
Also, $f\in L^2(\omega^n)$. By \eqref{eq:gauge-orthogonality 1}, we get
\begin{equation*}
\int_Mfh\,\omega^n=0,\quad h\in\overline{\mathbf{h}^D_{\parallelsum,\mathbb{R}}}(\omega).
\end{equation*}

According to Proposition~\ref{kernel approximation}, there exists $v_j\in C_c^\infty(M)$ such that
\begin{equation}\label{eq:residual}
e_j=L v_j-f,\qquad
\varepsilon_j=\lVert e_j\rVert_{L^2(\omega^n)}\rightarrow0,\quad j\rightarrow\infty.
\end{equation}

\subsubsection{Uniform geodesics for fixed $j$}
We choose a sufficiently large $m$.
Fix a small number $\delta>0$. For each fixed $j$, choose $s_j>0$ so small that
\begin{equation}\label{svj}
0<s\le s_j
\quad\Longrightarrow\quad
\lVert s v_j\rVert_{C_0^{m,\alpha}}\le\delta,
\qquad
\frac12\omega\le\omega+sdd^c v_j\le\frac32\omega.
\end{equation}

Let $u_t^{s,j}$ be the weak geodesic with endpoints
\begin{equation*}
u_0^{s,j}=s v_j,\qquad u_1^{s,j}=b.
\end{equation*}
All these endpoint pairs lie in one set with bounded $C_0^{m,\alpha}$ norms and fixed positive metric bounds, for $0<s\le s_j$. Fix a sufficiently large integer $m$. Therefore Lemma~\ref{C11 geodesic} gives constants independent of $j$ and of such $s$ with
\begin{equation}\label{eq:uniform-geodesic}
0\le\omega_{u_t^{s,j,\eps}}\le C\omega,
\qquad |\ddot u_t^{s,j,\eps}|\le C,
\end{equation}
the same spatial upper bound and 
\begin{equation}\label{eq:velocity-bound}
\sup_{j,\ 0<s\le s_j,\ t\in[0,1]}
\left(\lVert \dot u_t^{s,j}\rVert_{L^\infty}
+\lVert \dot u_t^{s,j}\rVert_{L^2(\omega^n)}\right)\leq C<\infty,
\end{equation}
see also Lemma \ref{vdot estimate}.

\subsubsection{Bound $\lVert G_j(s)\rVert_{L^2(\omega^n)} $}
Set
\[
\mathcal{M}_s=\mathcal{M}+s\mathcal{F}_\mu.
\]
At the first endpoint the first variations give
\begin{equation}\label{eq:dMs-density}
|d\mathcal{M}_s|_{s v_j}(w) |= |\int_Mw\,G_j(s)\,\omega^n |
\le\lVert w\rVert_{L^2(\omega^n)} \lVert G_j(s)\rVert_{L^2(\omega^n)} .
\end{equation}
where
\begin{equation*}
G_j(s)=(\bar R-R_{\omega_{s v_j}})\frac{\omega_{s v_j}^n}{\omega^n}
+s\left(\frac{d\mu}{\omega^n}-\frac{\omega_{s v_j}^n}{\omega^n}\right).
\end{equation*}

Because $\omega$ is cscK,
\[
G_j(0)=0,\quad 
\left.\frac{d}{ds}\right|_{s=0}R_{\omega_{s v_j}}=-L v_j.
\]
Using it and $\omega_{s v_j}|_{s=0}=\omega$, we have
\begin{equation*}
\left.\frac{d}{ds}\right|_{s=0}G_j(s)
=L v_j+\frac{d\mu}{\omega^n}-1=L v_j-f=e_j.
\end{equation*}

For every $x\in M$ and $0<s\le s_j$, the fundamental theorem of calculus twice gives
\begin{align*}
G_j(s,x)-G_j(0,x)-sG_j'(0,x)
=\int_0^s\int_0^\tau G_j''(\sigma,x)\,d\sigma\,d\tau
=\int_0^s(s-\sigma)G_j''(\sigma,x)\,d\sigma.
\end{align*}
Using the substitution $\sigma=rs$, we obtain
\[
G_j(s,x)=se_j(x)+s^2Q_{s,j}(x),\quad Q_{s,j}(x)=\int_0^1(1-r)G_j''(rs,x)\,dr.
\] 
We estimate $Q_{s,j}$ in $L^2$ norm 
\begin{align*}
\|Q_{s,j}\|_{L^2(\omega^n)}
\le\int_0^1(1-r)\|G_j''(rs)\|_{L^2(\omega^n)}\,dr
\le\frac12\sup_{0\le\sigma\le s_j}
\|G_j''(\sigma)\|_{L^2(\omega^n)}.
\end{align*}
The right hand side is finite, because of \eqref{svj}.

In conclusion, we have obtained that
\begin{equation}\label{eq:taylor-estimate}
\begin{gathered}
G_j(s)=se_j+s^2Q_{s,j},\quad
\|Q_{s,j}\|_{L^2(\omega^n)}\le C_j,
\quad 0<s\le s_j.
\end{gathered}
\end{equation}

\subsubsection{Bound $\mathcal{M}_s$}
We insert \eqref{eq:residual} and \eqref{eq:taylor-estimate} into \eqref{eq:dMs-density}, using Cauchy-Schwarz,
\begin{equation*}
\left|d\mathcal{M}_s|_{s v_j}(w)\right|
\le\lVert w\rVert_{L^2(\omega^n)}(s\varepsilon_j+C_js^2).
\end{equation*}


Taking $w=\dot u_0^{s,j}$ and using \eqref{eq:velocity-bound} yields
\begin{equation}\label{eq:endpoint-residual-bound}
\left|d\mathcal{M}_s|_{s v_j}(\dot u_0^{s,j})\right|
\le C s\varepsilon_j+C_js^2.
\end{equation}

At $b$, $\omega_b=\omega_2$ is cscK and $d\mu=\omega_b^n$. So,
\begin{equation}\label{eq:second-endpoint-zero}
d\mathcal{M}_s|_b(w)=0.
\end{equation}

Lemma \ref{lem 6.9}, together with \eqref{eq:second-endpoint-zero} and \eqref{eq:endpoint-residual-bound}, give
\begin{equation}\label{eq:slope-upper}
\begin{split}
&\frac{d}{dt}\Big\vert^-_1 (\mathcal{M}_s\circ u^{s,j})
-
\frac{d}{dt}\Big\vert^+_0(\mathcal{M}_s\circ u^{s,j})\\
&\le d\mathcal{M}_s|_b(\dot u_1^{s,j})
-d\mathcal{M}_s|_{s v_j}(\dot u_0^{s,j})\\
&\le C s\varepsilon_j+C_js^2.
\end{split}
\end{equation}

\subsubsection{Bound $I_\mu$}
Since $\mathcal{M}$ is convex (Theorem \ref{thm convex}) and $\mathcal{E}$ is affine $\partial_t^2\mathcal{E}(u_t^{s,j})=0$ along geodesic $u_t^{s,j}$ (Lemma \ref{ddce}), we get that 
\begin{equation*}
\begin{split}
0&\le s\Bigl[\frac{d}{dt}\Big\vert^-_1(I_\mu\circ u^{s,j})-\frac{d}{dt}\Big\vert^+_0(I_\mu\circ u^{s,j})\Bigr]\\
&\le\frac{d}{dt}\Big\vert^-_1(\mathcal{M}_s\circ u^{s,j})-\frac{d}{dt}\Big\vert^+_0(\mathcal{M}_s\circ u^{s,j}).
\end{split}
\end{equation*}
Dividing by $s>0$ and inserting \eqref{eq:slope-upper}, we have
\begin{equation}\label{eq:slope-final}
\frac{d}{dt}\Big\vert^-_1 (I_\mu\circ u^{s,j})-\frac{d}{dt}\Big\vert^+_0(I_\mu\circ u^{s,j})
\le C\varepsilon_j+C_js.
\end{equation}

\subsubsection{Passage to the limit}
By \eqref{eq:uniform-geodesic}, the constant in Lemma~\ref{lem 6.2} is uniform for the family. Put
\[
W_{s,j}=u_1^{s,j}-u_0^{s,j}=b-sv_j,
\qquad
\overline W_{s,j}=\frac1V\int_MW_{s,j}\,d\mu.
\]
Combining \eqref{eq:endpoint-variance} with \eqref{eq:slope-final} gives
\begin{equation}\label{W convergence}
\lVert W_{s,j}-\overline W_{s,j}\rVert_{L^2(\mu)}^2
\le C\varepsilon_j+C_js.
\end{equation}
The constant $C$ is independent of $j$.

Note that $j$ is fixed. We take $s\rightarrow0$. Then 
$$sv_j\to0\text{ in }L^2(\mu),\quad W_{s,j}\to b\text{ in }L^2(\mu).$$ Thus, \eqref{W convergence} gives 
\begin{equation*}
\lVert b-\overline b_\mu\rVert_{L^2(\mu)}^2\le C\varepsilon_j.
\end{equation*}

Finally, let $j\to\infty$, then $\epsilon_j\rightarrow 0$. \eqref{eq:residual} gives
\begin{equation*}
\lVert b-\overline b_\mu\rVert_{L^2(\mu)}^2=0.
\end{equation*}
Thus $b$ is constant and $\omega_2=\omega$.
From the gauge choice 
$ g^*\omega_1=\omega$,
we have
$ \omega_1=(g^{-1})^*\omega_2. $
This completes the proof of Theorem~\ref{thm:main}.


\begin{thebibliography}{9}

\bibitem{Ao}
T. Aoi.
\newblock
A conical approximation of constant scalar curvature K\"ahler metrics of Poincar\'e type.
\newblock
 arXiv preprint arXiv:2210.11862 (2022).

\bibitem{MR4215262}
V. Apostolov, H. Auvray, L. M. Sektnan. 
\newblock
Extremal K\"ahler Poincar\'e type metrics on toric varieties. 
\newblock
J. Geom. Anal. 31 (2021), no. 2, 1223-1290.


\bibitem{A}
H. Auvray. 
\newblock
The space of Poincar\'e type K\"ahler metrics on the complement of a divisor.
\newblock
J. Reine Angew. Math. 722 (2017), 1-64.

\bibitem{A2}
H. Auvray. 
\newblock
Metrics of Poincar\'e type with constant scalar curvature: a topological constraint.
\newblock
J. Lond. Math. Soc. (2) 87 (2013), no. 2, 607-621.

\bibitem{A3}
H. Auvray. 
\newblock
Asymptotic properties of extremal K\"ahler metrics of Poincar\'e type.
\newblock
Proc. Lond. Math. Soc. (3) 115 (2017), no. 4, 813-853.

\bibitem{A4}
H. Auvray. 
\newblock
Note on Poincar\'e type Futaki characters.
\newblock
Ann. Inst. Fourier (Grenoble) 68 (2018), no. 1, 319-344.

\bibitem{BB}
R.Berman, and B. Berndtsson. 
\newblock
Convexity of the $K$-energy on the space of K\"ahler metrics and uniqueness of extremal metrics.
\newblock
J. Amer. Math. Soc. 30 (2017), no. 4, 1165-1196.

\bibitem{B}
H.  Brezis.
\newblock
 Functional analysis, Sobolev spaces and partial differential equations.
\newblock
Universitext. Springer, New York, 2011. xiv+599 pp. ISBN: 978-0-387-70913-0

\bibitem{C}
E. Calabi. 
\newblock
Extremal K\"ahler metrics II.
\newblock
 Differential geometry and complex analysis: a volume dedicated to the memory of Harry Ernest Rauch. Berlin, Heidelberg: Springer Berlin Heidelberg, 1985. 95-114.

\bibitem{C2}
E. Calabi. 
\newblock
The space of K\"ahler metrics.
\newblock
Proceedings of the International Congress of Mathematicians 1954, Vol. 2 page 206-207 (1954), Amsterdam.

\bibitem{C3}
E. Calabi.
\newblock
 Extremal K\"ahler metrics.
\newblock
 In Seminar on Differential Geometry, volume 102 of Ann. of Math. Stud., pages 259–290, Princeton Univ. Press, Princeton, N.J., 1982.

\bibitem{CFH}
X.X. Chen, M. Feldman, and J. Hu. 
\newblock
Geodesic convexity of small neighborhood in the space of K\"ahler potentials.
\newblock
J. Funct. Anal. 279 (2020), no. 7, 108603, 65 pp.

\bibitem{CLP}
X.X. Chen, L. Li, M. Paun. 
\newblock
Approximation of weak geodesics and subharmonicity of Mabuchi energy.
\newblock
Ann. Fac. Sci. Toulouse Math. (6) 25 (2016), no. 5, 935-957.

\bibitem{CPZ}
X.X. Chen, M. Paun, Y. Zeng. 
\newblock
On deformation of extremal metrics.
\newblock
arXiv preprint arXiv:1506.01290 (2015).

\bibitem{CTW}
J. Chu, V. Tosatti and B. Weinkove. 
\newblock
$C^{1,1}$ regularity for degenerate complex Monge-Amp\`ere equations and geodesic rays.
\newblock
Comm. Partial Differential Equations 43 (2018), no. 2, 292-312.

\bibitem{DLempert}
T.  Darvas and L.  Lempert. 
\newblock
Weak geodesics in the space of K\"ahler metrics.
\newblock
Math. Res. Lett. 19 (2012), no. 5, 1127-1135.

\bibitem{D5}
T.  Darvas. 
\newblock
Morse theory and geodesics in the space of K\"ahler metrics.
\newblock
Proc. Amer. Math. Soc. 142 (2014), no. 8, 2775-2782.

 
 \bibitem{D}
 S. K. Donaldson.
 \newblock
  Symmetric spaces, K\"ahler geometry and Hamiltonian dynamics
  \newblock  Northern California Symplectic Geometry Seminar, Amer. Math. Soc. Transl. Ser. 2, vol. 196, Amer. Math. Soc., Providence, RI, 1999, pp. 13–33, DOI 10.1090/trans2/196/02. MR1736211
  
  \bibitem{F}
  Y. Feng. 
  \newblock
  A gluing construction of constant scalar curvature K\"ahler metrics of Poincar\'e type.
  \newblock
   arXiv preprint arXiv:2405.11952 (2024).


\bibitem{MR0062490}
M. P. Gaffney
A special Stokes's theorem for complete Riemannian manifolds.
 \newblock
Ann. of Math. (2) 60 (1954), 140-145.

\bibitem{GT}
D. Gilbarg, N.S. Trudinger
\newblock
Elliptic partial differential equations of second order
\newblock 
Reprint of the 1998 edition. Classics in Mathematics. Springer-Verlag, Berlin, 2001. xiv+517 pp. ISBN: 3-540-41160-7

 
\bibitem{LV}
L. Lempert, and L. Vivas. 
\newblock
Geodesics in the space of K\"ahler metrics.
\newblock
Duke Math.  J.  162(7)(2013), 1369-1381.


\bibitem{M}
T. Mabuchi.
\newblock
 Some symplectic geometry on compact K\"ahler manifolds. I
 \newblock
  Osaka J. Math. 24 (1987), no. 2, 227–252. MR909015

\bibitem{M2}
T. Mabuchi.
\newblock
$K$-energy maps integrating Futaki invariants.
\newblock Tohoku Math. J. (2) 38 (1986), no. 4, 575–593, DOI 10.2748/tmj/1178228410. MR867064

\bibitem{R}
R. M. Range.
\newblock
 Holomorphic functions and integral representations in several complex variables.
 \newblock
Graduate Texts in Mathematics, 108. Springer-Verlag, New York, 1986. xx+386 pp. ISBN: 0-387-96259-X

\bibitem{S}
L. M. Sektnan. 
\newblock
 Blowing up extremal Poincar\'e type manifolds. 
 \newblock
Math. Res. Lett. 30 (2023), no. 1, 185-238.

\bibitem{TY}
G.Tian, S.T. Yau. 
\newblock
Existence of K\"ahler-Einstein metrics on complete K\"ahler manifolds and their applications to algebraic geometry.
\newblock
Mathematical aspects of string theory (San Diego, Calif., 1986), 574-628,
Adv. Ser. Math. Phys., 1, World Sci. Publishing, Singapore, 1987.

\bibitem{Z}
K. Zheng. 
\newblock
Geodesics in the space of K\"ahler cone metrics II: Uniqueness of constant scalar curvature K\"ahler cone metrics.
\newblock
Comm. Pure Appl. Math. 72 (2019), no. 12, 2621-2701.

\end{thebibliography}
\end{document}